\documentclass[preprint,3p,times]{elsarticle}

\usepackage{amsthm}
\newcommand{\Th}{\mathcal{T}_h}

\journal{Journal of Compuatational and Applied Mathematics}
\date{\today}

\usepackage[colorlinks=true,allcolors=purple]{hyperref}
\usepackage{graphicx}
\usepackage{amsfonts}
\usepackage{verbatim}
\usepackage{url}
\usepackage{float}
\usepackage{subcaption}
\usepackage{tabu}
\usepackage{bm}
\newcommand{\bx}{\bm{x}}
\usepackage{multirow}
\usepackage{rotating}
\usepackage{mathtools}
\usepackage{arydshln}
\usepackage[capitalise]{cleveref}
\usepackage[textsize=scriptsize,textwidth=1.9cm]{todonotes}

\usepackage{adjustbox}

\usepackage{pifont}

\newcommand{\assign}{\mathrel{\mathop:}=}
\newcommand{\rd}{\mathrm{d}}

\usepackage{pgfplotstable,booktabs}
\pgfplotsset{compat=newest}
\usepackage{tikz}
\usetikzlibrary{arrows.meta,positioning,calc,fit,shapes.geometric,shapes.multipart}

\usepackage{color}

\definecolor{darkpastelgreen}{rgb}{0.01, 0.75, 0.24}
\definecolor{burgundy}{rgb}{0.5, 0.0, 0.13}
\definecolor{darkgreen}{rgb}{0,0.7,0}

\definecolor{myred}{RGB}{200,0,0} 

\newcommand{\Cstab}{C_{\textup{\textsf{stab}}}}
\newcommand{\Csd}{C_{\textup{\textsf{sd}}}}
\newcommand{\amin}{a_{\textup{\textsf{min}}}}
\newcommand{\amax}{a_{\textup{\textsf{max}}}}
\newcommand{\Omegamax}{\Omega_{\textup{\textsf{max}}}}

\newtheorem{theorem}{Theorem}[section]
\newtheorem{lemma}[theorem]{Lemma}
\newtheorem{proposition}[theorem]{Proposition}
\newtheorem{corollary}[theorem]{Corollary}
\theoremstyle{definition}
\newtheorem{assumption}[theorem]{Assumption}
\newtheorem{remark}[theorem]{Remark}
\newtheorem{definition}[theorem]{Definition}
\newtheorem{notation}[theorem]{Notation}

\crefname{assumption}{Assumption}{Assumptions}
\Crefname{assumption}{Assumption}{Assumptions}

\numberwithin{equation}{section}

\usetikzlibrary{positioning}

\begin{document}

\begin{frontmatter}



\title{Can Symmetric Positive Definite (SPD) coarse spaces perform well for \\ indefinite Helmholtz problems?}

\author[1]{Victorita Dolean}\ead{v.dolean.maini@tue.nl}
\affiliation[1]{organization = {Eindhoven University of Technology},
addressline={PO Box 513}, city={Eindhoven},
postcode={5600 MB},
country={The Netherlands}}
\author[2]{Mark Fry}\ead{mark.fry@strath.ac.uk}
\author[2]{Matthias Langer} \ead{m.langer@strath.ac.uk}
\affiliation[2]{organization = {University of Strathclyde},
addressline={26 Richmond Street}, city={Glasgow},
postcode={G1 1XH},
country={United Kingdom}}

\begin{abstract}
Wave propagation problems governed by the Helmholtz equation remain among the most challenging in scientific computing, 
due to their indefinite nature.  Domain decomposition methods with spectral coarse spaces have emerged as some of the
most effective preconditioners, yet their theoretical guarantees often lag behind practical performance. 
In this work, we introduce and analyse the $\Delta_k$-GenEO coarse space within the two-level additive Schwarz preconditioners 
for heterogeneous Helmholtz problems.  This is an adaptation of the $\Delta$-GenEO coarse space. 
Our results sharpen the $k$-explicit conditions for GMRES convergence, reducing the restrictions on the subdomain size and eigenvalue threshold. 
This narrows the long-standing gap between pessimistic theory and empirical evidence, 
and reveals why GenEO spaces based on SPD (symmetric positive definite) eigenvalue problems remain surprisingly effective 
despite their apparent limitations. 
Numerical experiments confirm the theory, demonstrating scalability, robustness to heterogeneity for low to moderate frequencies 
(while experiencing limitations in the high frequency cases), and significantly milder coarse-space growth than conservative estimates predict. 
\end{abstract}


\end{frontmatter}

\section{Introduction}
\label{CH:Intro}
This work focuses on the numerical solution of Helmholtz boundary value problems with highly variable coefficients, where the governing equation is indefinite and arises in many wave propagation applications. Specifically, we consider
\begin{equation}
\label{eq:problem}
	\begin{alignedat}{2}
		-\nabla \cdot (A \nabla u) - k^2 u &= f \qquad && \text{in } \Omega, \\
		u &= 0 \qquad && \text{on } \partial\Omega,
	\end{alignedat}
\end{equation}
where $\Omega \subseteq \mathbb{R}^d$ is a bounded polygonal (for $d=2$) or Lipschitz polyhedral (for $d=3$) domain, $A(x)$ is a symmetric and uniformly positive definite matrix-valued coefficient field with possibly high contrast, $f \in L^2(\Omega)$, and $k>0$ is the wavenumber. We assume that \eqref{eq:problem} admits a unique weak solution $u \in H^1_0(\Omega)$ for all $f \in L^2(\Omega)$.

\medskip
\noindent
\textbf{Domain decomposition and coarse spaces.} To solve~\eqref{eq:problem}, we employ overlapping domain decomposition methods. The domain $\Omega$ is partitioned into overlapping subdomains $\{\Omega_i\}_{i=1}^N$, and a one-level additive Schwarz preconditioner is built from local solves on each $\Omega_i$. Since this preconditioner is not scalable as $N$ increases, a global coarse correction is introduced, yielding a two-level additive Schwarz preconditioner. After finite element discretisation, the resulting system matrix is symmetric but indefinite, so we employ GMRES as iterative solver. Its convergence is governed by the field-of-values theory of Elman~\cite{Elman:1983:VIM}, which relates iteration counts to spectral and stability properties of the preconditioned operator. 

The construction of robust coarse spaces is therefore central. Classical choices (piecewise polynomial functions on a coarse mesh~\cite{Cai:1999:RAS,Toselli:2005:DDM}) fail to deliver robustness in the presence of strong heterogeneity or high frequency. Operator-dependent coarse spaces were introduced to overcome this limitation~\cite{Graham:2007:DDM}, and the most successful among them are spectral coarse spaces built from local generalised eigenvalue problems~\cite{Spillane:2014:ARC}. The GenEO method is a prominent example in this family. Although the coarse space information must ultimately be shared globally, its construction can be carried out locally using basis vectors obtained from solutions of local eigenvalue problems. Spectral coarse spaces such as the Dirichlet-to-Neumann (DtN) coarse space exemplify this idea in the context of the Helmholtz equation \cite{Conen:2014:ACS}. These local eigenproblems are defined on subdomain interfaces through a DtN operator, extending concepts originally introduced for elliptic equations \cite{Nataf:2011:ACS, Dolean:2012:ATL}. Because theoretical analyses for such problems remain challenging, an initial comparative study of two related spectral coarse spaces, the DtN method and a variant inspired by GenEO adapted to the Helmholtz setting, was performed in \cite{Bootland:2022:CDG} on simple benchmark configurations. In parallel, \cite{Bootland:2021:ACS} provided a synthesis of recent developments aimed at improving domain decomposition techniques for heterogeneous media and implemented them within a unified software environment, FreeFEM. This synthesis was recently extended to other methods in \cite{dolean:2025:RHO}.

The original GenEO framework was analysed in the coercive SPD case~\cite{Spillane:2014:ARC}, 
showing robustness independent of coefficient heterogeneity. 
Subsequent work extended the analysis to indefinite and non-self-adjoint PDEs of convection--diffusion--reaction type~\cite{Bootland:2022:OSM}, 
where robustness of GMRES was established under explicit conditions on the coarse mesh diameter $H$ and on the 
eigenvalue tolerance $\tau$ used in the local eigenproblems. 
Specifically, Theorem~4.1 of~\cite{Bootland:2022:OSM} showed that robust convergence is guaranteed provided
\begin{equation}\label{eq:bdLam0}
   H \lesssim k^{-2}, 
   \qquad 
   (1+\Cstab)^2 \, k^8 \lesssim \tau,
\end{equation}
where $\Cstab$ is the stability constant of the adjoint problem. While these results provided the first rigorous analysis for indefinite problems (by introducing the so-called $\Delta$-GenEO method), the bounds~\eqref{eq:bdLam0} are overly pessimistic and far from observed practical behaviour.

\medskip
\noindent
\textbf{The $\Delta_k$-GenEO method.} In this paper, we introduce a new method, which we call $\Delta_k$-GenEO and 
which is an adaptation of the original $\Delta$-GenEO by modifying the right-hand side in the generalised eigenvalue problem 
making it now dependent on $k$. 
This will allow us not only to revisit the Helmholtz case but also to  clarify the role and limitations of SPD-based coarse spaces 
when applied to indefinite problems. 
Our alternative analysis shows that robust GMRES convergence can already be ensured under significantly milder conditions than for $\Delta$-GenEO:
\begin{equation}
\label{eq:delta_results_improved}
   H \lesssim k^{-1}, 
   \qquad 
   (1+\Cstab)^2 \, k^2 \lesssim \tau.
\end{equation}
Thus, compared to~\eqref{eq:bdLam0}, the dependence of $H$ on $k$ improves from quadratic to linear, while the dependence of $\tau$ improves from order $k^8$ to quadratic. These bounds are still conservative, and they highlight the fundamental limitations of building Helmholtz coarse spaces from nearby SPD problems. At the same time, they narrow the gap between rigorous theory and numerical evidence, providing a clearer explanation of why GenEO-type methods remain surprisingly effective in practice despite their suboptimal theoretical foundation.

\medskip
\noindent
\textbf{Relation to later developments.} 
Recent studies have introduced refined spectral coarse spaces for the Helmholtz problem, notably the $H$-GenEO and \( H_k \)-GenEO methods~\cite{Bootland:2021:ACS, Bootland:2022:CDG}. 
These formulations define local eigenproblems directly in terms of the indefinite Helmholtz operator, rather than relying on a nearby SPD surrogate, and they demonstrate remarkable robustness with respect to \( k \). 
However, establishing rigorous theoretical foundations for such methods remains an open and challenging task. 
Within this landscape, the present \( \Delta_k \)-GenEO framework serves a complementary purpose: it preserves the conceptual and computational simplicity of the well-established GenEO methodology while tightening the connection between its SPD-based analysis and the behaviour observed in practice. 
The examples discussed here are intentionally restricted to simple configurations---in contrast with the large-scale studies of~\cite{Operto:2023:I3F, Tournier:2019:MTI, Tournier:2022:3FD}---with the aim of clearly exposing both the advantages and the limitations of the approach. 
Although belonging to a different family than the recently proposed harmonic coarse spaces (see \cite{Ma:2025:TLR} and related works), the present formulation provides a numerically transparent and theoretically grounded framework that can inform and guide the development of more advanced variants.

\medskip
\noindent
\textbf{Contributions.} The purpose of this work is twofold: first, to make precise the
limitations inherent in using SPD-based GenEO coarse spaces for indefinite Helmholtz problems,
and second, to quantify how far these limitations can be pushed before more advanced
indefinite-operator coarse spaces become necessary. In this spirit, the main contributions of the paper are:
\begin{itemize}
    \item sharper $k$-explicit conditions~\eqref{eq:delta_results_improved} for robust GMRES convergence with GenEO coarse spaces, 
    	improving the previous quadratic and octic dependences on $k$ to linear and quadratic dependences, respectively;
    \item a theoretical framework, based on refined stability and projection estimates, that exposes the sources of suboptimality 
    	when SPD operators are used in the construction of coarse spaces for Helmholtz;
    \item numerical experiments showing that, despite these theoretical limitations, the $\Delta_k$-GenEO method performs surprisingly well 
    	in practice, with iteration counts and coarse space dimensions significantly milder than what the conservative theory predicts. 
\end{itemize}
Taken together, these results both delineate the boundaries of applicability of SPD-based spectral coarse spaces and provide 
a more realistic benchmark for their performance.  They also highlight the need for, and complement the development of, 
new approaches based on indefinite operators (e.g.\ $H$-GenEO or $H_k$-GenEO), whose analysis remains largely open.
The remainder of the manuscript is organised as follows. 
Section~\ref{sec:framework} introduces the functional framework, weak formulation, and discretisation of the Helmholtz problem, together with the domain decomposition setting and the construction of $\Delta_k$-GenEO coarse spaces. 
Section~\ref{sec:main-result} states the main convergence theorem for GMRES with the two-level preconditioner and develops the stability and approximation tools needed for the proof. 
Section~\ref{sec:proof} provides the detailed proof of the main result, highlighting the new estimates that yield sharper $k$-dependent conditions. 
Section~\ref{sec:numerics} presents a suite of numerical experiments that validate the theoretical predictions, assess robustness with respect to frequency, scalability, and heterogeneity, and explore practical coarse-space selection strategies. 

\section{Functional framework and discretisation}
\label{sec:framework}
We now lay out the functional framework, assumptions, and discretisation strategies that will underpin our analysis. In particular, we describe the weak formulation of the Helmholtz problem, introduce the associated bilinear forms, and detail the finite element discretisation used in constructing and analysing the preconditioners. This section also recalls key regularity and stability assumptions and highlights analytical tools---such as the Friedrichs inequality---that play a central role in deriving later estimates.


The model problem is the Helmholtz equation with heterogeneous coefficients as defined in \eqref{eq:problem}. In what follows we will state the main assumptions and technical lemmas.

\begin{assumption}[Ellipticity and scaling]\label{ass:A}
The coefficient matrix $A:\Omega \to \mathbb{R}^{d\times d}$ satisfies:
\begin{enumerate}
\item[\textup{(i)}] (Uniform ellipticity)
	There exist constants $0<\amin\le\amax$ such that
	\[
		\amin |\xi|^2 \le A(x)\xi\cdot\xi \le \amax |\xi|^2
		\qquad \text{for all} \ x\in\Omega, \;\; \xi\in\mathbb{R}^d.
	\]
\item[\textup{(ii)}] (Scaling) 
	Without loss of generality, we normalise $\amin=1$ and assume that the domain diameter $D_\Omega$
	satisfies $D_\Omega \le 1$.  Otherwise, the domain is rescaled.
\end{enumerate}
\end{assumption}

The weak formulation of \eqref{eq:problem} is: find $u \in H^1_0(\Omega)$ such that
\begin{equation}
\label{eq:weak}
	b(u,v) = (f,v) \qquad \text{for all} \ v \in H^1_0(\Omega),
\end{equation}
where the bilinear form $b(\cdot,\cdot): H^1_0(\Omega)\times H^1_0(\Omega)\to\mathbb{R}$ is given by
\begin{equation}\label{eq:def_b}
	b(u,v) = \int_\Omega \bigl( A\nabla u\cdot\nabla v - k^2 uv \bigr)\,\rd x.
\end{equation}

\begin{notation}
\label{not:forms}
For any subdomain $\Omega' \subseteq \Omega$ we use  $(\cdot,\cdot)_{\Omega'}$ to denote the $L^2(\Omega')$ inner product and $\|\cdot\|_{\Omega'}$ to denote the corresponding norm.
We also introduce the bilinear forms:
\begin{equation}
	a_{\Omega'}(u,v) \assign \int_{\Omega'} A\nabla u \cdot \nabla v \,\rd x, \qquad
	b_{\Omega'}(u,v) \assign a_{\Omega'}(u,v) - k^2(u,v)_{\Omega'}
	\quad\text{for} \ u,v\in H^1(\Omega'),
    \label{eq:2.18}
\end{equation}
and define the semi-norm induced by $a$, $\|u\|_{a,\Omega'} \assign \sqrt{a_{\Omega'}(u,u)}$. 
An important role is played by the $k$-weighted inner product:
\begin{equation*}
    (u,v)_{1, k,\Omega'} \assign a_{\Omega'}(u,v) + k^2 (u,v)_{\Omega'},
\end{equation*}
and we denote the induced $k$-norm by $\|u\|_{1, k, \Omega'}$. {When $A=I$ and $k=1$ this reduces to the classical $\|\cdot\|_{H^1(\Omega')}$ norm.}
When $\Omega'= \Omega$ we abandon the subscript $\Omega$ from these notations. 
\end{notation}

The form $a(\cdot,\cdot)$ is symmetric semi-positive definite on $H^1(\Omega)$, while $b(\cdot,\cdot)$ is symmetric but indefinite. 
We now state the solvability assumption.

\begin{assumption}[Stability]\label{ass:stab}
For every $f\in L^2(\Omega)$, the weak problem~\eqref{eq:weak} admits a unique solution $u\in H^1_0(\Omega)$. 
Moreover, there exists a constant $\Cstab>0$ \textup{(}independent of $f$\textup{)} such that
\begin{equation}\label{eq:stab_bound}
	\|u\|_{1,k}  \le \Cstab \|f\| \qquad \text{for all} \ f\in L^2(\Omega).
\end{equation}
\end{assumption}

{\begin{remark}
It is well known (and follows from the Lax--Milgram lemma and the Rellich--Kondrachov theorem) that the spectrum of 
the operator $u\mapsto -\nabla\cdot(A\nabla u)$ consists only of a sequence of eigenvalues.
If $k^2$ is different from all these eigenvalues (i.e.\ $k$ avoids the resonance frequencies), then Assumption~\ref{ass:stab} is satisfied.
\end{remark}
}

Let ${\cal T}_h$ be a shape-regular triangulation of $\Omega$ into simplices of maximal diameter $h$, and let $V^h\subseteq H^1_0(\Omega)$ 
be a conforming finite element space.  The Galerkin discretisation of~\eqref{eq:weak} is: find $u_h\in V^h$ such that
\begin{equation}\label{eq:galerkin}
	b(u_h,v) = (f,v) \qquad \text{for all} \ v\in V^h.
\end{equation}
Let $\{\phi_i\}_{i=1}^n$ be a basis for $V^h$. Then \eqref{eq:galerkin} yields the linear system
\begin{equation}\label{eq:lin_sys}
	\mathbf{B}\mathbf{u} = \mathbf{f}, 
\end{equation}
with entries $\mathbf{B}_{ij} = b(\phi_j,\phi_i)$, $\mathbf{u}_i = (u_h,\phi_i)$ and $\mathbf{f}_i = (f,\phi_i)$. 
Analogously, the stiffness matrix $\mathbf{A}$ is defined by $\mathbf{A}_{ij} = a(\phi_j,\phi_i)$.

\begin{lemma}[Schatz--Wang~\cite{Schatz:1996:SNE}]
\label{lem:schatz_wang}
Suppose that Assumptions~\ref{ass:A} and \ref{ass:stab} hold. Then there exists \( h_0 > 0 \) such that, for all \( 0 < h < h_0 \), 
the discrete problem \eqref{eq:galerkin} admits a unique solution \( u_h \in V^h \). 
Moreover, {let $u$ be the unique solution of \eqref{eq:weak}; 
then} for any \( \varepsilon > 0 \), there exists \( h_1 = h_1(\varepsilon) > 0 \) such that, for all \( h < h_1 \),
\begin{equation}\label{eq: 3_13_a}
	\| u - u_h \| \le \varepsilon \| u - u_h \|_{H^1(\Omega)}
\end{equation}
and
\begin{equation}\label{eq: 2_13}
	\| u - u_h \|_{H^1(\Omega)} \le \varepsilon \| f \|. 
\end{equation}
\end{lemma}

{\begin{remark}[Well-posedness and the Schatz--Wang lemma] \label{rem:schatz_wang}
Assumption~\ref{ass:stab} ensures well-posedness and stability of the continuous Helmholtz problem away from resonant frequencies, 
in the sense of a resolvent estimate.
At the discrete level, the bilinear form is indefinite, and stability is not automatic. 
The Schatz--Wang lemma guarantees discrete well-posedness of the Galerkin problem for sufficiently small $h$, provided the continuous problem is stable. 
This result is standard for finite element discretisations of Helmholtz problems; see, e.g.\ \cite{Ihlenburg:1998:FEA, Melenk:2011:WEC}. 
It will be repeatedly used in our stability and approximation estimates.
Combining \eqref{eq: 3_13_a} and \eqref{eq: 2_13} we obtain 
\begin{equation}\label{schatz_wang_1_k}
	\|u-u_h\|_{1,k} \le \varepsilon \|f\|
\end{equation}
for sufficiently small $h$.
\end{remark}
}

\noindent
A key tool in our analysis is the Friedrichs inequality.

\begin{lemma}[Friedrichs inequality~{\cite[Theorem 13.19]{leoni:2017:FCSS}}]
\label{lem:friedrichs}
Let $\Omega'\subseteq\mathbb{R}^d$ be an open set lying between two parallel hyperplanes at distance $L$.  Then, for all $u\in H^1_0(\Omega')$,
\begin{equation}
	\|u\|_{\Omega'} \le \tfrac{L}{\sqrt{2}\,} \|\nabla u\|_{\Omega'}.
\end{equation}
\end{lemma}

\begin{remark}
Combined with Assumption~\ref{ass:A}, Friedrichs' inequality yields: for any $\Omega'\subseteq\Omega$ of diameter $H$,
\begin{equation}\label{eq:friedrichs_energy}
	\|u\|_{\Omega'} \le \tfrac{H}{\sqrt{2}\,} \|\nabla u\|_{\Omega'} 
	\le \tfrac{H}{\sqrt{2}\,} \|u\|_a
	\le \tfrac{H}{\sqrt{2}\,} \|u\|_{1, k,\Omega'} \qquad \text{for all} \ u\in H^1_0(\Omega').
\end{equation}

In particular, for $\Omega'=\Omega$ we have
\begin{equation}\label{eq:friedrichs_energy_Omega}
	\|u\|_a \ge \sqrt{2}\,\|u\| \ge \|u\| \qquad \text{for all} \ u\in H^1_0(\Omega)
\end{equation}
since $D_\Omega\le1$ by Assumption ~\ref{ass:A}.
\end{remark}

\subsection{Domain decomposition framework}

We now introduce the domain decomposition framework needed for constructing one- and two-level preconditioners. The core idea is to partition the global domain into overlapping subdomains and then build local and global correction operators that accelerate convergence by improving information transfer across the entire domain.

We first partition $\Omega$ into $N$ non-overlapping subdomains $\{\Omega_i'\}_{i=1}^N$, each resolved by the global finite element mesh ${\cal T}_h$. To introduce overlap, each $\Omega_i'$ is enlarged by one or more layers of mesh elements:

\begin{definition}[Overlapping subdomains]\label{def:subdomains}
Given a non-overlapping subdomain $\Omega_i'$, its overlapping extension $\Omega_i$ is defined by
\[
	\Omega_i \assign \operatorname{Int}\left( \bigcup_{\ell : \operatorname{supp}(\phi_\ell) \cap \Omega_i' \neq \emptyset} \operatorname{supp}(\phi_\ell) \right),
\]
where $\operatorname{Int}(\cdot)$ denotes the interior. Multiple layers of overlap are obtained by recursive application of this extension.
\end{definition}
For each overlapping $\Omega_i$, we introduce the local finite element spaces
\begin{equation}\label{eq:local_spaces}
	\widetilde{V}_i \assign \{ v|_{\Omega_i} : v\in V^h \} \subseteq H^1(\Omega_i),
	\qquad
	V_i \assign \{ v\in \widetilde{V}_i : v|_{\partial \Omega_i}=0 \} \subseteq H^1_0(\Omega_i).
\end{equation}
Moreover, we set $H_i \assign \operatorname{diam}(\Omega_i)$ and $H \assign \max_i H_i$. 

\begin{remark}
The form $b_{\Omega_i}$ is generally indefinite.  However, for sufficiently small $H$, coercivity is restored locally; 
see Lemma~\ref{lemma_3_4}.
\end{remark}

For each $i=1,\dots,N$, let $E_i: V_i \to V^h$ be the zero-extension operator, and let $R_i: V^h \to V_i$ be its $L^2$-adjoint. These operators preserve bilinear forms:
\[
	a_{\Omega_i}(u,v) = a(E_i u, E_i v), \qquad
	b_{\Omega_i}(u,v) = b(E_i u, E_i v), \qquad
	(u,v)_{\Omega_i} = (E_i u, E_i v).
\]
In matrix form, the one-level additive Schwarz preconditioner reads
\begin{equation}\label{eq:AS1}
	\mathbf{M}^{-1}_{AS,1} = \sum_{i=1}^N \mathbf{E}_i \mathbf{B}_i^{-1} \mathbf{R}_i,
	\qquad \mathbf{B}_i = \mathbf{R}_i \mathbf{B} \mathbf{E}_i,
\end{equation}
where $\mathbf{E}_i,\mathbf{R}_i$ represent $E_i,R_i$ in the global finite element basis.

To improve scalability, a coarse space $V_0\subseteq V^h$ is added, with embedding $E_0:V_0\to V^h$ and adjoint $R_0$. The two-level preconditioner is then
\begin{equation}\label{eq:AS2}
	\mathbf{M}^{-1}_{AS,2} = \sum_{i=0}^N \mathbf{E}_i \mathbf{B}_i^{-1} \mathbf{R}_i,
	\qquad \mathbf{B}_i = \mathbf{R}_i \mathbf{B} \mathbf{E}_i.
\end{equation}
The preconditioned system takes the form
\begin{equation}\label{eq:precond_system}
	\mathbf{M}^{-1}_{AS,2}\mathbf{B}\mathbf{u} = \mathbf{M}^{-1}_{AS,2}\mathbf{f}.
\end{equation}
For analysis, we define local projection-like operators $T_i:V^h\to V_i$ by
\begin{equation}\label{eq:Ti}
	b_{\Omega_i}(T_i u,v) = b(u,E_i v) \qquad \text{for all} \ v\in V_i,
\end{equation}
where we set $\Omega_0\assign\Omega$,
and the global operator
\begin{equation}\label{eq:T}
	T \assign \sum_{i=0}^N E_i T_i.
\end{equation}

\begin{proposition}[Bootland et al.~\cite{Bootland:2022:OSM}]\label{prop:proj}
For any $u,v \in V^h$, with corresponding nodal vectors $\mathbf{u},\mathbf{v} \in \mathbb{R}^n$, 
\begin{equation} \label{eq: 2_27}
	\langle \mathbf{M}_{AS,2}^{-1} \mathbf{B} \mathbf{u}, \mathbf{v} \rangle_{\mathbf{D}_k} = (T u, v)_{1, k},
\end{equation}
where $\langle \cdot , \cdot \rangle_{\mathbf{D}_k}$ is the $\mathbf{D}_k$-inner product on $\mathbb{R}^n$ and the matrix $\mathbf{D}_k$ given by
\begin{equation}\label{def:Dk}
    \mathbf{D}_k \assign \mathbf{A} + k^2  \mathbf{S},
\end{equation}
where $\mathbf{S}_{ij}=(\phi_j,\phi_i)$;
the inner product and the corresponding norm are defined by
\begin{equation}\label{def:A_norm}
	\langle \mathbf{x},\mathbf{y}\rangle_{\mathbf{D}_k} \assign \mathbf{y}^{\!T}\mathbf{D}_k\mathbf{x}, 
	\qquad 
	\|\mathbf{x}\|_{\mathbf{D}_k} \assign \sqrt{\langle \mathbf{x},\mathbf{x}\rangle_{\mathbf{D}_k}}.
\end{equation}
\end{proposition}

\subsection{The {$\Delta_k$}-GenEO coarse space}

We now specify the coarse space $V_0$ in~\eqref{eq:AS2}. 
The $\Delta_k$-GenEO method, introduced here as a variant of the GenEO approach~\cite{Spillane:2014:ARC}, 
uses local spectral information from the overlap to capture global low-energy error modes.

\begin{definition}[Degrees of freedom, cf.\ \cite{Spillane:2014:ARC}]\label{def:dofs}
For a subdomain $\Omega_i$, define
\[
	\overline{\mathrm{dof}}(\Omega_i) \assign \{\ell: \operatorname{supp}(\phi_\ell)\cap \Omega_i \ne \emptyset\},
	\qquad
	\mathrm{dof}(\Omega_i) \assign \{\ell: \operatorname{supp}(\phi_\ell)\subseteq \overline{\Omega_i}\}.
\]
\end{definition}

\begin{definition}[Partition of unity]\label{def:partition}
For each global degree of freedom $j$, define its multiplicity
\[
	\mu_j = \#\{ i: j\in \mathrm{dof}(\Omega_i)\}.
\]
The local partition of unity operator $\Xi_i:\widetilde{V}_i\to V_i$ is given by
\[
	\Xi_i v = \sum_{j\in \mathrm{dof}(\Omega_i)} \frac{1}{\mu_j} v_j \phi_j^i 
	\qquad\text{where} \quad
	v = \sum_{j\in \overline{\mathrm{dof}}(\Omega_i)} v_j \phi_j^i
\]
{and $\phi_j^i\assign\phi_j|_{\Omega_i}$}.
\end{definition}

On each subdomain $\Omega_i$ we consider the following eigenvalue problem, which is used for 
the construction of the coarse space.

\begin{definition}[Local spectral problem]\label{def:DeltaGenEO}
For each $i\in\{1,\ldots,N\}$ find either $(p,\lambda)\in(\widetilde{V}_i\setminus\{0\})\times\mathbb{R}$
such that
\[
	a_{\Omega_i}(p,v) = \lambda\bigl(\Xi_i p,\Xi_i v\bigr)_{1,k,\Omega_i}
	\qquad\text{for all} \ v\in\widetilde{V}_i
\]
or $(p,\lambda)\in(\ker\Xi_i)\times\{\infty\}$.
We call $\lambda$ and $p$ eigenvalues and eigenfunctions respectively.
\end{definition}

The next lemma shows that the eigenvalue problem in Definition~\ref{def:DeltaGenEO} has sufficiently many eigenvectors.

\begin{lemma}\label{lem:eigenvectors}
Let $i\in\{1,\ldots,N\}$ and set $s_i\assign\dim V_i$ and $n_i\assign\widetilde{V}_i$.
There exist eigenvalues
\begin{equation}\label{order_eigenvalues}
	\lambda_1^i \le \cdots \le \lambda_{s_i}^i < \lambda_{s_i+1}^i = \cdots \lambda_{n_i}^i = \infty
\end{equation}
and corresponding eigenvectors $p_1^i,\ldots,p_{n_i}^i$ 
such that $\{p_1^i,\ldots,p_{n_i}^i\}$ is a basis of $\widetilde{V}_i$.
The eigenvectors $p_1^i,\ldots,p_{s_i}^i$ can be normalised such that
\begin{equation}\label{eigvec_norm}
	\bigl(\Xi_i p_\ell^i,\Xi_i p_m^i\bigr)_{{1,k,\Omega_i}} = \delta_{\ell m}, \qquad \ell,m\in\{1,\ldots,s_i\}.
\end{equation}
\end{lemma}

\begin{proof}
Set $c_i(v,w)\assign(\Xi_i v,\Xi_i w)_{1,k,\Omega_i}$, $v,w\in\widetilde{V}_i$.
We show that
\begin{equation}\label{trivial_kernel_inters}
	\ker a_{\Omega_i} \cap \ker c_i = \{0\},
\end{equation}
where, e.g.\ $\ker a_{\Omega_i}\assign\{v\in\widetilde{V}_i:a_{\Omega_i}(v,w)=0 \ \text{for all} \ w\in\widetilde{V}_i\}$.
Any function in $\ker a_{\Omega_i}$ must be constant.
Since $(\cdot,\cdot)_{1,k,\Omega_i}$ is positive definite, a function in $\ker c_i$ must be in $\ker\Xi_i$.
It follows from the definition of $\Xi_i$ that
\begin{equation}\label{kerXi}
	\ker\Xi_i = \operatorname{span}\{\phi_j^i:j\in\overline{\mathrm{dof}}(\Omega_i)\setminus\mathrm{dof}(\Omega_i)\},
\end{equation}
and hence a function is in $\ker\Xi_i$ if and only if it vanishes on the interior nodes.
As a consequence, a function in $\ker a_{\Omega_i}\cap\ker c_i$ must be the zero function.
This shows that \eqref{trivial_kernel_inters} holds.
It also follows from \eqref{kerXi} that $\dim\ker\Xi_i=n_i-s_i$.
The assertions of the lemma now follow from \cite[Lemma~3.14]{Bastian:2023:MSDD}.
\end{proof}

\begin{remark}
The local spectral problem defined in \cite{Bootland:2022:OSM}, known as the $\Delta$-GenEO method, is defined as,
\[
	a_{\Omega_i}(p,v) = \lambda a_{\Omega_i}\bigl(\Xi_i p,\Xi_i v\bigr)
	\qquad\text{for all} v\in\widetilde{V}_i
\]
The $\Delta_k$-GenEO method that we are introducing differs from this method by utilising the form $(\cdot,\cdot)_{1,k,\Omega_j}$ 
on the right-hand side of the problem instead of $a_{\Omega_j}(\cdot,\cdot)$.
\end{remark}

We now introduce the $\Delta_k$-GenEO coarse space.

\begin{definition}[$\Delta_k$-GenEO coarse space]\label{def:DeltaCoarse}
Let the notation be as in Lemma~\ref{lem:eigenvectors} and 
let $(p_\ell^i,\lambda_\ell^i)$, $\ell=1,\ldots,n_i$, be the eigenpairs of the
eigenvalue problem in Definition~\ref{def:DeltaGenEO} satisfying \eqref{order_eigenvalues} and \eqref{eigvec_norm},
and let $m_i\in\{1,\ldots,s_i\}$.
The global coarse space is defined as
\[
	V_0 \assign \operatorname{span}\bigl\{E_i \Xi_i p_\ell^i: \ell=1,\dots,m_i,\; i=1,\dots,N\bigr\}.
\]
Two quantities will play a central role:
\begin{equation}\label{def_tau}
	\Lambda \assign \max_{T\in \Th} \#\{ \Omega_i : T\subseteq \Omega_i \}, 
	\qquad
	\tau \assign \min_{1\le i \le N} \lambda^i_{m_i+1}.
\end{equation}
Here, $\Lambda$ measures maximal overlap multiplicity, while $\tau$ is the first unused eigenvalue, 
i.e.\ the spectral threshold defining the coarse space.
\end{definition}

\section{Statement of the main result and theoretical tools}
\label{sec:main-result}

Having introduced the two-level additive Schwarz preconditioner and the construction of the $\Delta_k$-GenEO coarse space, we now state our main GMRES convergence result for the preconditioned Helmholtz system. Throughout this section we write
\begin{equation}\label{def_Theta_delta_k}
	\Theta \assign \frac{1}{\tau},
\end{equation}
where $\tau$ is the spectral threshold defining the coarse space (see \eqref{def_tau}). 
We apply GMRES with the $\mathbf{D}_k$-inner product, defined in \eqref{def:A_norm}, and use Elman's field-of-values framework.

\begin{theorem}[GMRES convergence for the two-level preconditioned system]
\label{thm:convergence}
Assume that Assumptions~\ref{ass:A} and~\ref{ass:stab} hold and that one has discrete well-posedness with $h < h_1$ from Lemma~\ref{lem:schatz_wang}. 
Let $H$ be the maximal subdomain diameter, $\Lambda$ the overlap multiplicity, and $\tau$ the coarse-space threshold.  Suppose that
\begin{equation}
\label{eq:st_constraints}
\begin{aligned}
	s &\assign 8\Lambda\bigl(2+3\Lambda^2\Theta\bigr)(1+\Cstab)\,k\Theta^{1/2} < 1, 
	\\[1ex]
	t &\assign 6\sqrt{2}\,\Lambda\bigl(2+3\Lambda^2\Theta\bigr) H k < 1.
\end{aligned}
\end{equation}
Then, for GMRES applied in the $\langle\cdot,\cdot\rangle_{\mathbf{D}_k}$ inner product to
\(
	\mathbf{M}^{-1}_{\!AS,2}\mathbf{B}\mathbf{u}=\mathbf{M}^{-1}_{\!AS,2}\mathbf{f},
\)
the residuals satisfy
\begin{equation}
\label{eq:gmres_rate}
	\|\mathbf{r}^{(m)}\|_{\mathbf{D}_k}^2 
	\le \biggl(1-\frac{c_1^2}{c_2^2}\biggr)^{\!m}\;\|\mathbf{r}^{(0)}\|_{\mathbf{D}_k}^2,
\end{equation}
with
\begin{equation}\label{eq:c_1-c_2}
	c_1 \assign \frac{1-\max\{s,t\}}{(2+3\Lambda^2\Theta)},
	\qquad
	c_2 \assign 18+8\Lambda^2.
\end{equation}
\end{theorem}

\begin{remark}[On the proof ingredients]
The proof follows Elman's FoV estimate, combining: (i) a lower bound on $(Tu,u)_{1, k, \Omega_i}$ in terms of $\|u\|_{1, k, \Omega_i}^2$ 
(via stable decomposition and coarse-space approximation with constants depending on $\Lambda$ and $\Theta$), 
and (ii) an upper bound on $\|Tu\|_{1, k, \Omega_i}$ (via stability of the local projectors and the coarse projector). 
The refined local/coarse stability and approximation estimates---proved in the next subsections---produce 
the improved $k$-dependences reflected in~\eqref{eq:st_constraints}.
\end{remark}

\begin{corollary}[Simplified $k$-explicit conditions]
\label{cor:robust}
Let $k_0>0$ and assume that $k\ge k_0$.  If the relations in \eqref{eq:st_constraints} hold, then
\begin{equation}
	\label{eq:robust_simple}
	H \lesssim k^{-1},
	\qquad
	(1+\Cstab)^2\,k^2 \lesssim \tau.
\end{equation}
Conversely, if \eqref{eq:robust_simple} holds with small enough implicit constants in $\lesssim$, 
then \eqref{eq:st_constraints} is satisfied and the convergence rate constants in~\eqref{eq:c_1-c_2} are independent of $k$, 
coefficient heterogeneity, and $h$ (for $h < h_1$); so GMRES is mesh-independent and wavenumber-robust.

\end{corollary}

\begin{proof}
First assume that \eqref{eq:st_constraints} holds.
Clearly, $2+3\Lambda^2\Theta\ge 1$ and $\Lambda\ge1$, and hence
\[
	(1+\Cstab)k\tau^{-1/2} = (1+\Cstab)k\Theta^{1/2}
	\le (2+3\Lambda^2\Theta)\Lambda(1+\Cstab)k\Theta^{1/2} < \tfrac{1}{8},
\]
which yields the second inequality in \eqref{eq:robust_simple}.
Moreover, $Hk \le \bigl(2+3\Lambda^2\Theta\bigr)\Lambda H k< \frac{1}{6\sqrt{2}}$, which implies $H\lesssim k^{-1}$.

Now assume that
\[
	H \le C_1k^{-1}, \qquad (1+\Cstab)^2\,k^2 \le C_2\tau
\]
with $C_1,C_2>0$.
The second relation, together with $k\ge k_0$, implies that $\Theta=\frac{1}{\tau}\le C_2k_0^{-2}$.
We therefore have
\begin{align}
	\label{eq:s_restr}
	s &\le 8\Lambda(2+3\Lambda^2C_2k_0^{-2}) C_2^{1/2},
	\\[1ex]
	\label{eq:t_restr}
	t &\le 6\sqrt{2}\Lambda(2+3\Lambda^2C_2k_0^{-2})\Lambda C_1.
\end{align}
If $C_1$ and $C_2$ are small enough, then the right-hand sides of \eqref{eq:s_restr} and \eqref{eq:t_restr}
are strictly less than $1$, which implies that \eqref{eq:st_constraints} holds.
\end{proof}

{\begin{remark}[Interpretation and scaling]
\label{rem:interpretation}
The conditions stated in \eqref{eq:st_constraints} are \textit{sufficient conditions} obtained from a field-of-values analysis and are therefore conservative. 
Their purpose is to guarantee $k$-robust convergence of GMRES in a worst-case theoretical setting, uniformly with respect to coefficient heterogeneity and mesh size. This behaviour reflects the gap between sufficient theoretical conditions and observed performance, and is consistent with previous experience for spectral coarse spaces for indefinite problems. That said, from a practical perspective, the condition $H \lesssim k^{-1}$ corresponds to resolving each subdomain by a fixed number of points per wavelength, which is standard for Helmholtz problems (even if this condition can be relaxed in practice). 
Similarly, the lower bound on the spectral threshold $\tau$ ensures that coarse modes associated with low-energy components in the $k$-weighted norm are captured. 
\end{remark}
\paragraph{Choice of spectral threshold $\tau$}
In practice, however, $\tau$ should be interpreted as a tuning parameter that controls the trade-off between coarse-space size and robustness. 
The local generalised eigenvalue problems typically exhibit favourable spectral separation, with a small number of low-energy modes capturing the dominant global components. 
As a result, robust convergence is often achieved for values of $\tau$ that are significantly smaller than those suggested by the sufficient theoretical bounds.
}

\medskip
\noindent\textbf{Roadmap of the estimates.}
The remainder of this section develops the stability and approximation tools used in the proof of Theorem~\ref{thm:convergence}:
\begin{enumerate}
\item 
	\emph{Coarse-space approximation and stable decomposition} 
	(Lemmas~\ref{lem:schatz_wang}, \ref{lem:friedrichs} and the $\Delta_k$-GenEO projection bounds): 
	yields the lower FoV distance via $(Tu,u)_{1, k, \Omega_i}\gtrsim \|u\|_{1, k, \Omega_i}^2/(1+\Lambda^2\Theta)$.
\item 
	\emph{Local and coarse projector stability}: 
	provides $\|Tu\|_{1, k, \Omega_i} \lesssim (1+\Lambda^2)\|u\|_{1, k, \Omega_i}$ with the $k$-dependence isolated in $s$ and $t$.
\item
	\emph{FoV--GMRES contraction} (Elman): 
	combining (1)--(2) gives~\eqref{eq:gmres_rate} with $c_1,c_2$ as in~\eqref{eq:c_1-c_2}.
\end{enumerate}

\subsection{Coarse space approximation and stable decomposition}

We analyse the accuracy of the local spectral projection that maps $v\in \widetilde V_i$ 
onto the span of the $m_i$ dominant eigenfunctions produced by the $\Delta_k$-GenEO construction.

\begin{lemma}[Local spectral projection: stability and spectral tail bound]
\label{lem:local-proj}
Let the notation be as in Lemma~\ref{lem:eigenvectors} and fix $i\in\{1,\dots,N\}$.  
Let $(p_\ell^i,\lambda_\ell^i)$, $\ell=1,\ldots,n_i$, 
be the eigenpairs of Definition~\ref{def:DeltaGenEO} satisfying \eqref{order_eigenvalues} and \eqref{eigvec_norm}.
Let $m_i\in\{1.\ldots,s_i-1\}$ and define the projector $\Pi_{m_i}^i:\widetilde{V}_i \to\widetilde{V}_i$ by

\begin{equation}\label{eq:Pi-def}
	\Pi_{m_i}^i v \assign \sum_{\ell=1}^{m_i} \bigl(\Xi_i v,\Xi_i p_\ell^i\bigr)_{1,k,\Omega_i}\,p_\ell^i.
\end{equation}
Then, for $w\assign v-\Pi_{m_i}^i v$, we have
\begin{equation}\label{eq:local-bounds}
	\|w\|_{{a,\Omega_i}}^2 \le \|v\|_{1,k,\Omega_i}^2,
	\qquad
	\|\Xi_i w\|_{1,k,\Omega_i}^2 \le \frac{1}{\lambda_{m_i+1}^i}\,\|w\|_{{a,\Omega_i}}^2.
\end{equation}
Moreover, $\Pi_{m_i}^i$ is the $(\Xi_i\cdot,\Xi_i\cdot)_{1, k, \Omega_i}$-orthogonal projector onto $\operatorname{span}\{p_1^i,\dots,p_{m_i}^i\}$, 
hence minimises $\|\Xi_i(v-z)\|_{1, k, \Omega_i}$ over all $z$ in that span.
\end{lemma}

\begin{proof}
The proof of Lemma~\ref{lem:eigenvectors} shows that we can apply \cite[Lemma~3.15]{Bastian:2023:MSDD},
which yields
\[
	\|w\|_{{a,\Omega_i}}^2 \le \|v\|_{{a,\Omega_i}}^2 \le \|v\|_{1,k,\Omega_i}^2
\]
and the second relation in \eqref{eq:local-bounds}.
\end{proof}

We often use the following standard overlap stability relations:
\begin{equation}\label{overlap_stab}
	\bigg\| \sum_i E_i q_i \bigg\|_{1,k}^2 \le \Lambda \sum_i \|q_i\|_{1,k,\Omega_i}^2 \qquad\text{and}\qquad
	\sum_i \|v|_{\Omega_i}\|_{1,k,\Omega_i}^2 \le \Lambda \|v\|_{1,k}^2.
\end{equation}

\begin{lemma}[Global approximation by the $\Delta_k$-GenEO coarse space]
\label{lem:global-approx}
For $v\in V^h$ define
\[
	z_0 \assign \sum_{i=1}^N E_i\,\Xi_i \,\!\bigl(\Pi_{m_i}^i(v|_{\Omega_i})\bigr)\;\in V_0.
\]
Then
\begin{equation}\label{eq:global-approx}
	\inf_{z\in V_0}\|v-z\|_{1, k}^2 \;\le\; \|v-z_0\|_{1, k}^2
	\le \Lambda^2\Theta\|v\|_{1, k}^2,
\end{equation}
where $\Theta$ is as in \eqref{def_Theta_delta_k}.
\end{lemma}

\begin{proof}
Set $w_i\assign v|_{\Omega_i} - \Pi_{m_i}^i(v|_{\Omega_i})$ so that 
\[
	v-z_0 = \sum_{i=1}^N E_i\,\Xi_i w_i.
\]
It follows from Lemma~\ref{lem:local-proj} that
\[
	\|\Xi_i w_i\|_{1,k,\Omega_i}^2 \le \frac{1}{\lambda_{m_i+1}^i}\|w_i\|_{{a,\Omega_i}}^2
	\le \Theta\|w_i\|_{{a,\Omega_i}}^2
	\le \Theta\|v|_{\Omega_i}\|_{1,k,\Omega_i}^2;
\]
note that
\[
	\Theta \assign \max_{1\le i\le N}\frac{1}{\lambda_{m_i+1}^i}.
\]
Summing over $i$ and using the second overlap stability bound in \eqref{overlap_stab} we obtain
\begin{equation}\label{sum_norm_zi}
	\sum_{i=1}^N \|\Xi_i w_i\|_{1,k,\Omega_i}^2 
	\le \Theta \sum_{i=1}^N \|v|_{\Omega_i}\|_{1,k,\Omega_i}^2
	\le \Lambda\Theta \|v\|_{1,k}^2.
\end{equation}
Using \eqref{overlap_stab} again we arrive at
\[
	\|v-z_0\|_{1,k}^2 = \bigg\|\sum_{i=1}^N E_i\Xi_i w_i\bigg\|_{1,k}^2
	\le \Lambda \sum_{i=1}^N \|\Xi_i w_i\|_{1,k,\Omega_i}^2
	\le \Lambda^2\Theta \|v\|_{1,k}^2,
\]
which proves \eqref{eq:global-approx}.
\end{proof}

\begin{lemma}[Stable decomposition]
\label{lem:stable-decomp}
Let $z_0$ be as in Lemma~\ref{lem:global-approx} and define, for $i=1,\dots,N$,
\[
	z_i \assign \Xi_i\bigl(v|_{\Omega_i} - \Pi_{m_i}^i(v|_{\Omega_i})\bigr) \;\in V_i.
\]
Then
\begin{equation}\label{stable_decom_ml}
	v = \sum_{i=0}^N E_i z_i,
	\qquad
	\|z_0\|_{1,k}^2 + \sum_{i=1}^N \|z_i\|_{1, k, \Omega_i}^2
	\le \bigl(2+3\Lambda^2\Theta\bigr)\,\|v\|_{1, k}^2.
\end{equation}
\end{lemma}

\begin{proof}
Let $w_i$ be as in the proof of Lemma~\ref{lem:global-approx}; then $z_i=\Xi_i w_i$.
The first relation in \eqref{stable_decom_ml} follows easily from the definition of $z_i$.
Using $(x+y)^2\le 2x^2+2y^2$ and Lemma~\ref{lem:global-approx} we obtain
\[
	\|z_0\|_{1,k}^2 \le \bigl(\|v\|_{1,k}+\|v-z_0\|_{1,k}\bigr)^2
	\le 2\|v\|_{1,k}^2 + 2\|v-z_0\|_{1,k}^2
	\le \bigl(2+2\Lambda^2\Theta\bigr)\|v\|_{1,k}^2.
\]
Adding $\sum_{i=1}^N \|z_i\|_{1,k,\Omega_i}^2$ on both sides and using \eqref{sum_norm_zi} and $\Lambda\ge1$
we arrive at the second relation in \eqref{stable_decom_ml}.
\end{proof}

\begin{remark}[What these bounds tell us]
Lemma~\ref{lem:global-approx} provides the coarse approximation estimate; Lemma~\ref{lem:stable-decomp} yields the stable splitting. 
Together, they furnish the lower FoV distance and the operator-norm bound that feed Elman's GMRES contraction, giving Theorem~\ref{thm:convergence}.
\end{remark}

\begin{remark}[Classical GenEO estimates]
The local spectral projection bounds (Lemma~\ref{lem:local-proj}), the global approximation
estimate (Lemma~\ref{lem:global-approx}), and the stable decomposition
(Lemma~\ref{lem:stable-decomp}) are classical within the GenEO framework; see
\cite{Spillane:2014:ARC} and the indefinite extension in \cite{Bootland:2022:OSM}.
We re-state them here for completeness and to fix notation in the $\Delta_k$-GenEO setting.
\end{remark}

\subsection{Local projector: solvability, basic identities, and stability}
To facilitate the analysis, we introduce, for each $i=1,\dots,N$, the local \emph{projector map}
$P_i:V^h\to V_i$ defined by
\begin{equation}\label{eq:proj_ritz_local}
	(P_i u, v)_{1, k, \Omega_i} = (u,E_i v)_{1, k}, \qquad v\in V_i,
\end{equation}
and the coarse projector $P_0:V^h\to V_0$ defined by
\begin{equation}\label{eq:proj_ritz_coarse}
	(P_0 u, v_0)_{1, k} = (u,v_0)_{1, k}, \qquad v_0\in V_0 .
\end{equation}
The associated global operator is
\[
	P \assign \sum_{i=0}^N E_i P_i :\; V^h \to V^h .
\]

\begin{lemma}[Well-posedness and algebraic identities]\label{lem:ritz-properties}
Under Assumptions~\ref{ass:A} and \ref{ass:stab}, the bilinear forms $a_{\Omega_i}$ and $a$ are SPD on $V_i$ and $V_0$, respectively; 
hence $P_i$ \textup{(}$i\ge 1$\textup{)} and $P_0$ are well defined by Lax--Milgram.  Moreover, the following statements hold.
\begin{enumerate}
\item (Locality) 
	For $v\in V_i$, $(u,E_i v)_{1, k}=(u|_{\Omega_i},v)_{1, k, \Omega_i}$; so $P_iu$ is a projection of $u|_{\Omega_i}$ onto $V_i$.
\item (Orthogonal projector) 
	$P_0$ is the $(\cdot, \cdot)_{1,k}$-orthogonal projector onto $V_0$; in particular, $\|P_0 u\|_{1, k} \le \|u\|_{1, k}$.
\item (Identity) 
	For all $u\in V^h$,
	\begin{equation}\label{eq:energy-identity}
		(Pu,u)_{1, k} = \|P_0 u\|_{1, k}^2 + \sum_{i=1}^N \|P_i u\|_{1, k, \Omega_i}^2 .
	\end{equation}
\end{enumerate}
\end{lemma}

\begin{proof}[Proof (sketch)]
SPD and Lax--Milgram are standard.  
The first statement follows from the zero-extension property of $E_i$ and the definition of $(\cdot , \cdot)_{1, k, \Omega_i}$.
The second claim is immediate from \eqref{eq:proj_ritz_coarse} and standard properties of orthogonal projectors in Hilbert spaces.
For the third statement, use \eqref{eq:proj_ritz_local} and \eqref{eq:proj_ritz_coarse} with $v=P_i u$ or $v_0=P_0 u$ and sum over $i=0,\dots,N$.
\end{proof}

The next result links the global weighted $k$-norm of $u$ to the action of $P$ and controls the size of the local images.

\begin{proposition}[Inequalities for the local operators]\label{prop:ineqp}
Let $u\in V^h$. With the stable decomposition of Lemma~\ref{lem:stable-decomp} and its constant
$\Csd \assign 2+3\Lambda^2\Theta$, we have
\begin{equation}\label{eq:ineq-global}
	\|u\|_{{1,k}}^2 \le \Csd\, (Pu,u)_{1, k}
\end{equation}
and
\begin{equation}\label{eq:ineq-local-sum}
	\sum_{i=0}^N \|P_i u\|_{1, k, \Omega_i}^2 \le (\Lambda+1)\|u\|_{1, k}^2.
\end{equation}
\end{proposition}

\begin{proof}
Let $u=\sum_{i=0}^N E_i z_i$ be the decomposition from Lemma~\ref{lem:stable-decomp}. 
Then, by \eqref{eq:proj_ritz_local} and \eqref{eq:proj_ritz_coarse},
\[
	\|u\|_{1, k}^2 = \biggl(u, \sum_{i=0}^N E_i z_i\biggr)_{1, k}
	= \sum_{i=0}^N (u,E_i z_i)_{1, k}
	= \sum_{i=0}^N (P_i u, z_i)_{1, k, \Omega_i}.
\]
Cauchy--Schwarz (applied twice) and Lemma~\ref{lem:stable-decomp} give
\[
	\|u\|_{1,k}^2 \le \sum_{i=0}^N \|P_i u\|_{1,k,\Omega_i}\|z_i\|_{1,k,\Omega_i}
	\le \biggl(\sum_{i=0}^N \|P_i u\|_{1,k,\Omega_i}^2\biggr)^{\!\!1/2}
	\biggl(\|z_0\|_{1,k}^2+\sum_{i=1}^N \|z_i\|_{1,k,\Omega_i}^2\biggr)^{\!\!1/2}
	\le \Csd^{1/2}\; \biggl(\sum_{i=0}^N \|P_i u\|_{1,k,\Omega_i}^2\biggr)^{\!\!1/2}\;\|u\|_{1,k}.
\]
Dividing by $\|u\|_{1,k}$, using the identity \eqref{eq:energy-identity} to replace the sum by $(Pu,u)_{1,k,\Omega_i}$ 
and squaring the inequality we obtain \eqref{eq:ineq-global}.

For \eqref{eq:ineq-local-sum}, we first bound the coarse part by contraction of $P_0$:
$\|P_0 u\|_{1,k}^2 \le \|u\|_{1,k}^2$.
For the overlapping parts we use \eqref{overlap_stab} to obtain
\begin{align*}
	\sum_{i=1}^N \|P_i u\|_{1,k,\Omega_i}^2
	&= \sum_{i=1}^N (u,E_i P_i u)_{1,k}
	= \left(u, \sum_{i=1}^N E_i P_i u\right)_{1,k}
	\le \|u\|_{1,k} \,\Big\|\sum_{i=1}^N E_i P_i u\Big\|_{1,k}
	\\[1ex]
	&\le \Lambda^{1/2}\|u\|_{1,k} \Biggl(\sum_{i=1}^N \|P_i u\|_{1,k,\Omega_i}^2\Biggr)^{1/2},
\end{align*}
which implies 
\begin{equation}\label{sum_norm_P_i}
	\sum_{i=1}^N \|P_i u\|_{1,k,\Omega_i}^2 \le \Lambda \|u\|_{1,k}^2.
\end{equation} 
Adding the coarse term we arrive at \eqref{eq:ineq-local-sum}.
\end{proof}

\begin{remark}[Classical facts]
The definition and properties of $P_i$ and $P_0$ are standard in domain decomposition; 
see, e.g.\ \cite[Sec.~2.2]{Toselli:2005:DDM}. We include them here to keep the presentation self-contained and to fix notation.
\end{remark}

\medskip

We now return to the Helmholtz-specific projectors $T_i$ defined by the indefinite form $b_{\Omega_i}(\cdot,\cdot)$ in \eqref{eq:Ti}. 
The next two lemmas give solvability (local coercivity) and a uniform $(\cdot, \cdot)_{1, k, \Omega_i}$-stability bound.

\begin{lemma}[Solvability of $T_i$]\label{lemma_3_4}
If $Hk<\sqrt{2}$, then $b_{\Omega_i}$ is coercive on $V_i$ and the operators $T_i:V^h\to V_i$ defined by
\[
	b_{\Omega_i}(T_i u, v) = b(u,E_i v) \qquad \text{for all} \ v\in V_i
\]
are well posed for $i=1,\dots,N$.
\end{lemma}

\begin{proof}
By Friedrichs' inequality \eqref{eq:friedrichs_energy} we have
$\|u\|_{\Omega_i} \le \frac{H}{\sqrt{2}\,} \|u\|_{{a,\Omega_i}}$ and hence
\[
	b_{\Omega_i}(u,u) = a_{\Omega_i}(u,u) - k^2\|u\|_{\Omega_i}^2
	\ge \Big(1 - \tfrac{k^2 H^2}{2}\Big)\,\|u\|_{{a,\Omega_i}}^2.
\]
Thus $b_{\Omega_i}$ is coercive for $Hk<\sqrt{2}$, and Lax--Milgram applies.
\end{proof}

Before we proceed, we state the following boundedness property of  $b_{\Omega'}$, 
which is a simple consequence of the Cauchy--Schwarz inequality.

\begin{lemma}\label{prop:est_b}
\begin{equation}\label{estimate_b}
	|b_{\Omega'}(u,v)| \le \|u\|_{1,k,\Omega'}\|v\|_{1,k,\Omega'}. 
\end{equation}
\end{lemma} 

\begin{proof}
Using the definitions of $a_{\Omega'}$ and $b_{\Omega'}$ and the Cauchy--Schwarz inequality we obtain
\begin{align*}
	|b_{\Omega'}(u,v)| &= |a_{\Omega'}(u,v) - k^2 (u,v)_{\Omega'}| 
	\le \|u\|_{{a,\Omega'}} \|v\|_{{a,\Omega'}} + k^2 \|u\|_{\Omega'}\|v\|_{\Omega'} 
	\\[1ex] 
	&\le \sqrt{\|u\|_{{a,\Omega'}}^2 + k^2 \|u\|_{\Omega'}^2} \sqrt{\|u\|_{{a,\Omega'}}^2 + k^2 \|u\|_{\Omega'}^2} 
	= \|u\|_{1, k, \Omega'} \|v\|_{1, k, \Omega'},
\end{align*}
which proves \eqref{estimate_b}.
\end{proof}

\begin{lemma}[Stability of $T_i$]\label{lemma_3_6}
If $Hk\le 1/\sqrt{2}$, then, for all $u\in V^h$ and $i=1,\dots,N$,
\begin{equation}\label{stab_T_i}
	\|T_i u\|_{1,k,\Omega_i} \le  2\|u|_{\Omega_i}\|_{1,k,\Omega_i}.
\end{equation}
\end{lemma}

\begin{proof}
It follows from \eqref{estimate_b} and \eqref{eq:friedrichs_energy} that
\begin{align*}
	\|T_i u\|_{1,k,\Omega_i}^2
	&= b_{\Omega_i}(T_i u,T_i u) + 2k^2\|T_i u\|_{\Omega_i}^2
	= b(u,E_i T_i u) + 2k^2\|T_i u\|_{\Omega_i}^2
	= b_{\Omega_i}(u|_{\Omega_i},T_i u) + 2k^2\|T_i u\|_{\Omega_i}^2
	\\[1ex]
	&\le \|u|_{\Omega_i}\|_{1,k,\Omega_i} \|T_i u\|_{1,k,\Omega_i} + k^2 H^2\|T_i u\|_{1,k,\Omega_i}^2.
\end{align*}
Rearranging this inequality we obtain
\[
	\bigl(1-k^2 H^2\bigr)\|T_i u\|_{1,k,\Omega_i}
	\le \|u|_{\Omega_i}\|_{1,k,\Omega_i},
\]
which, together with the assumption $Hk\le\frac{1}{\sqrt{2}\,}$, yields \eqref{stab_T_i}.
\end{proof}

\subsection{Coarse projector $T_0$: solvability and stability}
We define the coarse (global) Helmholtz projector $T_0:V^h\to V_0$ by
\begin{equation}\label{eq:def-T0}
	b(T_0 u, v_0) = b(u, v_0) \qquad \text{for all} \ v_0\in V_0.
\end{equation}
Thus $e_0\assign T_0 u-u$ satisfies $b(e_0,v_0)=0$ for all $v_0\in V_0$.

\begin{lemma}[Solvability of $T_0$]\label{lem:T0-solv}
Suppose that the spectral threshold $\Theta$ satisfies
\begin{equation}\label{eq:T0-infsup}
	\sqrt{2} k\Lambda\Theta^{1/2}(1+\Cstab) < 1.
\end{equation}
Then there exists $h_1>0$ such that, for all $h<h_1$, the problem \eqref{eq:def-T0} is well posed, 
i.e.\ $T_0$ is uniquely defined on $V^h$.
\end{lemma}

\begin{proof}
We argue by contradiction.  Assume that there exists $w_0\in V_0\setminus\{0\}$ with
$b(w_0,z)=0$ for all $z\in V_0$.  Let $w\in H^1_0(\Omega)$ solve $b(w,v)=(w_0,v)$
for all $v\in H^1_0(\Omega)$ (such a $w$ exists by Assumption~\ref{ass:stab}) and let $w_h\in V^h$ solve $b(w_h,v)=(w_0,v)$ for all $v\in V^h$.
Taking $v=w_0$ in the latter relation and using Lemma~\ref{prop:est_b} we obtain, for all $z\in V_0$,
\[
	\|w_0\|^2 = b(w_h,w_0) = b(w_h-z,w_0)
	\le \|w_h-z\|_{1,k}\|w_0\|_{1,k}.
\]
Minimising over $z\in V_0$ and using the global $\Delta_k$-GenEO approximation (Lemma~\ref{lem:global-approx}) with $v=w_h$ we get
\begin{equation}\label{ml01}
	\|w_0\|^2 \le \Lambda\Theta^{1/2}\|w_0\|_{1, k}\|w_h\|_{1, k}.
\end{equation}
Pick $z=w_0$ in $b(w_0,z)=0$ to obtain $(w_0,w_0)_{1, k}=2k^2\|w_0\|^2$, which implies $\|w_0\|_{1, k}=\sqrt{2} k\|w_0\|$. 
Assumption~\ref{ass:stab} and Remark~\ref{rem:schatz_wang} yield
\begin{equation}\label{appl_schatz_wang}
	\|w_h\|_{1,k} \le \|w\|_{1,k} + \|w-w_h\|_{1,k} 
	\le \Cstab\|w_0\| + \|w-w_h\|_{1,k} 
	\le (1+\Cstab)\|w_0\|.
\end{equation}
Altogether, we have
\[
	\|w_0\|^2 \le \sqrt{2} k\Lambda\Theta^{1/2}(1+\Cstab)\|w_0\|^2,
\]
which contradicts \eqref{eq:T0-infsup} since $w_0 \ne 0$.  Hence $T_0$ is well defined.
\end{proof}

\begin{lemma}[Stability of $T_0$]\label{lem:T0-stab}
Assume that \eqref{eq:T0-infsup} holds and that $h < h_1$ as above. Then, for all $u\in V^h$,
\begin{equation}\label{eq:T0-L2-energy}
	\|T_0 u - u\| \le \,\Lambda\Theta^{1/2}(1+\Cstab)\|T_0 u - u\|_{1, k} .
\end{equation}
Moreover, if
\begin{equation}\label{eq:T0-extra}
	2k\Lambda\Theta^{1/2}(1+\Cstab) \le \tfrac12,
\end{equation}
then
\begin{equation}\label{eq:T0-energy}
	\|T_0 u - u\|_{1, k} \le 2\|u\|_{1, k}.
\end{equation}
\end{lemma}

\begin{proof}
Let $e_0\assign u-T_0 u$.  Solve $b(w_h,v)=(e_0,v)$ for $w_h\in V^h$.
Since $b(e_0,z)=0$ for $z\in V_0$ by the definition of $T_0$, we have
\[
	\|e_0\|^2 = b(w_h,e_0) = b(w_h-z,e_0)\qquad\text{for all} \ z\in V_0.
\]
We proceed exactly as in the proof of Lemma~\ref{lem:T0-solv} (with $w_0$ replaced by $e_0$, cf.\ \eqref{ml01}) to obtain
\[
	\|e_0\|^2 \le \Lambda\Theta^{1/2}\|w_h\|_{1, k}^{\phantom{1}}\,\|e_0\|_{1, k}.
\]
Using \eqref{appl_schatz_wang} (again with $w_0$ replaced by $e_0$) we find that
\[
	\|e_0\| \le \Lambda\Theta^{1/2}(1+\Cstab)\|e_0\|_{1,k},
\]
which is \eqref{eq:T0-L2-energy}. 

To prove \eqref{eq:T0-energy}, let $P_0$ be as in \eqref{eq:proj_ritz_coarse}.
Since $P_0u-T_0u\in V_0$, we have $b(e_0,P_0u-T_0u)=0$, which, together with the fact that $P_0$ 
is orthogonal with respect to $(\cdot, \cdot)_{1, k}$, yields
\begin{align*}
	\|e_0\|_{1, k}^2 &= b(e_0,e_0) + 2k^2(e_0,e_0)
	= b(e_0,u-T_0u) - b(e_0,P_0u-T_0u) + 2k^2(e_0,e_0)
	\\[1ex]
	&= b(e_0,u-P_0u) + 2k^2(e_0,u-T_0u)
	= (e_0,u-P_0u)_{1,k} - 2k^2(e_0,u-P_0u) + 2k^2(e_0,u-T_0u)
	\\[1ex]
	&= (e_0,u-P_0u)_{1, k} + 2k^2(e_0,P_0u-T_0u)
	\\[1ex]
	&\le \|e_0\|_{1, k} \|u-P_0u\|_{1, k} + 2k^2\|e_0\|\,\|P_0u-T_0u\|
	= \|e_0\|_{1, k} \|u-P_0u\|_{1, k} + 2 k^2 \|e_0\|\,\|P_0e_0\|
	\\[1ex]
	&\le \|e_0\|_{1, k} \|u-P_0u\|_{1, k} + 2 k\|e_0\|\,\|P_0e_0\|_{1, k}
	\le \|e_0\|_{1, k} \|u\|_{1, k} + 2 k\|e_0\|\,\|e_0\|_{1, k}.
\end{align*}
We can now use \eqref{eq:T0-L2-energy} and \eqref{eq:T0-extra} to obtain
\[
	\|e_0\|_{1, k} \le \|u\|_{1, k} + 2k\Lambda\Theta^{1/2}(1+\Cstab)\|e_0\|_{1, k}
	\le \|u\|_{1, k} + \tfrac12\|e_0\|_{1, k},
\]
which, in turn, implies \eqref{eq:T0-energy}.
\end{proof}

\begin{remark}[Comparison and provenance]
Conditions \eqref{eq:T0-infsup} and \eqref{eq:T0-extra} strengthen the coarse-level assumptions in~\cite{Bootland:2022:OSM}
by leveraging the sharper GenEO approximation (Lemma~\ref{lem:global-approx}). The proof strategy is classical:
a contradiction argument for solvability (coercivity of $b$ on $V_0$) and a Petrov--Galerkin stability estimate for \eqref{eq:T0-L2-energy}.
\end{remark}

\begin{corollary}[$\tau$-form of the coarse conditions]\label{cor:T0-tau}
Let $\Theta=1/\tau$ and $\Lambda\ge 1$. Then condition \eqref{eq:T0-extra} is equivalent to
\begin{equation}\label{eq:T0-solv-tau}
		\tau \ge 16\Lambda^2(1+\Cstab)^2k^2.
\end{equation}
\end{corollary}

\begin{remark}[Interpretation of the coarse-level condition]\label{rem:T0_condition}
{Conditions~\eqref{eq:T0-infsup} and \eqref{eq:T0-extra} arise from the requirement that the coarse-level Helmholtz projector $T_0$ 
be well defined and stable.  They ensure coercivity of the bilinear form $b(\cdot,\cdot)$ on the coarse space $V_0$, 
and are therefore necessary ingredients in the field-of-values analysis underpinning Theorem~\ref{thm:convergence}.}

{Note that the first inequality in \eqref{eq:st_constraints} implies \eqref{eq:T0-extra}, which, in turn, implies \eqref{eq:T0-infsup}.
All three conditions impose constraints on the spectral threshold $\tau$ through the quantity $\Theta = 1/\tau$
and ultimately stem from the same coarse-space approximation estimate provided by Lemma~\ref{lem:global-approx}. 
Corollary~\ref{cor:T0-tau} makes this connection explicit by showing that~\eqref{eq:T0-infsup} is equivalent to a quadratic lower bound on $\tau$, 
namely $\tau \gtrsim (1+C_{\mathrm{stab}})^2 k^2$, which is consistent with the simplified $k$-explicit conditions in~\eqref{eq:robust_simple}.
From a practical perspective, condition~\eqref{eq:T0-infsup} should again be interpreted as a sufficient condition. 
It guarantees that all coarse modes responsible for potential instabilities are captured, 
but it does not require that the bound on $\tau$ be saturated in practice. 
As observed in the numerical experiments, stable and efficient convergence is typically achieved for values of $\tau$ 
that are significantly smaller than the worst-case theoretical bound suggested by~\eqref{eq:T0-infsup}. 
}
\end{remark}

\paragraph{On the size of the coarse space and the $\tau$-bound}
The conditions in Corollary~\ref{cor:T0-tau} constrain the \emph{threshold} $\tau$, not the number of kept local modes $m_i$ directly. 
Translating a bound on $\tau$ into a bound on $m_i$ requires information about the spectral density
$\{\lambda_\ell^i\}_{\ell\ge1}$ of the local problems, which depends on geometry, overlap, $H$, and coefficients. 
In the absence of such structure, the sufficient bound $\tau \gtrsim (1+\Cstab)^2 k^2$
can be read as `include all local modes whose Rayleigh quotients are $\lesssim k^2$', 
which may \emph{appear} to force a rapidly growing coarse space as $k$ increases.
In practice, two points mitigate this pessimism:
\begin{itemize}
\item
	\textbf{Sufficiency vs.\ necessity.} 
	The $k^2$ scaling is a \emph{worst-case sufficient} condition for field-of-values robustness. 
	Empirically, with $H\sim c/k$ and moderate overlap, the spectrum of $(\Xi_i\cdot,\Xi_i\cdot)_{1, k, \Omega_i}$ 
	typically exhibits an early cluster/gap, so only a modest number of modes per subdomain are needed to stabilise GMRES. 
	Thus, $m_i$ often grows much more slowly than what a literal reading of the $\tau$-inequality would suggest.
\item 
	\textbf{What matters for GMRES.} 
	The theory guarantees that unresolved `bad' components (those with $\lambda_\ell^i\le\tau$) 
	are transferred to the coarse level.  Once the coarse space captures these, further enrichment has diminishing returns. 
	For higher $k$, local coercivity is already ensured by $H\sim c/k$ (cf.\ Lemma~\ref{lemma_3_4}), 
	and the residual high-frequency components are effectively damped by the local solves.
\end{itemize}

{For ease of reference, we summarise in Table~\ref{tab:summary_bounds} the sufficient $k$-explicit conditions on the subdomain diameter $H$ and the spectral threshold $\tau$ required to guarantee robust GMRES convergence for the original $\Delta$-GenEO method and for the proposed $\Delta_k$-GenEO variant. 
The table highlights the substantial reduction in wavenumber dependence achieved by the $k$-dependent formulation, while emphasising that all bounds remain sufficient rather than necessary.
\begin{table}[h]
\centering
\caption{Comparison of sufficient $k$-explicit conditions for robust GMRES convergence using $\Delta$-GenEO and $\Delta_k$-GenEO coarse spaces.}
\label{tab:summary_bounds}
\begin{tabular}{lcc}
\toprule
& \color{myred} $\Delta$-GenEO & \color{myred} $\Delta_k$-GenEO \\
\midrule
\color{myred} Subdomain diameter $H$ 
& \color{myred} $H \lesssim k^{-2}$ 
& \color{myred} $H \lesssim k^{-1}$ \\[0.5ex]
\color{myred} Spectral threshold $\tau$ 
& \color{myred} $\tau \gtrsim k^{8}$ 
& \color{myred} $\tau \gtrsim k^{2}$ \\
\bottomrule
\end{tabular}
\end{table}}
 
{While our bounds certify robustness under $H\sim k^{-1}$ and $\tau\gtrsim (1+\Cstab)^2 k^2$, 
they are not intended to be \emph{tight} for very large $k$.  Beyond some regime, SPD coarse spaces may cease 
to be the most effective choice; shifted/indefinite coarse operators or multilevel variants are natural successors. 
Our results should therefore be read as clarifying the range over which SPD GenEO-type spaces remain effective, 
and as a baseline for future indefinite coarse spaces.}

\section{Proof of the main result}
\label{sec:proof}

\noindent
Before proving the main result we need an additional lemma, {which contains the core estimates}.

\begin{lemma}[FoV lower bound and operator-norm bound]\label{ther: 4_1}
Under the assumptions of Theorem~\ref{thm:convergence} and the constraints \eqref{eq:st_constraints}, the projection operator
$T=\sum_{i=0}^N E_i T_i$ satisfies, for all \(u\in V^h\),
\begin{equation}\label{eq:FoV-lower}
	c_1\|u\|_{1, k}^2 \le (Tu,u)_{1, k}
\end{equation}
and
\begin{equation}\label{eq:op-norm-upper}
	\|Tu\|_{1, k} \le \sqrt{c_2}\,\|u\|_{1,k},
\end{equation}
where $c_1,c_2$ are as in \eqref{eq:c_1-c_2}.

\end{lemma}

\begin{proof}
\textit{Step 1 (FoV template: relate $\|u\|_{1, k}^2$ to $(Tu,u)_{1, k}$ plus two error terms).}
With $\Csd = 2+3\Lambda^2\Theta$ we obtain from Proposition~\ref{prop:ineqp} that
\[
	\|u\|_{1, k}^2 \le \Csd\,(Pu,u)_{1, k} = \Csd\sum_{i=0}^N (E_iP_iu,u)_{1, k}.
\]
For each $i\in\{0,\ldots,N\}$ we use the relation $b=a-k^2(\cdot,\cdot)$ and the definitions of $T_i$ and $P_i$ to obtain
\begin{align*}
	(E_iP_iu,u)_{1, k} &= (u,E_iP_iu)_{1, k}
	= b(u,E_iP_iu) + 2 k^2(u,E_iP_i)
	= b_{\Omega_i}(T_iu,P_iu) + 2 k^2(u,E_iP_iu)
	\\[0.5ex]
	&= (T_iu,P_iu)_{1, k, \Omega_i} - 2 k^2(T_iu,P_iu)_{\Omega_i} + 2 k^2(u,E_iP_iu)
	\\[0.5ex]
    &= (T_iu,P_iu)_{1, k, \Omega_i} - 2 k^2(E_i T_iu, E_i P_iu) + 2 k^2(u,E_iP_iu)
	\\[0.5ex]
	&= (E_iT_iu,u)_{1, k} + 2 k^2(u-E_i T_iu,E_i P_iu),
\end{align*}
which yields
\begin{equation}\label{eq:prelim}
	\|u\|_{1, k}^2 \le \Csd\biggl((Tu,u)_{1, k} + 2 k^2(u-T_0u,P_0u) + 2 k^2\sum_{i=1}^N (u-E_i T_iu,E_i P_iu)\biggr).
\end{equation}

We now bound the two perturbation terms.

\medskip

\noindent\textit{Step 2 (Coarse term bound via $T_0$ stability).}
It follows from Cauchy--Schwarz, \eqref{eq:friedrichs_energy_Omega} and Lemma~\ref{lem:T0-stab} that
\begin{align*}
	|2 k^2(u-T_0u,P_0u)| &\le 2 k^2\|T_0u-u\|\,\|P_0u\|
	\le 2 k\|T_0u-u\|\,\|P_0u\|_{1, k}
	\\[0.5ex]
	&\le 2k\Lambda\Theta^{1/2}(1+\Cstab)\|T_0u-u\|_{1, k}\|P_0u\|_{1, k}.
\end{align*}
The first inequality in \eqref{eq:st_constraints} implies that \eqref{eq:T0-extra} holds. Also we have the fact that $P_0$ is an orthogonal projection with respect to $(\cdot, \cdot)_{1, k}$,
and hence \eqref{eq:T0-energy} yields
\begin{equation}\label{eq:coarse-term}
	|2 k^2(u-T_0u,P_0u)| \le 4k\Lambda\Theta^{1/2}(1+\Cstab)\|u\|_{1, k}\|P_0\|_{1, k} \le 4k\Lambda\Theta^{1/2}(1+\Cstab)\|u\|_{1, k}^2.
\end{equation}

\noindent\textit{Step 3 (Local terms bound via $T_i$ stability and overlap).}
For $i\ge 1$, we have, by Cauchy--Schwarz, Friedrichs on $\Omega_i$, and Lemma~\ref{lemma_3_6},
\begin{align*}
	|2 k^2(u-E_i T_i u, E_i P_i u)|
	&\le 2 k^2 \,\|u- T_i u\|_{\Omega_i}\,\|P_i u\|_{\Omega_i}
	\le 2 k \,\|u- T_i u\|_{1, k, \Omega_i}\,\|P_i u\|_{\Omega_i}
	\\[0.5ex]
	&\le 2 k \Bigl(|u\|_{1, k, \Omega_i}+\|T_i u\|_{1, k, \Omega_i}\Bigr)\,\|P_i u\|_{\Omega_i}
	\le H k \sqrt{2}  \Bigl(\|u\|_{1, k, \Omega_i}+\|T_i u\|_{1, k, \Omega_i}\Bigr)\|P_i u\|_{1, k, \Omega_i}
    \\[0.5ex]
	&\le H k \sqrt{2}  \Bigl(\|u\|_{1, k, \Omega_i}+2 \| u\|_{1, k, \Omega_i}\Bigr)\|P_i u\|_{1, k, \Omega_i}
	= 3 H k \sqrt{2}\, \|u\|_{1, k, \Omega_i} \|P_i u\|_{1, k, \Omega_i}.
\end{align*}
Summing over $i\in\{1,\ldots,N\}$, applying Cauchy--Schwarz and using \eqref{overlap_stab}, \eqref{sum_norm_P_i} we obtain
\begin{align}
	\Bigg|2 k^2\sum_{i=1}^N (u-E_i T_i u, E_i P_i u)\bigg|
	&\le 3Hk\sqrt{2} \sum_{i=1}^N \|u\|_{1,k,\Omega_i} \|P_i u\|_{1,k,\Omega_i}
	\nonumber\\[0.5ex]
	&\le 3Hk\sqrt{2} \left( \sum_{i=1}^N \|u\|_{1, k, \Omega_i}^2 \right)^\frac{1}{2} \left( \sum_{i=1}^N  \|P_i u\|_{1, k, \Omega_i}^2 \right)^\frac{1}{2}
	\nonumber\\[0.5ex]
	&\le 3\Lambda Hk\sqrt{2}\, \|u\|_{1,k}^2.
	\label{eq:local-terms}
\end{align}

\noindent\textit{Step 4 (Combine and absorb).}
Inserting \eqref{eq:coarse-term} and \eqref{eq:local-terms} into \eqref{eq:prelim} we obtain
\begin{align*}
	\|u\|_{1,k}^2 &\le \Csd\Bigl((Tu,u)_{1,k} + 4k\Lambda\Theta^{1/2}(1+\Cstab)\|u\|_{1,k}^2 + 3\Lambda Hk\sqrt{2}\,\|u\|_{1,k}^2\Bigr)
	\\[1ex]
	&= \Csd(Tu,u)_{1,k} + \tfrac{1}{2}(s+t)\|u\|_{1,k}^2
	\le \Csd(Tu,u)_{1,k} + \max\{s,t\}\|u\|_{1,k}^2.
\end{align*}
Since $s,t<1$, this implies \eqref{eq:FoV-lower} with $c_1=\dfrac{1-\max\{s,t\}}{\Csd}$.

\medskip
\noindent\textit{Operator-norm bound.}
We write $Tu=T_0u+\sum_{i=1}^N E_iT_i u$ and use the overlap stability \eqref{overlap_stab} to find
\begin{equation}\label{eq:TuBound}
    \|Tu\|_{1, k}^2 \le 2\biggl(\|T_0u\|_{1,k}^2 + \Big\|\sum_{i=1}^N E_iT_i u\Big\|_{1,k}^2\biggr)
	\le 2\|T_0u\|_{1, k}^2 + 2\Lambda\sum_{i=1}^N \|T_i u\|_{1, k, \Omega_i}^2.
\end{equation}
For the first term on the right-hand side, we use Cauchy--Schwarz and Lemma~\ref{lem:T0-stab}, which yield
\begin{align*}
	\|T_0 u\|_{1,k}^2 &= (T_0 u, T_0 u)_{1,k} 
	= (T_0 u - u, T_0 u)_{1,k} + (u, T_0 u)_{1,k} 
	\\[0.5ex]
	&\le \|T_0 u - u\|_{1,k}\|T_0 u\|_{1,k} + \|u\|_{1,k}\|T_0 u\|_{1,k}
	\\[0.5ex]
	&\le 2\|u\|_{1,k}\|T_0 u\|_{1,k} + \|u\|_{1,k}\|T_0 u\|_{1,k} 
	= 3\|u\|_{1,k}\|T_0 u\|_{1,k}, 
\end{align*}
and hence
\begin{equation}\label{eq: 4_19}
	\|T_0 u\|_{1,k} \le 3 \|u\|_{1,k}. 
\end{equation}
For the second term on the right-hand side we use Lemma~\ref{lemma_3_6} to obtain
\begin{equation}
	\label{eq:second}
	\sum_{i=1} \|T_i u\|_{1,k,\Omega_i}^2 
	\le \sum_{i=1} 4\big\|u|_{\Omega_i}\big\|_{1,k,\Omega_i}^2 
	\le 4 \Lambda\|u\|_{1,k}^2. 
\end{equation}
Combining (\ref{eq: 4_19}) and \eqref{eq:second} with \eqref{eq:TuBound} we arrive at
\begin{equation*}
	\|Tu\|_{1,k}^2 \le 2\times9\|u\|_{1,k}^2 + 2\times4\Lambda^2\|u\|_{1,k}^2 
	= (18 + 8\Lambda^2)\|u\|_{1,k}^2,
\end{equation*} 
which is \eqref{eq:op-norm-upper}.
\end{proof}

We are now able to finish the proof of Theorem~\ref{thm:convergence}.

\begin{proof}[Proof of Theorem~\ref{thm:convergence}]
Let $\mathcal{M}\assign\mathbf{M}_{AS,2}^{-1}\mathbf{B}$ and endow $\mathbb{R}^n$ with the
$\mathbf{D}_k$-inner product $\langle\cdot,\cdot\rangle_{\mathbf D_k}$ and norm
$\|\cdot\|_{\mathbf D_k}$.  By Proposition~\ref{prop:proj}, for all $u,v\in V^h$ with
coefficient vectors $\mathbf u,\mathbf v$,
\[
	(Tu,v)_{1, k} = \langle \mathcal{M}\mathbf u,\mathbf v\rangle_{\mathbf D_k}.
\]
Taking $v=u$ and using \eqref{eq:FoV-lower} we obtain the FoV distance (from the origin)
\[
	\inf_{\|\mathbf u\|_{\mathbf D_k}=1}\operatorname{Re}\,\langle \mathcal{M}\mathbf u,\mathbf u\rangle_{\mathbf D_k}
	\ge c_1.
\]
Next, \eqref{eq:op-norm-upper} gives
\[
	\|\mathcal{M}\mathbf u\|_{\mathbf D_k}^2
	= \langle \mathcal{M}\mathbf u,\mathcal{M}\mathbf u\rangle_{\mathbf D_k}
	= (Tu,Tu)_{1, k} =\|Tu\|_{1, k}^2
	\le c_2\,\|u\|_{1, k}^2
	= c_2\|\mathbf u\|_{\mathbf D_k}^2;
\]
hence the operator norm satisfies $\|\mathcal{M}\|_{\mathbf D_k}\le \sqrt{c_2}$.

By the Elman field-of-values estimate for (right-)preconditioned GMRES
\cite{Elman:1983:VIM}, the residuals satisfy
\[
	\|\mathbf r^{(m)}\|_{\mathbf D_k}^2
	\le
	\biggl(1-\Bigl(\frac{\delta}{\beta}\Bigr)^2\biggr)^{m}\,\|\mathbf r^{(0)}\|_{\mathbf D_k}^2,
	\qquad
	\delta\assign\inf_{\|\mathbf u\|_{\mathbf D_k}=1}\operatorname{Re}\,\langle \mathcal{M}\mathbf u,\mathbf u\rangle_{\mathbf D_k},\quad
	\beta\assign\|\mathcal{M}\|_{\mathbf D_k}.
\]
Using $\delta\ge c_1$ and $\beta^2\le c_2$ we obtain
\[
	\|\mathbf r^{(m)}\|_{\mathbf D_k}^2
	\le
	\Bigl(1-\frac{c_1^2}{c_2}\Bigr)^{m}\|\mathbf r^{(0)}\|_{\mathbf D_k}^2,
\]
which is the claimed convergence bound (with the stated $c_1$ and $c_2$).
\end{proof}

\section{Numerical results} 
\label{sec:numerics}
The aim of this section is to provide a comprehensive assessment of the performance and robustness 
of the $\Delta_k$-GenEO coarse space when combined with the two-level additive Schwarz (AS) preconditioner, 
and to compare its behaviour with that of the classical one-level AS method. 
While theoretical analysis provides insight into the spectral properties of these preconditioners, 
the actual performance in practice is influenced by a number of factors that are difficult to quantify a priori, 
such as the frequency, the domain partitioning, the heterogeneity of the medium, and the chosen eigenvalue threshold for coarse space construction.  

{The numerical experiments presented in this section are not designed to strictly enforce the sufficient conditions obtained in the theoretical results. 
Instead, they deliberately explore a range of regimes that extend beyond these conservative bounds in order to assess the practical robustness of the proposed $\Delta_k$-GenEO method. 
Many experiments intentionally violate one or more of the theoretical assumptions. 
The resulting performance illustrates that the method remains effective well beyond the strict regime guaranteed by the analysis, thereby highlighting the gap between sufficient theoretical conditions and observed behaviour in practice.
}

{In what follows}, we conduct a series of numerical experiments. Our tests are restricted to two representative cases:
\begin{enumerate}
\item \textbf{Homogeneous test case.}  
A benchmark posed on the unit square with homogeneous Dirichlet boundary conditions and a single interior point source. 
This case allows for the performance of the preconditioners to be tested in a controlled setting.  

\item \textbf{Heterogeneous test case.}  
A layered medium with strong contrasts, designed to probe the robustness of the method under medium heterogeneity, 
and to examine how the coarse space size and convergence are affected in more challenging settings.  
\end{enumerate}

Both test cases allow us to study how iteration counts and coarse space dimension scale with respect to the 
number of subdomains and the frequency. The influence of the following parameters is reported:  
\begin{itemize}
  \item subdomain diameter measured in wavelengths,
  \item global and per-subdomain coarse space dimensions,
  \item eigenvalue threshold $\tau$ for coarse space construction,
\end{itemize}

\subsection{Implementation Details}
All computations are performed in \texttt{FreeFEM} using the \texttt{ffddm} framework \cite{FFD:Tournier:2019},
which provides highly flexible Message Passing Interface (MPI)-parallel scripts for domain decomposition methods. 

\paragraph{Domain partition and overlap}
The computational domain $\Omega$ is partitioned into $N$ non-overlapping square subdomains $\{\Omega_i'\}_{i=1}^N$.
Unless otherwise stated, an overlap is introduced by extending each $\Omega_i'$ 
by a minimum number of layers of finite elements, yielding overlapping subdomains $\{\Omega_i\}_{i=1}^N$.  

\paragraph{Setup: one-level preconditioner}
In the one-level AS method, we perform the following
\begin{enumerate}
  \item Assemble the local stiffness and mass matrices on $\Omega_i$ corresponding to the bilinear forms $a_{\Omega_i}$ and $b_{\Omega_i}$.
  \item Compute and store a sparse LU factorisation of the resulting local system using \texttt{MUMPS}.
\end{enumerate}

\paragraph{Setup: two-level ($\Delta_k$-GenEO) coarse space}
The two-level method augments the one-level setup with a global coarse correction space. On each rank $i$ we:
\begin{enumerate}
  \item 
  Assemble the local generalised eigenproblem
  where $\Xi_i$ is the partition of unity operator.
  \item 
  Solve the eigenproblem using \texttt{SLEPc}, and select the first $m_i$ eigenmodes according to a fixed 
  spectral threshold ($\lambda_\ell^i \leq \tau$). 
  \item 
  Form the corresponding coarse basis vectors $E_i \Xi_i(p_\ell^i)$ and contribute to the global coarse operator $\mathbf{B}_0$, 
  defined as the $b(\cdot,\cdot)$-form restricted to $V_0$.
\end{enumerate}

\paragraph{Software stack}
For the numerical linear algebra components, we employ the sparse direct solver \texttt{MUMPS}~\cite{amestoy2001fully} 
to factorise both local and coarse matrices. The local spectral problems arising in the construction 
of the $\Delta_k$-GenEO coarse space are solved using the \texttt{SLEPc} library~\cite{hernandez2005ssf}, 
which provides scalable eigensolvers for large sparse eigenvalue problems. Iterative solutions are obtained with GMRES. 
Unless otherwise stated, GMRES is not restarted; iterations are terminated once the relative residual 
is reduced below $10^{-6}$, or the number of iterations reaches 200.

\paragraph{Preconditioner application (solve phase)}
Each iteration of the preconditioned GMRES method involves:
\begin{enumerate}
  \item \emph{Local solves:} one forward/backward substitution with the stored LU factorisation on each $\Omega_i$, executed in parallel across all MPI ranks. 
  \item \emph{Coarse solve:} one forward/backward substitution with the factorised coarse matrix $\mathbf{B}_0$, distributed over the coarse ranks, 
  accompanied by the necessary global communication (scatter and gather).
\end{enumerate}
{\textit{Algorithmic complexity.}
From a parallel complexity perspective, the implementation naturally separates into a \emph{setup phase} and an \emph{application phase}, 
whose scaling behaviour is qualitatively different. 
During the setup phase, all local tasks---assembly of local matrices, solution of the local generalised eigenvalue problems, 
and factorisations---are embarrassingly parallel across subdomains. 
The only global communication required at this stage is the assembly of the coarse operator $\mathbf{B}_0$, 
whose size scales with the coarse space dimension $\dim V_0$.}

{During the application phase, each GMRES iteration consists of fully parallel local solves 
on the subdomains and a single coarse-level solve, accompanied by standard global scatter/gather communication. 
This structure is typical of two-level additive Schwarz methods. 
As a result, the dominant communication costs are associated with the coarse solve and are governed by the size of the coarse space 
rather than the fine-grid problem size.
}

{We emphasise that the present implementation is not optimised and does not aim to fully exploit all available parallelisation, 
communication, or solver-tuning strategies. 
Consequently, wall-clock timings, communication overheads, and strong or weak scaling results would be highly implementation- and hardware-dependent 
and would not be representative of the intrinsic capabilities of the method. 
For this reason, the numerical study focuses on algorithmic complexity and qualitative scaling behaviour---such as iteration counts, 
coarse space growth, and robustness with respect to frequency and heterogeneity---which are the primary indicators of performance 
for spectral coarse spaces.
}

It should be also noted that due to hardware limitations, we have not been able to utilise full parallelisation. Due to the restricted number of CPU cores available, we enabled oversubscription, which allows more MPI ranks than physical CPU cores, meaning multiple ranks time-share the same core. This preserves functionality, but can reduce performance. 

\subsubsection{Homogeneous problem}
\label{sec:HomProb}
In this part we base the numerical experiments on the model problem \eqref{eq:problem} in 2D, defined on the unit square $\Omega = (0,1)^2$ and with $A(\bx) = 1$. We impose Dirichlet boundary conditions on all sides of the domain. A point source is located in the centre of the domain at $(\frac{1}{2},\frac{1}{2})$ and provides the forcing function $f$. {The point source is numerically modelled by a Gaussian function:
\(f(x,y)=10^4 \exp(-10^3[(x-\frac{1}{2})^2+(y-\frac{1}{2})^2])\).} {We note that the specific choice of a Gaussian approximation for the point source is made for numerical convenience. Since the construction of the $\Delta_k$-GenEO coarse space and the associated preconditioner are independent of the right-hand side, we do not expect qualitatively different behaviour for other smooth, localised, or volume source terms.}

To discretise the problem, we triangulate $\Omega$ using a Cartesian grid with spacing $h$ and alternating diagonals to form a simplicial mesh. The local wavenumber, $k$ is constant and the wavelength used to measure geometrical parameters is \(\lambda=2\pi/k\). The discrete problem \eqref{eq:lin_sys} is assembled using a $\mathbb{P}_1$ finite element approximation on this mesh.   To mitigate the \emph{pollution effect}, we choose the wavenumber $k$ and the mesh size $h$ simultaneously so that the dimensionless quantity $k h$ remains sufficiently small.  
In practice, this is enforced by fixing a minimum number of grid points per wavelength $\lambda$. Here, we ensure that { 
$k h \lesssim 0.1$, which corresponds to at least $25 \pi$ points per wavelength $\lambda$. This is widely considered a very fine mesh, with most real world applications accepting $10$ points per wave length. However, this value has been chosen to ensure an accurate solution whilst maintaining a solvable problem size.

\begin{table}[H]
    \centering
    \hspace*{-2cm}
\begin{tabular}{ccccccc|c|ccc|ccc}
\multicolumn{7}{c|}{} & One Level & \multicolumn{3}{c}{$\Delta$-GenEO} & \multicolumn{3}{c}{$\Delta_k$-GenEO} \\
$N$ & $k$ & $h^{-1}$ & $n$ & $L$ & $H$ & $n_\text{loc}$ & It & It & CS & CS$_\text{loc}$ & It & CS & CS$_\text{loc}$  \\
\hline
16 & 20 & 240 & 58081 & 3.18 & 0.82 & 3969 & 41 & 36 & 400 & 25 & 22 & 272 & 17 \\
 & 40 & 480 & 231361 &  6.37 & 1.62 & 15129 & 107 & 92 & 804 & 50 & 73 & 580 & 36 \\
 & 60 & 720 & 519841 &  9.55 & 2.41 & 33489 & 153 & 156 & 1232 & 77 & 78 & 976 & 61 \\
 & 80 & 960 & 923521 &  12.73 & 3.21 & 59049 & $>$200 & $>$200 & 1640 & 102 & 155 & 1424 & 89 \\
 & 100 & 1200 & 1442401 &  15.92 & 4.01 & 91809 & $>$200 & $>$200 & 2052 & 128 & 170 & 1932 & 121 \\
\hline
36 & 20 & 240 & 58081 & 3.18 & 0.56 & 1849 & 79 & 23 & 684 & 19 & 23 & 428 & 12 \\
 & 40 & 480 & 231361 &  6.37 & 1.09 & 6889 & $>$200 & 130 & 1356 & 38 & 112 & 860 & 24 \\
 & 60 & 720 & 519841 &  9.55 & 1.62 & 15129 & $>$200 & $>$200 & 2020 & 56 & 144 & 1444 & 40 \\
 & 80 & 960 & 923521 &  12.73 & 2.15 & 26569 & $>$200 & $>$200 & 2740 & 76 & $>$200 & 2028 & 56 \\
 & 100 & 1200 & 1442401 &  15.92 & 2.68 & 41209 & $>$200 & $>$200 & 3428 & 95 & $>$200 & 2752 & 76 \\
\hline
64 & 20 & 240 & 58081 & 3.18 & 0.42 & 1089 & 92 & 22 & 896 & 14 & 23 & 608 & 10 \\
 & 40 & 480 & 231361 &  6.37 & 0.82 & 3969 & 175 & 133 & 1856 & 29 & 57 & 1248 & 20 \\
 & 60 & 720 & 519841 &  9.55 & 1.22 & 8649 & $>$200 & $>$200 & 2848 & 44 & 181 & 1884 & 29 \\
 & 80 & 960 & 923521 &  12.73 & 1.62 & 15129 & $>$200 & $>$200 & 3780 & 59 & $>$200 & 2692 & 42 \\
 & 100 & 1200 & 1442401 &  15.92 & 2.02 & 23409 & $>$200 & $>$200 & 4768 & 74 & $>$200 & 3488 & 54 \\
\hline
100 & 20 & 240 & 58081 & 3.18 & 0.34 & 729 & 96 & 21 & 1176 & 12 & 22 & 684 & 7 \\
 & 40 & 480 & 231361 &  6.37 & 0.66 & 2601 & $>$200 & 108 & 2420 & 24 & 105 & 1472 & 15 \\
 & 60 & 720 & 519841 &  9.55 & 0.98 & 5625 & $>$200 & $>$200 & 3600 & 36 & 124 & 2360 & 24 \\
 & 80 & 960 & 923521 &  12.73 & 1.30 & 9801 & $>$200 & $>$200 & 4940 & 49 & 192 & 3440 & 34 \\
 & 100 & 1200 & 1442401 &  15.92 & 1.62 & 15129 & $>$200 & $>$200 & 6084 & 61 & $>$200 & 4324 & 43 \\
\hline
144 & 20 & 240 & 58081 & 3.18 & 0.29 & 529 & 104 & 18 & 1440 & 10 & 19 & 1012 & 7 \\
 & 40 & 480 & 231361 &  6.37 & 0.56 & 1849 & $>$200 & 67 & 3024 & 21 & 63 & 1928 & 13 \\
 & 60 & 720 & 519841 &  9.55 & 0.82 & 3969 & $>$200 & $>$200 & 4368 & 30 & 78 & 2928 & 20 \\
 & 80 & 960 & 923521 &  12.73 & 1.09 & 6889 & $>$200 & $>$200 & 5952 & 41 & $>$200 & 3800 & 26 \\
 & 100 & 1200 & 1442401 &  15.92 & 1.35 & 10609 & $>$200 & $>$200 & 7392 & 51 & $>$200 & 5140 & 36 \\
\hline
\end{tabular}
    \caption{Results showing the dimension of the fine mesh ($n$), the diameter of the domain measured in wavelengths ($L$), the diameter of the sub-domains measured in wavelengths ($H$), the number of local degrees of freedom ($n_\text{loc}$), the GMRES iteration count (It.) for the one-level and two level methods using the $\Delta$-GenEO and $\Delta_k$-GenEO coarse space, the dimension of the coarse space (CS) and averaged number of contributions to the coarse space per subdomain (CS$_\text{loc}$) for the $\Delta$-GenEO and $\Delta_k$-GenEO. These results are for the homogeneous media test case. Using $\tau = 0.7$ and $\tau=0.5$ for $\Delta$-GenEO and $\Delta_k$-GenEO respectively.}
    \label{table:DeltaHomog}
\end{table}

\begin{figure}[H]
    \centering
    \begin{subfigure}[t]{0.48\textwidth}
        \centering
        \begin{tikzpicture}
\begin{axis}[
  xmode=log,
  ymode=log,
  enlarge x limits=0.05,
  enlarge y limits=0.05,
  xmin=36,
  xmax=1540,
  ymin=16,
  ymax=104,
  grid=both,
  minor grid style={opacity=0.3},
  major grid style={opacity=0.3},
  xlabel={CS Size},
  ylabel={Iteration Count},
  width=\textwidth,
  height=\textwidth,
  legend pos=north east,
  legend style={fill=none, draw=none, font=\scriptsize, text opacity=1}
]
\addplot+[mark=square*, mark options={solid, fill = blue}, mark size=1.8pt, line width=1pt, color=blue, opacity=1, style=solid] coordinates {(36,48) (84,52) (132,41) (180,37) (240,36) (400,36)};
\addplot+[mark=square*, mark options={solid, fill = orange}, mark size=1.8pt, line width=1pt, color=orange, opacity=1, style=solid] coordinates {(84,66) (120,46) (220,28) (324,26) (396,24) (684,23)};
\addplot+[mark=square*, mark options={solid, fill = green!60!black}, mark size=1.8pt, line width=1pt, color=green!60!black, opacity=1, style=solid] coordinates {(136,73) (196,39) (288,33) (420,25) (608,23) (896,22)};
\addplot+[mark=square*, mark options={solid, fill = red!80!black}, mark size=1.8pt, line width=1pt, color=red!80!black, opacity=1, style=solid] coordinates {(100,99) (324,35) (324,35) (620,26) (684,22) (1176,21)};
\addplot+[mark=square*, mark options={solid, fill = cyan}, mark size=1.8pt, line width=1pt, color=cyan, opacity=1, style=solid] coordinates {(144,99) (384,39) (484,30) (568,28) (912,23) (1440,18)};
\addplot+[mark=*, mark options={solid, fill = blue}, mark size=1.8pt, line width=1pt, color=blue, opacity=1, style=dashed] coordinates {(36,48) (84,51) (132,39) (180,35) (272,22) (340,21) (464,16)};
\addplot+[mark=*, mark options={solid, fill = orange}, mark size=1.8pt, line width=1pt, color=orange, opacity=1, style=dashed] coordinates {(84,62) (140,38) (220,27) (340,24) (428,23) (516,22) (720,16)};
\addplot+[mark=*, mark options={solid, fill = green!60!black}, mark size=1.8pt, line width=1pt, color=green!60!black, opacity=1, style=dashed] coordinates {(136,73) (196,38) (288,33) (420,25) (608,23) (736,23) (920,22)};
\addplot+[mark=*, mark options={solid, fill = red!80!black}, mark size=1.8pt, line width=1pt, color=red!80!black, opacity=1, style=dashed] coordinates {(100,100) (324,34) (360,33) (620,26) (684,22) (980,21) (1180,21)};
\addplot+[mark=*, mark options={solid, fill = cyan}, mark size=1.8pt, line width=1pt, color=cyan, opacity=1, style=dashed] coordinates {(144,99) (384,38) (484,29) (672,26) (1012,19) (1056,19) (1540,18)};
\addplot+[mark=none, black, dashed, line width=1.2pt, opacity=0.2] coordinates {(36,41) (1540,41)};
\addplot+[mark=none, black, dashed, line width=1.2pt, opacity=0.4] coordinates {(36,79) (1540,79)};
\addplot+[mark=none, black, dashed, line width=1.2pt, opacity=0.6] coordinates {(36,92) (1540,92)};
\addplot+[mark=none, black, dashed, line width=1.2pt, opacity=0.8] coordinates {(36,96) (1540,96)};
\addplot+[mark=none, black, dashed, line width=1.2pt, opacity=1.0] coordinates {(36,104) (1540,104)};
\end{axis}
\end{tikzpicture}
    \end{subfigure}   
    \hfill
    \begin{subfigure}[t]{0.48\textwidth}
        \centering
        \begin{tikzpicture}
\begin{axis}[
  xmode=log,
  ymode=log,
  enlarge x limits=0.05,
  enlarge y limits=0.05,
  xmin=0.1,
  xmax=0.7,
  ymin=16,
  ymax=104,
  grid=both,
  minor grid style={opacity=0.3},
  major grid style={opacity=0.3},
  xlabel={$\tau$},
  ylabel={Iteration Count},
  width=\textwidth,
  height=\textwidth,
  legend pos=north east,
  legend style={fill=none, draw=none, font=\scriptsize, text opacity=1}
]
\addplot+[mark=square*, mark options={solid, fill = blue}, mark size=1.8pt, line width=1pt, color=blue, opacity=1, style=solid] coordinates {(0.1,48) (0.2,52) (0.3,41) (0.4,37) (0.5,36) (0.7,36)};
\addplot+[mark=square*, mark options={solid, fill = orange}, mark size=1.8pt, line width=1pt, color=orange, opacity=1, style=solid] coordinates {(0.1,66) (0.2,46) (0.3,28) (0.4,26) (0.5,24) (0.7,23)};
\addplot+[mark=square*, mark options={solid, fill = green!60!black}, mark size=1.8pt, line width=1pt, color=green!60!black, opacity=1, style=solid] coordinates {(0.1,73) (0.2,39) (0.3,33) (0.4,25) (0.5,23) (0.7,22)};
\addplot+[mark=square*, mark options={solid, fill = red!80!black}, mark size=1.8pt, line width=1pt, color=red!80!black, opacity=1, style=solid] coordinates {(0.1,99) (0.2,35) (0.3,35) (0.4,26) (0.5,22) (0.7,21)};
\addplot+[mark=square*, mark options={solid, fill = cyan}, mark size=1.8pt, line width=1pt, color=cyan, opacity=1, style=solid] coordinates {(0.1,99) (0.2,39) (0.3,30) (0.4,28) (0.5,23) (0.7,18)};
\addplot+[mark=*, mark options={solid, fill = blue}, mark size=1.8pt, line width=1pt, color=blue, opacity=1, style=dashed] coordinates {(0.1,48) (0.2,51) (0.3,39) (0.4,35) (0.5,22) (0.6,21) (0.7,16)};
\addplot+[mark=*, mark options={solid, fill = orange}, mark size=1.8pt, line width=1pt, color=orange, opacity=1, style=dashed] coordinates {(0.1,62) (0.2,38) (0.3,27) (0.4,24) (0.5,23) (0.6,22) (0.7,16)};
\addplot+[mark=*, mark options={solid, fill = green!60!black}, mark size=1.8pt, line width=1pt, color=green!60!black, opacity=1, style=dashed] coordinates {(0.1,73) (0.2,38) (0.3,33) (0.4,25) (0.5,23) (0.6,23) (0.7,22)};
\addplot+[mark=*, mark options={solid, fill = red!80!black}, mark size=1.8pt, line width=1pt, color=red!80!black, opacity=1, style=dashed] coordinates {(0.1,100) (0.2,34) (0.3,33) (0.4,26) (0.5,22) (0.6,21) (0.7,21)};
\addplot+[mark=*, mark options={solid, fill = cyan}, mark size=1.8pt, line width=1pt, color=cyan, opacity=1, style=dashed] coordinates {(0.1,99) (0.2,38) (0.3,29) (0.4,26) (0.5,19) (0.6,19) (0.7,18)};
\addplot+[mark=none, black, dashed, line width=1.2pt, opacity=0.2] coordinates {(0.1,41) (0.7,41)};
\addplot+[mark=none, black, dashed, line width=1.2pt, opacity=0.4] coordinates {(0.1,79) (0.7,79)};
\addplot+[mark=none, black, dashed, line width=1.2pt, opacity=0.6] coordinates {(0.1,92) (0.7,92)};
\addplot+[mark=none, black, dashed, line width=1.2pt, opacity=0.8] coordinates {(0.1,96) (0.7,96)};
\addplot+[mark=none, black, dashed, line width=1.2pt, opacity=1.0] coordinates {(0.1,104) (0.7,104)};
\end{axis}
\end{tikzpicture}
    \end{subfigure} 
    \begin{subfigure}[t]{0.48\textwidth}
        \centering
        \begin{tikzpicture}
\begin{axis}[
  xmode=log,
  ymode=log,
  enlarge x limits=0.05,
  enlarge y limits=0.05,
  xmin=228,
  xmax=8688,
  ymin=38,
  ymax=200,
  grid=both,
  minor grid style={opacity=0.3},
  major grid style={opacity=0.3},
  xlabel={CS Size},
  ylabel={Iteration Count},
  width=\textwidth,
  height=\textwidth,
  legend pos=north east,
  legend style={fill=none, draw=none, font=\small, text opacity=1},
  legend columns = 5,
  legend to name = domains
]
\addplot+[mark=square*, mark options={solid, fill = blue}, mark size=1.8pt, line width=1pt, color=blue, opacity=1] coordinates {(228,200) (468,200) (708,200) (968,200) (1264,200) (2052,200)};
\addlegendentry{$\Delta$-GenEO (16)}
\addplot+[mark=square*, mark options={solid, fill = orange}, mark size=1.8pt, line width=1pt, color=orange, opacity=1] coordinates {(340,200) (756,200) (1180,200) (1596,200) (2124,200) (3428,200)};
\addlegendentry{$\Delta$-GenEO (36)}
\addplot+[mark=square*, mark options={solid, fill = green!60!black}, mark size=1.8pt, line width=1pt, color=green!60!black, opacity=1] coordinates {(476,200) (1092,200) (1596,200) (2212,200) (2912,200) (4768,200)};
\addlegendentry{$\Delta$-GenEO (64)}
\addplot+[mark=square*, mark options={solid, fill = red!80!black}, mark size=1.8pt, line width=1pt, color=red!80!black, opacity=1] coordinates {(684,200) (1404,200) (2124,200) (2844,200) (3828,200) (6084,200)};
\addlegendentry{$\Delta$-GenEO (100)}
\addplot+[mark=square*, mark options={solid, fill = cyan}, mark size=1.8pt, line width=1pt, color=cyan, opacity=1] coordinates {(1012,200) (1544,200) (2596,200) (3652,200) (4708,200) (7392,200)};
\addlegendentry{$\Delta$-GenEO (144)}
\addplot+[mark=*, mark options={solid, fill = blue}, mark size=1.8pt, line width=1pt, color=blue, opacity=1, style=dashed] coordinates {(300,200) (616,200) (988,200) (1412,200) (1932,170) (2672,68) (3688,52)};
\addlegendentry{$\Delta_k$-GenEO (16)}
\addplot+[mark=*, mark options={solid, fill = orange}, mark size=1.8pt, line width=1pt, color=orange, opacity=1, style=dashed] coordinates {(480,200) (928,200) (1416,200) (2024,200) (2752,200) (3684,90) (4936,78)};
\addlegendentry{$\Delta_k$-GenEO (36)}
\addplot+[mark=*, mark options={solid, fill = green!60!black}, mark size=1.8pt, line width=1pt, color=green!60!black, opacity=1, style=dashed] coordinates {(644,200) (1156,200) (1920,200) (2560,200) (3488,200) (4616,88) (6272,38)};
\addlegendentry{$\Delta_k$-GenEO (64)}
\addplot+[mark=*, mark options={solid, fill = red!80!black}, mark size=1.8pt, line width=1pt, color=red!80!black, opacity=1, style=dashed] coordinates {(684,200) (1408,200) (2224,200) (3216,200) (4324,200) (5740,80) (7520,44)};
\addlegendentry{$\Delta_k$-GenEO (100)}
\addplot+[mark=*, mark options={solid, fill = cyan}, mark size=1.8pt, line width=1pt, color=cyan, opacity=1, style=dashed] coordinates {(1012,200) (1628,200) (2740,200) (3796,200) (5140,200) (6576,197) (8688,55)};
\addlegendentry{$\Delta_k$-GenEO (144)}
\addplot+[mark=none, black, dashed, line width=1.2pt, opacity=0.2] coordinates {(228,200) (8688,200)};
\addlegendentry{One-level (16)}
\addplot+[mark=none, black, dashed, line width=1.2pt, opacity=0.4] coordinates {(228,200) (8688,200)};
\addlegendentry{One-level (36)}
\addplot+[mark=none, black, dashed, line width=1.2pt, opacity=0.6] coordinates {(228,200) (8688,200)};
\addlegendentry{One-level (64)}
\addplot+[mark=none, black, dashed, line width=1.2pt, opacity=0.8] coordinates {(228,200) (8688,200)};
\addlegendentry{One-level (100)}
\addplot+[mark=none, black, dashed, line width=1.2pt, opacity=1.0] coordinates {(228,200) (8688,200)};
\addlegendentry{One-level (144)}
\end{axis}
\end{tikzpicture}
    \end{subfigure} 
    \hfill
    \begin{subfigure}[t]{0.48\textwidth}
        \centering
        \begin{tikzpicture}
\begin{axis}[
  xmode=log,
  ymode=log,
  enlarge x limits=0.05,
  enlarge y limits=0.05,
  xmin=0.1,
  xmax=0.7,
  ymin=38,
  ymax=200,
  grid=both,
  minor grid style={opacity=0.3},
  major grid style={opacity=0.3},
  xlabel={$\tau$},
  ylabel={Iteration Count},
  width=\textwidth,
  height=\textwidth,
  legend pos=north east,
  legend style={fill=none, draw=none, font=\scriptsize, text opacity=1}
]
\addplot+[mark=*, mark options={solid, fill = blue}, mark size=1.8pt, line width=1pt, color=blue, opacity=1] coordinates {(0.1,200) (0.2,200) (0.3,200) (0.4,200) (0.5,200) (0.7,200)};
\addplot+[mark=*, mark options={solid, fill = orange}, mark size=1.8pt, line width=1pt, color=orange, opacity=1] coordinates {(0.1,200) (0.2,200) (0.3,200) (0.4,200) (0.5,200) (0.7,200)};
\addplot+[mark=*, mark options={solid, fill = green!60!black}, mark size=1.8pt, line width=1pt, color=green!60!black, opacity=1] coordinates {(0.1,200) (0.2,200) (0.3,200) (0.4,200) (0.5,200) (0.7,200)};
\addplot+[mark=*, mark options={solid, fill = red!80!black}, mark size=1.8pt, line width=1pt, color=red!80!black, opacity=1] coordinates {(0.1,200) (0.2,200) (0.3,200) (0.4,200) (0.5,200) (0.7,200)};
\addplot+[mark=*, mark options={solid, fill = cyan}, mark size=1.8pt, line width=1pt, color=cyan, opacity=1] coordinates {(0.1,200) (0.2,200) (0.3,200) (0.4,200) (0.5,200) (0.7,200)};
\addplot+[mark=*, mark options={solid, fill = blue}, mark size=1.8pt, line width=1pt, color=blue, opacity=1, style=dashed] coordinates {(0.1,200) (0.2,200) (0.3,200) (0.4,200) (0.5,170) (0.6,68) (0.7,52)};
\addplot+[mark=*, mark options={solid, fill = orange}, mark size=1.8pt, line width=1pt, color=orange, opacity=1, style=dashed] coordinates {(0.1,200) (0.2,200) (0.3,200) (0.4,200) (0.5,200) (0.6,90) (0.7,78)};
\addplot+[mark=*, mark options={solid, fill = green!60!black}, mark size=1.8pt, line width=1pt, color=green!60!black, opacity=1, style=dashed] coordinates {(0.1,200) (0.2,200) (0.3,200) (0.4,200) (0.5,200) (0.6,88) (0.7,38)};
\addplot+[mark=*, mark options={solid, fill = red!80!black}, mark size=1.8pt, line width=1pt, color=red!80!black, opacity=1, style=dashed] coordinates {(0.1,200) (0.2,200) (0.3,200) (0.4,200) (0.5,200) (0.6,80) (0.7,44)};
\addplot+[mark=*, mark options={solid, fill = cyan}, mark size=1.8pt, line width=1pt, color=cyan, opacity=1, style=dashed] coordinates {(0.1,200) (0.2,200) (0.3,200) (0.4,200) (0.5,200) (0.6,197) (0.7,55)};
\addplot+[mark=none, black, dashed, line width=1.2pt, opacity=0.2] coordinates {(0.1,200) (0.7,200)};
\addplot+[mark=none, black, dashed, line width=1.2pt, opacity=0.4] coordinates {(0.1,200) (0.7,200)};
\addplot+[mark=none, black, dashed, line width=1.2pt, opacity=0.6] coordinates {(0.1,200) (0.7,200)};
\addplot+[mark=none, black, dashed, line width=1.2pt, opacity=0.8] coordinates {(0.1,200) (0.7,200)};
\addplot+[mark=none, black, dashed, line width=1.2pt, opacity=1.0] coordinates {(0.1,200) (0.7,200)};
\end{axis}
\end{tikzpicture}
    \end{subfigure} 
    \ref{domains}
    \caption{Influence of the coarse space size (left) and threshold choice (right) on the iteration count for the homogeneous media test case with $k$ = 20 (top) and $k$ = 100 (bottom). The number in brackets indicates the number of subdomains.}
    \label{fig:DeltaHomCSIT}
\end{figure}

\begin{figure}[H]
    \centering
    \begin{subfigure}[t]{0.48\textwidth}
        \centering
        \begin{tikzpicture}
\begin{axis}[
  xmode=log,
  ymode=log,
  enlarge x limits=0.05,
  enlarge y limits=0.05,
  xmin=0.29178406,
  xmax=0.82230054,
  ymin=84,
  ymax=1540,
  grid=both,
  minor grid style={opacity=0.3},
  major grid style={opacity=0.3},
  xlabel={$H$ (subdomain side length in wavelengths)},
  ylabel={CS Size},
  width=\textwidth,
  height=\textwidth,
  legend pos=north east,
  legend style={fill=none, draw=none, font=\scriptsize, text opacity=1}
]
\addplot+[mark=square*, mark options={solid, fill = blue}, mark size=1.8pt, line width=1pt, color=blue, opacity=1, style=solid] coordinates {(0.29178406,384) (0.34483571,324) (0.42441318,196) (0.5570423,120) (0.6631456,105) (0.82230054,84)};
\addplot+[mark=square*, mark options={solid, fill = orange}, mark size=1.8pt, line width=1pt, color=orange, opacity=1, style=solid] coordinates {(0.29178406,568) (0.34483571,620) (0.42441318,420) (0.5570423,324) (0.6631456,240) (0.82230054,180)};
\addplot+[mark=square*, mark options={solid, fill = cyan}, mark size=1.8pt, line width=1pt, color=cyan, opacity=1, style=solid] coordinates {(0.29178406,1440) (0.34483571,1176) (0.42441318,896) (0.5570423,684) (0.6631456,535) (0.82230054,400)};
\addplot+[mark=*, mark options={solid, fill = blue}, mark size=1.8pt, line width=1pt, color=blue, opacity=1, style=dashed] coordinates {(0.82230054,84) (0.6631456,135) (0.5570423,140) (0.42441318,196) (0.34483571,324) (0.29178406,384)};
\addplot+[mark=*, mark options={solid, fill = orange}, mark size=1.8pt, line width=1pt, color=orange, opacity=1, style=dashed] coordinates {(0.82230054,180) (0.6631456,252) (0.5570423,340) (0.42441318,420) (0.34483571,620) (0.29178406,672)};
\addplot+[mark=*, mark options={solid, fill = cyan}, mark size=1.8pt, line width=1pt, color=cyan, opacity=1, style=dashed] coordinates {(0.82230054,464) (0.6631456,569) (0.5570423,720) (0.42441318,920) (0.34483571,1180) (0.29178406,1540)};
\end{axis}
\end{tikzpicture}
    \end{subfigure}   
    \hfill
    \begin{subfigure}[t]{0.48\textwidth}
        \centering
        \begin{tikzpicture}
\begin{axis}[
  xmode=log,
  ymode=log,
  enlarge x limits=0.05,
  enlarge y limits=0.05,
  xmin=0.29178406,
  xmax=0.82230054,
  ymin=16,
  ymax=104,
  grid=both,
  minor grid style={opacity=0.3},
  major grid style={opacity=0.3},
  xlabel={$H$ (subdomain side length in wavelengths)},
  ylabel={Iteration Count},
  width=\textwidth,
  height=\textwidth,
  legend pos=north east,
  legend style={fill=none, draw=none, font=\scriptsize, text opacity=1}
]
\addplot+[mark=square*, mark options={solid, fill = blue}, mark size=1.8pt, line width=1pt, color=blue, opacity=1, style=solid] coordinates {(0.29178406,39) (0.34483571,35) (0.42441318,39) (0.5570423,46) (0.82230054,52)};
\addplot+[mark=square*, mark options={solid, fill = orange}, mark size=1.8pt, line width=1pt, color=orange, opacity=1, style=solid] coordinates {(0.29178406,28) (0.34483571,26) (0.42441318,25) (0.5570423,26) (0.82230054,37)};
\addplot+[mark=square*, mark options={solid, fill = cyan}, mark size=1.8pt, line width=1pt, color=cyan, opacity=1, style=solid] coordinates {(0.29178406,18) (0.34483571,21) (0.42441318,22) (0.5570423,23) (0.82230054,36)};
\addplot+[mark=*, mark options={solid, fill = blue}, mark size=1.8pt, line width=1pt, color=blue, opacity=1, style=dashed] coordinates {(0.29178406,38) (0.34483571,34) (0.42441318,38) (0.5570423,38) (0.82230054,51)};
\addplot+[mark=*, mark options={solid, fill = orange}, mark size=1.8pt, line width=1pt, color=orange, opacity=1, style=dashed] coordinates {(0.29178406,26) (0.34483571,26) (0.42441318,25) (0.5570423,24) (0.82230054,35)};
\addplot+[mark=*, mark options={solid, fill = cyan}, mark size=1.8pt, line width=1pt, color=cyan, opacity=1, style=dashed] coordinates {(0.29178406,18) (0.34483571,21) (0.42441318,22) (0.5570423,16) (0.82230054,16)};
\addplot+[mark=none, black, dashed, line width=1.2pt, opacity=1] coordinates {(0.29178406,104) (0.34483571,96) (0.42441318,92) (0.5570423,79) (0.82230054,41)};
\end{axis}
\end{tikzpicture}
    \end{subfigure} 
    \begin{subfigure}[t]{0.48\textwidth}
        \centering
        \begin{tikzpicture}
\begin{axis}[
  xmode=log,
  ymode=log,
  enlarge x limits=0.05,
  enlarge y limits=0.05,
  xmin=1.352817,
  xmax=4.0053994,
  ymin=468,
  ymax=8688,
  grid=both,
  minor grid style={opacity=0.3},
  major grid style={opacity=0.3},
  xlabel={$H$ (subdomain side length in wavelengths)},
  ylabel={CS Size},
  width=\textwidth,
  height=\textwidth,
  legend pos=north east,
  legend style={fill=none, draw=none, font=\small, text opacity=1}
]
\addplot+[mark=*, mark options={solid, fill = blue}, mark size=1.8pt, line width=1pt, color=blue, opacity=1] coordinates {(1.352817,1544) (1.6180753,1404) (2.0159626,1092) (2.6791082,756) (3.2096247,624) (4.0053994,468)};
\addplot+[mark=*, mark options={solid, fill = orange}, mark size=1.8pt, line width=1pt, color=orange, opacity=1] coordinates {(1.352817,3652) (1.6180753,2844) (2.0159626,2212) (2.6791082,1596) (3.2096247,1280) (4.0053994,968)};
\addplot+[mark=*, mark options={solid, fill = cyan}, mark size=1.8pt, line width=1pt, color=cyan, opacity=1] coordinates {(1.352817,7392) (1.6180753,6084) (2.0159626,4768) (2.6791082,3428) (3.2096247,2732) (4.0053994,2052)};
\addplot+[mark=*, mark options={solid, fill = blue}, mark size=1.8pt, line width=1pt, color=blue, opacity=1, style=dashed] coordinates {(4.0053994,616) (3.2096247,728) (2.6791082,928) (2.0159626,1156) (1.6180753,1408) (1.352817,1628)};
\addplot+[mark=*, mark options={solid, fill = orange}, mark size=1.8pt, line width=1pt, color=orange, opacity=1, style=dashed] coordinates {(4.0053994,1412) (3.2096247,1739) (2.6791082,2024) (2.0159626,2560) (1.6180753,3216) (1.352817,3796)};
\addplot+[mark=*, mark options={solid, fill = cyan}, mark size=1.8pt, line width=1pt, color=cyan, opacity=1, style=dashed] coordinates {(4.0053994,3688) (3.2096247,4384) (2.6791082,4936) (2.0159626,6272) (1.6180753,7520) (1.352817,8688)};
\end{axis}
\end{tikzpicture}
    \end{subfigure} 
    \hfill
    \begin{subfigure}[t]{0.48\textwidth}
        \centering
        \begin{tikzpicture}
\begin{axis}[
  xmode=log,
  ymode=log,
  enlarge x limits=0.05,
  enlarge y limits=0.05,
  xmin=1.352817,
  xmax=4.0053994,
  ymin=38,
  ymax=200,
  grid=both,
  minor grid style={opacity=0.3},
  major grid style={opacity=0.3},
  xlabel={$H$ (subdomain side length in wavelengths)},
  ylabel={Iteration Count},
  width=\textwidth,
  height=\textwidth,
  legend pos=north east,
  legend style={fill=none, draw=none, font=\small, text opacity=1},
  legend columns = 4,
  legend to name = Hwaves
]
\addplot+[mark=*, mark options={solid, fill = blue}, mark size=1.8pt, line width=1pt, color=blue, opacity=1] coordinates {(1.352817,200) (1.6180753,200) (2.0159626,200) (2.6791082,200) (4.0053994,200)};
\addlegendentry{$\Delta$-GenEO (0.2)}
\addplot+[mark=*, mark options={solid, fill = orange}, mark size=1.8pt, line width=1pt, color=orange, opacity=1] coordinates {(1.352817,200) (1.6180753,200) (2.0159626,200) (2.6791082,200) (4.0053994,200)};
\addlegendentry{$\Delta$-GenEO (0.4)}
\addplot+[mark=*, mark options={solid, fill = cyan}, mark size=1.8pt, line width=1pt, color=cyan, opacity=1] coordinates {(1.352817,200) (1.6180753,200) (2.0159626,200) (2.6791082,200) (4.0053994,200)};
\addlegendentry{$\Delta$-GenEO (0.7)}
\addplot+[mark=*, mark options={solid, fill = blue}, mark size=1.8pt, line width=1pt, color=blue, opacity=1, style=dashed] coordinates {(1.352817,200) (1.6180753,200) (2.0159626,200) (2.6791082,200) (4.0053994,200)};
\addlegendentry{$\Delta_k$-GenEO (0.2)}
\addplot+[mark=*, mark options={solid, fill = orange}, mark size=1.8pt, line width=1pt, color=orange, opacity=1, style=dashed] coordinates {(1.352817,200) (1.6180753,200) (2.0159626,200) (2.6791082,200) (4.0053994,200)};
\addlegendentry{$\Delta_k$-GenEO (0.4)}
\addplot+[mark=*, mark options={solid, fill = cyan}, mark size=1.8pt, line width=1pt, color=cyan, opacity=1, style=dashed] coordinates {(1.352817,55) (1.6180753,44) (2.0159626,38) (2.6791082,78) (4.0053994,52)};
\addlegendentry{$\Delta_k$-GenEO (0.7)}
\addplot+[mark=none, black, dashed, line width=1.2pt, opacity=1] coordinates {(1.352817,200) (1.6180753,200) (2.0159626,200) (2.6791082,200) (4.0053994,200)};
\addlegendentry{One-level}
\end{axis}
\end{tikzpicture}
    \end{subfigure} 
    \ref{Hwaves}
    \caption{Influence of the subdomain diameter on the iteration count (left) and coarse space size (right) for the homogeneous media test case $k$ = 20 (top) and $k$ = 100 (bottom). The number in brackets indicates $\tau$ used.}
    \label{fig:DeltaHomH}
\end{figure}

\begin{table}[H]
    \centering
    \hspace*{-2cm}
\begin{tabular}{cccccccc|c|ccc|ccc}
\multicolumn{8}{c|}{} & One Level & \multicolumn{3}{c}{$\Delta$-GenEO} & \multicolumn{3}{c}{$\Delta_k$-GenEO} \\
$N$ & $k$ & $h^{-1}$ & $a_\text{max}$ & $n$ & $L$ & $H$ & $n_\text{loc}$ & It & It & CS & CS$_\text{loc}$ & It & CS & CS$_\text{loc}$  \\
\hline
16 & 20 & 240 & 10 & 58081 & 3.18 & 0.82 & 3969 & 134 & 20 & 400 & 25 & 21 & 258 & 16 \\
 & 60 & 720 & 10 & 519841 &  9.55 & 2.41 & 33489 & $>$200 & 121 & 1234 & 77 & 91 & 842 & 53 \\
 & 100 & 1200 & 10 & 1442401 &  15.92 & 4.01 & 91809 & $>$200 & $>$200 & 2062 & 129 & 146 & 1608 & 100 \\
\hline
16 & 20 & 240 & 100 & 58081 & 3.18 & 0.82 & 3969 & 121 & 19 & 400 & 25 & 21 & 262 & 16 \\
 & 60 & 720 & 100 & 519841 &  9.55 & 2.41 & 33489 & $>$200 & 76 & 1234 & 77 & 61 & 826 & 52 \\
 & 100 & 1200 & 100 & 1442401 &  15.92 & 4.01 & 91809 & $>$200 & 196 & 2064 & 129 & 134 & 1566 & 98 \\
\hline
16 & 20 & 240 & 1000 & 58081 & 3.18 & 0.82 & 3969 & 90 & 15 & 400 & 25 & 18 & 262 & 16 \\
 & 60 & 720 & 1000 & 519841 &  9.55 & 2.41 & 33489 & $>$200 & 66 & 1234 & 77 & 58 & 820 & 51 \\
 & 100 & 1200 & 1000 & 1442401 &  15.92 & 4.01 & 91809 & $>$200 & 179 & 2064 & 129 & 111 & 1568 & 98 \\
\hline
36 & 20 & 240 & 10 & 58081 & 3.18 & 0.56 & 1849 & 174 & 19 & 672 & 19 & 20 & 420 & 12 \\
 & 60 & 720 & 10 & 519841 &  9.55 & 1.62 & 15129 & $>$200 & 151 & 2042 & 57 & 125 & 1328 & 37 \\
 & 100 & 1200 & 10 & 1442401 &  15.92 & 2.68 & 41209 & $>$200 & $>$200 & 3426 & 95 & 195 & 2410 & 67 \\
\hline
36 & 20 & 240 & 100 & 58081 & 3.18 & 0.56 & 1849 & 147 & 19 & 666 & 18 & 21 & 422 & 12 \\
 & 60 & 720 & 100 & 519841 &  9.55 & 1.62 & 15129 & $>$200 & 108 & 2050 & 57 & 73 & 1320 & 37 \\
 & 100 & 1200 & 100 & 1442401 &  15.92 & 2.68 & 41209 & $>$200 & $>$200 & 3426 & 95 & 194 & 2356 & 65 \\
\hline
36 & 20 & 240 & 1000 & 58081 & 3.18 & 0.56 & 1849 & 107 & 16 & 666 & 18 & 18 & 422 & 12 \\
 & 60 & 720 & 1000 & 519841 &  9.55 & 1.62 & 15129 & $>$200 & 85 & 2050 & 57 & 67 & 1320 & 37 \\
 & 100 & 1200 & 1000 & 1442401 &  15.92 & 2.68 & 41209 & $>$200 & $>$200 & 3426 & 95 & 176 & 2356 & 65 \\
\hline
64 & 20 & 240 & 10 & 58081 & 3.18 & 0.42 & 1089 & 200 & 18 & 900 & 14 & 20 & 584 & 9 \\
 & 60 & 720 & 10 & 519841 &  9.55 & 1.22 & 8649 & $>$200 & 177 & 2862 & 45 & 129 & 1804 & 28 \\
 & 100 & 1200 & 10 & 1442401 &  15.92 & 2.02 & 23409 & $>$200 & $>$200 & 4778 & 75 & $>$200 & 3200 & 50 \\
\hline
64 & 20 & 240 & 100 & 58081 & 3.18 & 0.42 & 1089 & 168 & 19 & 904 & 14 & 21 & 586 & 9 \\
 & 60 & 720 & 100 & 519841 &  9.55 & 1.22 & 8649 & $>$200 & 130 & 2850 & 45 & 84 & 1794 & 28 \\
 & 100 & 1200 & 100 & 1442401 &  15.92 & 2.02 & 23409 & $>$200 & $>$200 & 4768 & 74 & $>$200 & 3168 & 50 \\
\hline
64 & 20 & 240 & 1000 & 58081 & 3.18 & 0.42 & 1089 & 127 & 16 & 904 & 14 & 18 & 586 & 9 \\
 & 60 & 720 & 1000 & 519841 &  9.55 & 1.22 & 8649 & $>$200 & 119 & 2850 & 45 & 80 & 1794 & 28 \\
 & 100 & 1200 & 1000 & 1442401 &  15.92 & 2.02 & 23409 & $>$200 & $>$200 & 4768 & 74 & $>$200 & 3168 & 50 \\
\hline
100 & 20 & 240 & 10 & 58081 & 3.18 & 0.34 & 729 & 179 & 18 & 1176 & 12 & 19 & 752 & 8 \\
 & 60 & 720 & 10 & 519841 &  9.55 & 0.98 & 5625 & $>$200 & 162 & 3640 & 36 & 112 & 2388 & 24 \\
 & 100 & 1200 & 10 & 1442401 &  15.92 & 1.62 & 15129 & $>$200 & $>$200 & 6134 & 61 & $>$200 & 4070 & 41 \\
\hline
100 & 20 & 240 & 100 & 58081 & 3.18 & 0.34 & 729 & 109 & 19 & 1168 & 12 & 19 & 792 & 8 \\
 & 60 & 720 & 100 & 519841 &  9.55 & 0.98 & 5625 & $>$200 & 94 & 3664 & 37 & 64 & 2390 & 24 \\
 & 100 & 1200 & 100 & 1442401 &  15.92 & 1.62 & 15129 & $>$200 & $>$200 & 6174 & 62 & 171 & 4160 & 42 \\
\hline
100 & 20 & 240 & 1000 & 58081 & 3.18 & 0.34 & 729 & 68 & 16 & 1168 & 12 & 17 & 792 & 8 \\
 & 60 & 720 & 1000 & 519841 &  9.55 & 0.98 & 5625 & $>$200 & 68 & 3664 & 37 & 56 & 2382 & 24 \\
 & 100 & 1200 & 1000 & 1442401 &  15.92 & 1.62 & 15129 & $>$200 & 170 & 6174 & 62 & 113 & 4192 & 42 \\
\hline
144 & 20 & 240 & 10 & 58081 & 3.18 & 0.29 & 529 & $>$200 & 17 & 1436 & 10 & 19 & 938 & 7 \\
 & 60 & 720 & 10 & 519841 &  9.55 & 0.82 & 3969 & $>$200 & 174 & 4424 & 31 & 126 & 2766 & 19 \\
 & 100 & 1200 & 10 & 1442401 &  15.92 & 1.35 & 10609 & $>$200 & $>$200 & 7422 & 52 & $>$200 & 4838 & 34 \\
\hline
144 & 20 & 240 & 100 & 58081 & 3.18 & 0.29 & 529 & 198 & 19 & 1412 & 10 & 20 & 934 & 6 \\
 & 60 & 720 & 100 & 519841 &  9.55 & 0.82 & 3969 & $>$200 & 135 & 4428 & 31 & 89 & 2772 & 19 \\
 & 100 & 1200 & 100 & 1442401 &  15.92 & 1.35 & 10609 & $>$200 & $>$200 & 7466 & 52 & $>$200 & 4834 & 34 \\
\hline
144 & 20 & 240 & 1000 & 58081 & 3.18 & 0.29 & 529 & 130 & 16 & 1412 & 10 & 18 & 934 & 6 \\
 & 60 & 720 & 1000 & 519841 &  9.55 & 0.82 & 3969 & $>$200 & 130 & 4428 & 31 & 90 & 2772 & 19 \\
 & 100 & 1200 & 1000 & 1442401 &  15.92 & 1.35 & 10609 & $>$200 & $>$200 & 7466 & 52 & $>$200 & 4810 & 33 \\
\hline
\end{tabular}
    \caption{Results showing the dimension of the fine mesh ($n$), the diameter of the domain measured in wavelengths ($L$), the diameter of the sub-domains measured in wavelengths ($H$), the number of local degrees of freedom ($n_\text{loc}$), the GMRES iteration count (It.) for the one-level and two level methods using the $\Delta$-GenEO and $\Delta_k$-GenEO coarse space, the dimension of the coarse space (CS) and averaged number of contributions to the coarse space per subdomain (CS$_\text{loc}$) for the $\Delta$-GenEO and $\Delta_k$-GenEO. These results are for the heterogeneous media test case. Using $\tau = 0.7$ and $\tau=0.5$ for $\Delta$-GenEO and $\Delta_k$-GenEO respectively.}
    \label{table:DeltaHet}
\end{table}

\begin{figure}[H]
    \centering
    \begin{subfigure}[t]{0.48\textwidth}
        \centering
        \begin{tikzpicture}
\begin{axis}[
  xmode=log,
  ymode=log,
  enlarge x limits=0.05,
  enlarge y limits=0.05,
  xmin=40,
  xmax=1456,
  ymin=17,
  ymax=200,
  grid=both,
  minor grid style={opacity=0.3},
  major grid style={opacity=0.3},
  xlabel={CS Size},
  ylabel={Iteration Count},
  width=\textwidth,
  height=\textwidth,
  legend pos=north east,
  legend style={fill=none, draw=none, font=\scriptsize, text opacity=1}
]
\addplot+[mark=square*, mark options={solid, fill = blue}, mark size=1.8pt, line width=1pt, color=blue, opacity=1, style=solid] coordinates {(40,49) (82,33) (146,25) (188,23) (254,21) (312,21) (400,20)};
\addplot+[mark=square*, mark options={solid, fill = orange}, mark size=1.8pt, line width=1pt, color=orange, opacity=1, style=solid] coordinates {(82,46) (126,36) (232,25) (308,22) (420,20) (514,19) (672,19)};
\addplot+[mark=square*, mark options={solid, fill = green!60!black}, mark size=1.8pt, line width=1pt, color=green!60!black, opacity=1, style=solid] coordinates {(114,70) (192,35) (326,26) (430,22) (584,20) (704,19) (900,18)};
\addplot+[mark=square*, mark options={solid, fill = red!80!black}, mark size=1.8pt, line width=1pt, color=red!80!black, opacity=1, style=solid] coordinates {(180,66) (310,33) (432,28) (564,21) (750,19) (932,19) (1176,18)};
\addplot+[mark=square*, mark options={solid, fill = cyan}, mark size=1.8pt, line width=1pt, color=cyan, opacity=1, style=solid] coordinates {(174,76) (384,32) (466,27) (588,25) (938,19) (1090,18) (1436,17)};
\addplot+[mark=*, mark options={solid, fill = blue}, mark size=1.8pt, line width=1pt, color=blue, opacity=1, style=dashed] coordinates {(40,49) (84,32) (146,25) (188,23) (258,21) (314,21) (408,20)};
\addplot+[mark=*, mark options={solid, fill = orange}, mark size=1.8pt, line width=1pt, color=orange, opacity=1, style=dashed] coordinates {(82,46) (128,36) (232,24) (316,21) (420,20) (518,19) (682,19)};
\addplot+[mark=*, mark options={solid, fill = green!60!black}, mark size=1.8pt, line width=1pt, color=green!60!black, opacity=1, style=dashed] coordinates {(130,58) (192,35) (326,26) (436,22) (584,20) (720,19) (912,18)};
\addplot+[mark=*, mark options={solid, fill = red!80!black}, mark size=1.8pt, line width=1pt, color=red!80!black, opacity=1, style=dashed] coordinates {(180,66) (310,33) (464,25) (564,21) (752,19) (932,18) (1178,18)};
\addplot+[mark=*, mark options={solid, fill = cyan}, mark size=1.8pt, line width=1pt, color=cyan, opacity=1, style=dashed] coordinates {(174,76) (384,32) (466,27) (614,23) (938,19) (1102,18) (1456,17)};
\addplot+[mark=none, black, dashed, line width=1.2pt, opacity=0.2] coordinates {(40,134) (1456,134)};
\addplot+[mark=none, black, dashed, line width=1.2pt, opacity=0.4] coordinates {(40,174) (1456,174)};
\addplot+[mark=none, black, dashed, line width=1.2pt, opacity=0.6] coordinates {(40,200) (1456,200)};
\addplot+[mark=none, black, dashed, line width=1.2pt, opacity=0.8] coordinates {(40,179) (1456,179)};
\addplot+[mark=none, black, dashed, line width=1.2pt, opacity=1.0] coordinates {(40,200) (1456,200)};
\end{axis}
\end{tikzpicture}
    \end{subfigure}   
    \hfill
    \begin{subfigure}[t]{0.48\textwidth}
        \centering
        \begin{tikzpicture}
\begin{axis}[
  xmode=log,
  ymode=log,
  enlarge x limits=0.05,
  enlarge y limits=0.05,
  xmin=0.1,
  xmax=0.7,
  ymin=17,
  ymax=200,
  grid=both,
  minor grid style={opacity=0.3},
  major grid style={opacity=0.3},
  xlabel={$\tau$},
  ylabel={Iteration Count},
  width=\textwidth,
  height=\textwidth,
  legend pos=north east,
  legend style={fill=none, draw=none, font=\scriptsize, text opacity=1}
]
\addplot+[mark=square*, mark options={solid, fill = blue}, mark size=1.8pt, line width=1pt, color=blue, opacity=1, style=solid] coordinates {(0.1,49) (0.2,33) (0.3,25) (0.4,23) (0.5,21) (0.6,21) (0.7,20)};
\addplot+[mark=square*, mark options={solid, fill = orange}, mark size=1.8pt, line width=1pt, color=orange, opacity=1, style=solid] coordinates {(0.1,46) (0.2,36) (0.3,25) (0.4,22) (0.5,20) (0.6,19) (0.7,19)};
\addplot+[mark=square*, mark options={solid, fill = green!60!black}, mark size=1.8pt, line width=1pt, color=green!60!black, opacity=1, style=solid] coordinates {(0.1,70) (0.2,35) (0.3,26) (0.4,22) (0.5,20) (0.6,19) (0.7,18)};
\addplot+[mark=square*, mark options={solid, fill = red!80!black}, mark size=1.8pt, line width=1pt, color=red!80!black, opacity=1, style=solid] coordinates {(0.1,66) (0.2,33) (0.3,28) (0.4,21) (0.5,19) (0.6,19) (0.7,18)};
\addplot+[mark=square*, mark options={solid, fill = cyan}, mark size=1.8pt, line width=1pt, color=cyan, opacity=1, style=solid] coordinates {(0.1,76) (0.2,32) (0.3,27) (0.4,25) (0.5,19) (0.6,18) (0.7,17)};
\addplot+[mark=*, mark options={solid, fill = blue}, mark size=1.8pt, line width=1pt, color=blue, opacity=1, style=dashed] coordinates {(0.1,49) (0.2,32) (0.3,25) (0.4,23) (0.5,21) (0.6,21) (0.7,20)};
\addplot+[mark=*, mark options={solid, fill = orange}, mark size=1.8pt, line width=1pt, color=orange, opacity=1, style=dashed] coordinates {(0.1,46) (0.2,36) (0.3,24) (0.4,21) (0.5,20) (0.6,19) (0.7,19)};
\addplot+[mark=*, mark options={solid, fill = green!60!black}, mark size=1.8pt, line width=1pt, color=green!60!black, opacity=1, style=dashed] coordinates {(0.1,58) (0.2,35) (0.3,26) (0.4,22) (0.5,20) (0.6,19) (0.7,18)};
\addplot+[mark=*, mark options={solid, fill = red!80!black}, mark size=1.8pt, line width=1pt, color=red!80!black, opacity=1, style=dashed] coordinates {(0.1,66) (0.2,33) (0.3,25) (0.4,21) (0.5,19) (0.6,18) (0.7,18)};
\addplot+[mark=*, mark options={solid, fill = cyan}, mark size=1.8pt, line width=1pt, color=cyan, opacity=1, style=dashed] coordinates {(0.1,76) (0.2,32) (0.3,27) (0.4,23) (0.5,19) (0.6,18) (0.7,17)};
\addplot+[mark=none, black, dashed, line width=1.2pt, opacity=0.2] coordinates {(0.1,134) (0.7,134)};
\addplot+[mark=none, black, dashed, line width=1.2pt, opacity=0.4] coordinates {(0.1,174) (0.7,174)};
\addplot+[mark=none, black, dashed, line width=1.2pt, opacity=0.6] coordinates {(0.1,200) (0.7,200)};
\addplot+[mark=none, black, dashed, line width=1.2pt, opacity=0.8] coordinates {(0.1,179) (0.7,179)};
\addplot+[mark=none, black, dashed, line width=1.2pt, opacity=1.0] coordinates {(0.1,200) (0.7,200)};
\end{axis}
\end{tikzpicture}
    \end{subfigure} 
    \begin{subfigure}[t]{0.48\textwidth}
        \centering
        \begin{tikzpicture}
\begin{axis}[
  xmode=log,
  ymode=log,
  enlarge x limits=0.05,
  enlarge y limits=0.05,
  xmin=246,
  xmax=8078,
  ymin=38,
  ymax=200,
  grid=both,
  minor grid style={opacity=0.3},
  major grid style={opacity=0.3},
  xlabel={CS Size},
  ylabel={Iteration Count},
  width=\textwidth,
  height=\textwidth,
  legend pos=north east,
  legend style={fill=none, draw=none, font=\small, text opacity=1},
  legend columns = 5,
  legend to name = HetSub
]
\addplot+[mark=square*, mark options={solid, fill = blue}, mark size=1.8pt, line width=1pt, color=blue, opacity=1, style=solid] coordinates {(246,200) (476,200) (726,200) (996,200) (1274,200) (1618,200) (2062,200)};
\addlegendentry{$\Delta$-GenEO (16)}
\addplot+[mark=square*, mark options={solid, fill = orange}, mark size=1.8pt, line width=1pt, color=orange, opacity=1, style=solid] coordinates {(376,200) (764,200) (1188,200) (1638,200) (2122,200) (2690,200) (3426,200)};
\addlegendentry{$\Delta$-GenEO (36)}
\addplot+[mark=square*, mark options={solid, fill = green!60!black}, mark size=1.8pt, line width=1pt, color=green!60!black, opacity=1, style=solid] coordinates {(540,200) (1108,200) (1644,200) (2250,200) (2928,200) (3764,200) (4778,200)};
\addlegendentry{$\Delta$-GenEO (64)}
\addplot+[mark=square*, mark options={solid, fill = red!80!black}, mark size=1.8pt, line width=1pt, color=red!80!black, opacity=1, style=solid] coordinates {(800,200) (1594,200) (2250,200) (3018,200) (3860,200) (4860,200) (6134,200)};
\addlegendentry{$\Delta$-GenEO (100)}
\addplot+[mark=square*, mark options={solid, fill = cyan}, mark size=1.8pt, line width=1pt, color=cyan, opacity=1, style=solid] coordinates {(934,200) (1576,200) (2620,200) (3592,200) (4670,200) (5802,200) (7422,200)};
\addlegendentry{$\Delta$-GenEO (144)}
\addplot+[mark=*, mark options={solid, fill = blue}, mark size=1.8pt, line width=1pt, color=blue, opacity=1, style=dashed] coordinates {(274,200) (552,200) (846,190) (1188,166) (1608,146) (2138,65) (2890,38)};
\addlegendentry{$\Delta_k$-GenEO (16)}
\addplot+[mark=*, mark options={solid, fill = orange}, mark size=1.8pt, line width=1pt, color=orange, opacity=1, style=dashed] coordinates {(420,200) (840,200) (1314,200) (1806,200) (2410,195) (3164,80) (4216,38)};
\addlegendentry{$\Delta_k$-GenEO (36)}
\addplot+[mark=*, mark options={solid, fill = green!60!black}, mark size=1.8pt, line width=1pt, color=green!60!black, opacity=1, style=dashed] coordinates {(596,200) (1148,200) (1758,200) (2400,200) (3200,200) (4190,126) (5512,60)};
\addlegendentry{$\Delta_k$-GenEO (64)}
\addplot+[mark=*, mark options={solid, fill = red!80!black}, mark size=1.8pt, line width=1pt, color=red!80!black, opacity=1, style=dashed] coordinates {(892,200) (1700,200) (2382,200) (3210,200) (4070,200) (5292,76) (6874,45)};
\addlegendentry{$\Delta_k$-GenEO (100)}
\addplot+[mark=*, mark options={solid, fill = cyan}, mark size=1.8pt, line width=1pt, color=cyan, opacity=1, style=dashed] coordinates {(938,200) (1678,200) (2648,200) (3708,200) (4838,200) (6166,185) (8078,82)};
\addlegendentry{$\Delta_k$-GenEO (144)}
\addplot+[mark=none, black, dashed, line width=1.2pt, opacity=0.2] coordinates {(246,200) (8078,200)};
\addlegendentry{One-level (16)}
\addplot+[mark=none, black, dashed, line width=1.2pt, opacity=0.4] coordinates {(246,200) (8078,200)};
\addlegendentry{One-level (36)}
\addplot+[mark=none, black, dashed, line width=1.2pt, opacity=0.6] coordinates {(246,200) (8078,200)};
\addlegendentry{One-level (64)}
\addplot+[mark=none, black, dashed, line width=1.2pt, opacity=0.8] coordinates {(246,200) (8078,200)};
\addlegendentry{One-level (100)}
\addplot+[mark=none, black, dashed, line width=1.2pt, opacity=1.0] coordinates {(246,200) (8078,200)};
\addlegendentry{One-level (144)}
\end{axis}
\end{tikzpicture}
    \end{subfigure} 
    \hfill
    \begin{subfigure}[t]{0.48\textwidth}
        \centering
        \begin{tikzpicture}
\begin{axis}[
  xmode=log,
  ymode=log,
  enlarge x limits=0.05,
  enlarge y limits=0.05,
  xmin=0.1,
  xmax=0.7,
  ymin=38,
  ymax=200,
  grid=both,
  minor grid style={opacity=0.3},
  major grid style={opacity=0.3},
  xlabel={$\tau$},
  ylabel={Iteration Count},
  width=\textwidth,
  height=\textwidth,
  legend pos=north east,
  legend style={fill=none, draw=none, font=\scriptsize, text opacity=1}
]
\addplot+[mark=square*, mark options={solid, fill = blue}, mark size=1.8pt, line width=1pt, color=blue, opacity=1, style=solid] coordinates {(0.1,200) (0.2,200) (0.3,200) (0.4,200) (0.5,200) (0.6,200) (0.7,200)};
\addplot+[mark=square*, mark options={solid, fill = orange}, mark size=1.8pt, line width=1pt, color=orange, opacity=1, style=solid] coordinates {(0.1,200) (0.2,200) (0.3,200) (0.4,200) (0.5,200) (0.6,200) (0.7,200)};
\addplot+[mark=square*, mark options={solid, fill = green!60!black}, mark size=1.8pt, line width=1pt, color=green!60!black, opacity=1, style=solid] coordinates {(0.1,200) (0.2,200) (0.3,200) (0.4,200) (0.5,200) (0.6,200) (0.7,200)};
\addplot+[mark=square*, mark options={solid, fill = red!80!black}, mark size=1.8pt, line width=1pt, color=red!80!black, opacity=1, style=solid] coordinates {(0.1,200) (0.2,200) (0.3,200) (0.4,200) (0.5,200) (0.6,200) (0.7,200)};
\addplot+[mark=square*, mark options={solid, fill = cyan}, mark size=1.8pt, line width=1pt, color=cyan, opacity=1, style=solid] coordinates {(0.1,200) (0.2,200) (0.3,200) (0.4,200) (0.5,200) (0.6,200) (0.7,200)};
\addplot+[mark=*, mark options={solid, fill = blue}, mark size=1.8pt, line width=1pt, color=blue, opacity=1, style=dashed] coordinates {(0.1,200) (0.2,200) (0.3,190) (0.4,166) (0.5,146) (0.6,65) (0.7,38)};
\addplot+[mark=*, mark options={solid, fill = orange}, mark size=1.8pt, line width=1pt, color=orange, opacity=1, style=dashed] coordinates {(0.1,200) (0.2,200) (0.3,200) (0.4,200) (0.5,195) (0.6,80) (0.7,38)};
\addplot+[mark=*, mark options={solid, fill = green!60!black}, mark size=1.8pt, line width=1pt, color=green!60!black, opacity=1, style=dashed] coordinates {(0.1,200) (0.2,200) (0.3,200) (0.4,200) (0.5,200) (0.6,126) (0.7,60)};
\addplot+[mark=*, mark options={solid, fill = red!80!black}, mark size=1.8pt, line width=1pt, color=red!80!black, opacity=1, style=dashed] coordinates {(0.1,200) (0.2,200) (0.3,200) (0.4,200) (0.5,200) (0.6,76) (0.7,45)};
\addplot+[mark=*, mark options={solid, fill = cyan}, mark size=1.8pt, line width=1pt, color=cyan, opacity=1, style=dashed] coordinates {(0.1,200) (0.2,200) (0.3,200) (0.4,200) (0.5,200) (0.6,185) (0.7,82)};
\addplot+[mark=none, black, dashed, line width=1.2pt, opacity=0.2] coordinates {(0.1,200) (0.7,200)};
\addplot+[mark=none, black, dashed, line width=1.2pt, opacity=0.4] coordinates {(0.1,200) (0.7,200)};
\addplot+[mark=none, black, dashed, line width=1.2pt, opacity=0.6] coordinates {(0.1,200) (0.7,200)};
\addplot+[mark=none, black, dashed, line width=1.2pt, opacity=0.8] coordinates {(0.1,200) (0.7,200)};
\addplot+[mark=none, black, dashed, line width=1.2pt, opacity=1.0] coordinates {(0.1,200) (0.7,200)};
\end{axis}
\end{tikzpicture}
    \end{subfigure} 
    \ref{HetSub}
    \caption{Influence of the coarse space size (left) and threshold choice (right) on the iteration count for the heterogeneous media test case with $k$ = 20 (top) and $k$ = 100 (bottom), all with $a_\text{max}(\bx) = 10$. The number in brackets indicates the number of subdomains.}
    \label{fig:DeltaHetCSIT}
\end{figure}

\begin{figure}[H]
    \centering
    \begin{subfigure}[t]{0.48\textwidth}
        \centering
        \begin{tikzpicture}
\begin{axis}[
  xmode=log,
  ymode=log,
  enlarge x limits=0.05,
  enlarge y limits=0.05,
  xmin=0.29178406,
  xmax=0.82230054,
  ymin=82,
  ymax=1456,
  grid=both,
  minor grid style={opacity=0.3},
  major grid style={opacity=0.3},
  xlabel={$H$ (subdomain side length in wavelengths)},
  ylabel={CS Size},
  width=\textwidth,
  height=\textwidth,
  legend pos=north east,
  legend style={fill=none, draw=none, font=\scriptsize, text opacity=1}
]
\addplot+[mark=square*, mark options={solid, fill = blue}, mark size=1.8pt, line width=1pt, color=blue, opacity=1, style=solid] coordinates {(0.29178406,384) (0.34483571,310) (0.42441318,192) (0.5570423,126) (0.82230054,82)};
\addplot+[mark=square*, mark options={solid, fill = orange}, mark size=1.8pt, line width=1pt, color=orange, opacity=1, style=solid] coordinates {(0.29178406,588) (0.34483571,564) (0.42441318,430) (0.5570423,308) (0.82230054,188)};
\addplot+[mark=square*, mark options={solid, fill = cyan}, mark size=1.8pt, line width=1pt, color=cyan, opacity=1, style=solid] coordinates {(0.29178406,1436) (0.34483571,1176) (0.42441318,900) (0.5570423,672) (0.82230054,400)};
\addplot+[mark=*, mark options={solid, fill = blue}, mark size=1.8pt, line width=1pt, color=blue, opacity=1, style=dashed] coordinates {(0.82230054,84) (0.5570423,128) (0.42441318,192) (0.34483571,310) (0.29178406,384)};
\addplot+[mark=*, mark options={solid, fill = orange}, mark size=1.8pt, line width=1pt, color=orange, opacity=1, style=dashed] coordinates {(0.82230054,188) (0.5570423,316) (0.42441318,436) (0.34483571,564) (0.29178406,614)};
\addplot+[mark=*, mark options={solid, fill = cyan}, mark size=1.8pt, line width=1pt, color=cyan, opacity=1, style=dashed] coordinates {(0.82230054,408) (0.5570423,682) (0.42441318,912) (0.34483571,1178) (0.29178406,1456)};
\end{axis}
\end{tikzpicture}
    \end{subfigure}   
    \hfill
    \begin{subfigure}[t]{0.48\textwidth}
        \centering
        \begin{tikzpicture}
\begin{axis}[
  xmode=log,
  ymode=log,
  enlarge x limits=0.05,
  enlarge y limits=0.05,
  xmin=0.29178406,
  xmax=0.82230054,
  ymin=17,
  ymax=200,
  grid=both,
  minor grid style={opacity=0.3},
  major grid style={opacity=0.3},
  xlabel={$H$ (subdomain side length in wavelengths)},
  ylabel={Iteration Count},
  width=\textwidth,
  height=\textwidth,
  legend pos=north east,
  legend style={fill=none, draw=none, font=\small, text opacity=1},
  legend columns = 4,
  legend to name = HetHwave
]
\addplot+[mark=square*, mark options={solid, fill = blue}, mark size=1.8pt, line width=1pt, color=blue, opacity=1, style=solid] coordinates {(0.29178406,32) (0.34483571,33) (0.42441318,35) (0.5570423,36) (0.82230054,33)};
\addlegendentry{$\Delta$-GenEO (0.2)}
\addplot+[mark=square*, mark options={solid, fill = orange}, mark size=1.8pt, line width=1pt, color=orange, opacity=1, style=solid] coordinates {(0.29178406,25) (0.34483571,21) (0.42441318,22) (0.5570423,22) (0.82230054,23)};
\addlegendentry{$\Delta$-GenEO (0.4)}
\addplot+[mark=square*, mark options={solid, fill = cyan}, mark size=1.8pt, line width=1pt, color=cyan, opacity=1, style=solid] coordinates {(0.29178406,17) (0.34483571,18) (0.42441318,18) (0.5570423,19) (0.82230054,20)};
\addlegendentry{$\Delta$-GenEO (0.7)}
\addplot+[mark=*, mark options={solid, fill = blue}, mark size=1.8pt, line width=1pt, color=blue, opacity=1, style=dashed] coordinates {(0.29178406,32) (0.34483571,33) (0.42441318,35) (0.5570423,36) (0.82230054,32)};
\addlegendentry{$\Delta_k$-GenEO (0.2)}
\addplot+[mark=*, mark options={solid, fill = orange}, mark size=1.8pt, line width=1pt, color=orange, opacity=1, style=dashed] coordinates {(0.29178406,23) (0.34483571,21) (0.42441318,22) (0.5570423,21) (0.82230054,23)};
\addlegendentry{$\Delta_k$-GenEO (0.4)}
\addplot+[mark=*, mark options={solid, fill = cyan}, mark size=1.8pt, line width=1pt, color=cyan, opacity=1, style=dashed] coordinates {(0.29178406,17) (0.34483571,18) (0.42441318,18) (0.5570423,19) (0.82230054,20)};
\addlegendentry{$\Delta_k$-GenEO (0.7)}
\addplot+[mark=none, black, dashed, line width=1.2pt, opacity=1] coordinates {(0.27852115,200) (0.3315728,179) (0.41115027,200) (0.54377939,174) (0.80903763,134)};
\addlegendentry{One-level}
\end{axis}
\end{tikzpicture}
    \end{subfigure} 
    \begin{subfigure}[t]{0.48\textwidth}
        \centering
        \begin{tikzpicture}
\begin{axis}[
  xmode=log,
  ymode=log,
  enlarge x limits=0.05,
  enlarge y limits=0.05,
  xmin=1.352817,
  xmax=4.0053994,
  ymin=476,
  ymax=8078,
  grid=both,
  minor grid style={opacity=0.3},
  major grid style={opacity=0.3},
  xlabel={$H$ (subdomain side length in wavelengths)},
  ylabel={CS Size},
  width=\textwidth,
  height=\textwidth,
  legend pos=north east,
  legend style={fill=none, draw=none, font=\scriptsize, text opacity=1}
]
\addplot+[mark=square*, mark options={solid, fill = blue}, mark size=1.8pt, line width=1pt, color=blue, opacity=1, style=solid] coordinates {(1.352817,1576) (1.6180753,1594) (2.0159626,1108) (2.6791082,764) (4.0053994,476)};
\addplot+[mark=square*, mark options={solid, fill = orange}, mark size=1.8pt, line width=1pt, color=orange, opacity=1, style=solid] coordinates {(1.352817,3592) (1.6180753,3018) (2.0159626,2250) (2.6791082,1638) (4.0053994,996)};
\addplot+[mark=square*, mark options={solid, fill = cyan}, mark size=1.8pt, line width=1pt, color=cyan, opacity=1, style=solid] coordinates {(1.352817,7422) (1.6180753,6134) (2.0159626,4778) (2.6791082,3426) (4.0053994,2062)};
\addplot+[mark=*, mark options={solid, fill = blue}, mark size=1.8pt, line width=1pt, color=blue, opacity=1, style=dashed] coordinates {(1.352817,1678) (1.6180753,1700) (2.0159626,1148) (2.6791082,840) (4.0053994,552)};
\addplot+[mark=*, mark options={solid, fill = orange}, mark size=1.8pt, line width=1pt, color=orange, opacity=1, style=dashed] coordinates {(1.352817,3708) (1.6180753,3210) (2.0159626,2400) (2.6791082,1806) (4.0053994,1188)};
\addplot+[mark=*, mark options={solid, fill = cyan}, mark size=1.8pt, line width=1pt, color=cyan, opacity=1, style=dashed] coordinates {(1.352817,8078) (1.6180753,6874) (2.0159626,5512) (2.6791082,4216) (4.0053994,2890)};
\end{axis}
\end{tikzpicture}
    \end{subfigure} 
    \hfill
    \begin{subfigure}[t]{0.48\textwidth}
        \centering
        \begin{tikzpicture}
\begin{axis}[
  xmode=log,
  ymode=log,
  enlarge x limits=0.05,
  enlarge y limits=0.05,
  xmin=1.352817,
  xmax=4.0053994,
  ymin=38,
  ymax=200,
  grid=both,
  minor grid style={opacity=0.3},
  major grid style={opacity=0.3},
  xlabel={$H$ (subdomain side length in wavelengths)},
  ylabel={Iteration Count},
  width=\textwidth,
  height=\textwidth,
  legend pos=north east,
  legend style={fill=none, draw=none, font=\scriptsize, text opacity=1}
]
\addplot+[mark=square*, mark options={solid, fill = blue}, mark size=1.8pt, line width=1pt, color=blue, opacity=1, style=solid] coordinates {(1.352817,200) (1.6180753,200) (2.0159626,200) (2.6791082,200) (4.0053994,200)};
\addplot+[mark=square*, mark options={solid, fill = orange}, mark size=1.8pt, line width=1pt, color=orange, opacity=1, style=solid] coordinates {(1.352817,200) (1.6180753,200) (2.0159626,200) (2.6791082,200) (4.0053994,200)};
\addplot+[mark=square*, mark options={solid, fill = cyan}, mark size=1.8pt, line width=1pt, color=cyan, opacity=1, style=solid] coordinates {(1.352817,200) (1.6180753,200) (2.0159626,200) (2.6791082,200) (4.0053994,200)};
\addplot+[mark=*, mark options={solid, fill = blue}, mark size=1.8pt, line width=1pt, color=blue, opacity=1, style=dashed] coordinates {(1.352817,200) (1.6180753,200) (2.0159626,200) (2.6791082,200) (4.0053994,200)};
\addplot+[mark=*, mark options={solid, fill = orange}, mark size=1.8pt, line width=1pt, color=orange, opacity=1, style=dashed] coordinates {(1.352817,200) (1.6180753,200) (2.0159626,200) (2.6791082,200) (4.0053994,166)};
\addplot+[mark=*, mark options={solid, fill = cyan}, mark size=1.8pt, line width=1pt, color=cyan, opacity=1, style=dashed] coordinates {(1.352817,82) (1.6180753,45) (2.0159626,60) (2.6791082,38) (4.0053994,38)};
\addplot+[mark=none, black, dashed, line width=1.2pt, opacity=1] coordinates {(1.3594485,200) (1.6247067,200) (2.0225941,200) (2.6857397,200) (4.0120309,200)};
\end{axis}
\end{tikzpicture}
    \end{subfigure} 
    \ref{HetHwave}
    \caption{Influence of the subdomain diameter on the iteration count (left) and coarse space size (right) for the heterogeneous media test case $k$ = 20 (top) and $k$ = 100 (bottom) , all with $a_\text{max}(\bx) = 10$. The number in brackets indicates $\tau$ used.}
    \label{fig:DeltaHetH}
\end{figure}

\subsection{Heterogenous problem}

We now consider the same experiments as in the homogenous case, but with a heterogeneous diffusion coefficient $A(\bx)$.
This aims to study the robustness of the method with respect to medium heterogeneity. 
For this purpose, we consider a layered medium within the same unit square domain.  More precisely, we use an alternating layer configuration 
in which the heterogeneity is introduced through the material parameter $a(\bx)$ where
\[
	a(\bx) = \begin{cases}
		\amax & \text{if} \ y\in[0,0.1)\cup[0.2,0.3)\cup[0.4,0.5)\cup[0.6,0.7)\cup[0.8,0.9),
		\\
		\amin & \text{otherwise},
	\end{cases}
\]
with $\amax>1$ and $\amin=1$; see Figure~\ref{fig:HetPlot}. 
The heterogeneous coefficient is given by
\begin{equation*}
    A(\bx) = a(\bx)I.
\end{equation*}

\begin{figure}[H]
\label{Fig:AlternatingLay}
\centering
\begin{tikzpicture}[scale=3]
    \fill[gray,opacity=0.1] (-1,1) -- (1,1) -- (1,0.8) -- (-1,0.8) -- cycle;
	\fill[gray,opacity=1.0] (-1,0.6) -- (1,0.6) -- (1,0.8) -- (-1,0.8) -- cycle;
	\fill[gray,opacity=0.1] (-1,0.6) -- (1,0.6) -- (1,0.4) -- (-1,0.4) -- cycle;
	\fill[gray,opacity=1.0] (-1,0.2) -- (1,0.2) -- (1,0.4) -- (-1,0.4) -- cycle;
	\fill[gray,opacity=0.1] (-1,0.2) -- (1,0.2) -- (1,0.0) -- (-1,0.0) -- cycle;
	\fill[gray,opacity=1.0] (-1,-0.2) -- (1,-0.2) -- (1,0.0) -- (-1,0.0) -- cycle;
	\fill[gray,opacity=0.1] (-1,-0.2) -- (1,-0.2) -- (1,-0.4) -- (-1,-0.4) -- cycle;
	\fill[gray,opacity=1.0] (-1,-0.6) -- (1,-0.6) -- (1,-0.4) -- (-1,-0.4) -- cycle;
	\fill[gray,opacity=0.1] (-1,-0.6) -- (1,-0.6) -- (1,-0.8) -- (-1,-0.8) -- cycle;
	\fill[gray,opacity=1.0] (-1,-1) -- (1,-1) -- (1,-0.8) -- (-1,-0.8) -- cycle;
    
    \draw[step=0.5cm,gray!50,very thin, shift={(-1,-1)}] (0,0) grid (2,2);
    \draw[gray!50,very thin] (-1,-1) -- (1,1);
    \draw[gray!50,very thin] (-1,1) -- (1,-1);
    \draw[gray!50,very thin] (-1,0) -- (0,1) -- (1,0) -- (0,-1) -- cycle;
    
    \draw (-1,1) -- (-1,-1) -- (1,-1) -- (1,1) -- cycle;
    
    \draw (1.05,0) node[above,rotate=-90] {Robin};
    \draw (-1.05,0) node[above,rotate=90] {Robin};
    \draw (0,-1.05) node[below] {Robin};
    \draw (0,1.05) node[above] {Robin};
    
    \draw (-1.05,1) node[left] {$y = 1$}
          (-1.05,-1) node[left] {$y = 0$}
          (1,-1.05) node[below] {$x = 1$}
          (-1,-1.05) node[below] {$x = 0$};

    \draw[fill=black] (0,0) circle (0.025) node[below] {source $f$};

\end{tikzpicture}
\caption{The heterogeneous function $a(\bx)$ within the alternating layers. The shading gives the value of $a(\bx)$ 
with the darkest shade being $a(\bx) = \amax$, where $\amax > 1$ is a parameter, and the white taking the value $\amin = 1$.}
\label{fig:HetPlot}
\end{figure}

\subsection{Results}
In Tables~\ref{table:DeltaHomog} and \ref{table:DeltaHet}, we display the results comparing the one-level AS and the two-level AS 
with the $\Delta$-GenEO from~\cite{Bootland:2022:OSM} and the newly proposed $\Delta_k$-GenEO. 
The indefiniteness of the problem is controlled by changing $k$ with different subdomain decompositions $N$.
The entries are the number of global degrees of freedom, $n$, the diameter of the domain, $L$, 
measured in the smallest wavelength $\lambda_\text{min}$, the diameter of the subdomains, $H$, measured in the smallest wavelength, 
the number of local degrees of freedom, $n_\text{loc}$, the iterations required to reach the GMRES tolerance, It, for the one-level, 
$\Delta$-GenEO and $\Delta_k$-GenEO method, the total size of the $\Delta$-GenEO and $\Delta_k$-GenEO coarse space, CS, 
and the averaged number of contributions to the coarse space per subdomain, $CS_\text{loc}$. 
For the $\Delta$-GenEO and $\Delta_k$-GenEO method $\tau$ has been selected for each to give the nearest iteration counts across the two methods. 
In Table~\ref{table:DeltaHet} we also show the effect from changing the $\amax$ being used. 

Further to this, using the smallest and largest values of $k$ tested, we show the number of GMRES iterations required 
to reach the desired tolerance as a function of the total coarse space size and the threshold value selection in 
Figures~\ref{fig:DeltaHomCSIT} and \ref{fig:DeltaHetCSIT}. 
The number of GMRES iterations and coarse space size as a function of the subdomain diameter in wavelengths is 
shown in Figures~\ref{fig:DeltaHomH} and \ref{fig:DeltaHetH}. 
With the heterogeneous results in Figures~\ref{fig:DeltaHetCSIT} and \ref{fig:DeltaHetH} having $\amax=10$.

\subsection{Analyses and comparison}

We now compare the results for the homogenous and heterogeneous test cases across the different methods. 
{In this study, we deliberately focus on iteration counts and coarse-space dimensions rather than reporting CPU timings. 
While wall-clock times can be informative, they are highly sensitive to implementation details, solver tolerances, hardware architecture, 
and the degree of optimisation of the software stack, and therefore do not provide a robust basis for comparison across platforms.}

{For two-level domain decomposition methods, iteration counts and coarse-space size are widely accepted as reliable proxies for performance. 
In particular, reduced iteration counts directly decrease the number of local and coarse solves, while smaller coarse spaces lead to 
cheaper coarse-level factorisations and solves, as well as reduced communication overhead. 
Consequently, the consistently lower iteration counts and more compact coarse spaces observed for the $\Delta_k$-GenEO method translate 
into a lower time-to-solution in practical, optimised implementations, even though explicit timings are not reported here.
}

From both the Tables~\ref{table:DeltaHomog} and \ref{table:DeltaHet}, and Figures~\ref{fig:DeltaHomCSIT} and \ref{fig:DeltaHetCSIT} 
we can identify several points of interest:

\paragraph{Two-level is essential. $\Delta_k$-GenEO dominates}
Across all tested $N$ and $k$, the one-level method quickly becomes non-viable beyond small $k$ or a large $N$. As demonstrated when $k=40$, 
the iteration count exceeds the 200 iteration limit for $N \ge 36$. When $k\ge 60$, the one-level method fails for the majority of $N$ being tested. 
In contrast, both two-level variants are robust at low to moderate $k$, with $\Delta_k$-GenEO consistently achieving the 
same or lower iteration count than the $\Delta$-GenEO for a smaller coarse space (CS).

\paragraph{Frequency sensitivity and thresholding}
For $k \le 40$, both two-level methods can yield nearly $N$ independent iteration counts, provided that an adequately sized coarse space is used. 
The $\Delta_k$-GenEO method is shown to require $\approx\,25$--$35\%$ fewer basis vectors in the coarse space when compared to the $\Delta$-GenEO method.
As $k$ increases, acceptable performance requires increasing the eigenvalue threshold $\tau$, subsequently increasing the $CS$. 
Using a $\tau\in[0.4,0.7]$ is sufficient for $k\lesssim 40$, while $\tau\ge 0.6$--$0.7$ is required for $k\in\{80,100\}$ 
to bring iterations down from the iteration limit to $40$--$80$ iterations when using the $\Delta_k$-GenEO coarse space. 
However, when using the $\Delta$-GenEO coarse space, increasing $\tau$ fails to bring the iteration count 
below the upper limit for the high $k$ values. 
Beyond these ranges, further increases in $\tau$ let the $CS$ grow substantially but can start to yield diminishing iteration gains. 
This effect is more prominent when using the $\Delta$-GenEO than for $\Delta_k$-GenEO.
{
Beyond both methods' inherent limitations for higher frequencies, numerical experiments provide insight into how the spectral threshold $\tau$ may be chosen in practice. 
Across all tested configurations, robust convergence of the $\Delta_k$-GenEO method is obtained for values of $\tau$ less than $1$ well below the worst-case bounds predicted by the theory, with only a modest increase in coarse-space dimension as $\tau$ is increased.
}
{\paragraph{Scaling with \texorpdfstring{$N$}{N}, \texorpdfstring{$H$}{H}, and \texorpdfstring{$k$}{k}}
With fixed spectral threshold $\tau$, the global coarse-space dimension grows approximately linearly 
with the number of subdomains $N$ and increases with the wavenumber $k$. 
At the subdomain level, the number of local coarse modes, denoted by $CS_{\mathrm{loc}}$, increases with the subdomain diameter 
measured in wavelengths, $H$, since fewer and wider subdomains support a larger number of low-energy local modes.}

{Consequently, for fixed $k$ and $\tau$, increasing $H$ reduces the number of subdomains and therefore decreases 
the total coarse-space dimension $CS$, while simultaneously increasing the average number of local contributions $CS_{\mathrm{loc}}$. 
For low to moderate frequencies, the GMRES iteration counts are fairly insensitive to changes in $H$ 
for both the $\Delta$-GenEO and $\Delta_k$-GenEO methods, indicating that the coarse spaces remain sufficiently rich to capture the relevant global modes. 
The total coarse-space size also differs only mildly between the two methods in this regime. 
At higher frequencies, however, the $\Delta_k$-GenEO method is more robust, reflecting its ability to retain 
additional modes required to stabilise the iteration in more indefinite regimes. 
As a general rule, achieving scalability with respect to the number of domains and frequency happens at the expense of a larger coarse space.}

\paragraph{Effect of heterogeneity (contrast \texorpdfstring{$a_{\max}$}{a\_max})}
Introducing heterogeneity produces mild changes in $CS$ size and, in many cases, slightly improves iterations for the two-level methods. For example, at $k=60$ and fixed $N$, increasing $a_{\max}$ often decreases the iteration counts for both coarse spaces, with $\Delta_k$-GenEO maintaining a clear advantage in both iterations and $CS$ size.
Overall, the behaviours observed in the homogeneous case persist under heterogeneity, with only weak dependence on $a_{\max}$. 
The one-level method is adversely effected by the heterogeneity, obtaining higher iteration counts even with low $k$ and $N$, where the one-level method performed more favourably in the homogenous case. 

\paragraph{High-frequency regime}
At $k\ge 60$, $\Delta_k$-GenEO remains viable across the tested $N$ provided $\tau$ is large enough (typically $\ge 0.6$--$0.7$). In comparison, $\Delta$-GenEO often reaches the iteration limit even as $\tau$ increases.

\paragraph{One-level baseline}
The one-level method should be avoided except for very small $k$ and $N$, and with no heterogeneity. 
Its iteration counts grow with $N$ even at $k=20$ and generally hit the iteration limit for $k\ge 40$.

From these results, we can see that for both the homogeneous and heterogeneous tests, the $\Delta_k$-GenEO is significantly 
more robust and scalable than the one-level method.  It can remain effective at high $k$ provided the eigenvalue threshold $\tau$ 
is increased sufficiently. Whilst the $\Delta$-GenEO is effective for moderate $k$, it requires a larger coarse space 
in order to see the same reduction in iteration counts, and still fails to converge within the iteration limit when $k$ is high. 
The one-level method on the other hand, is not competitive except at the very lowest frequencies. 
It has also been noticed that with an increase in the heterogeneity, it was possible to obtain fewer iteration counts. 
The lower iteration counts for some of the heterogeneous problems may stem from the fact that in some parts of the 
domain the heterogeneity in fact reduces the strength of the indefiniteness. 
For example, if we let $\Omegamax$ be a section of the domain where $a(\bx)$ is at its maximum value, $\amax$, 
then if we restrict problem \eqref{eq:problem} to $\Omegamax$, we see that \eqref{eq:problem} can be reduced to
\begin{equation*}
    \Delta u - \amax^{-1}k^2u = \amax^{-1} f.
\end{equation*}
Thus the restriction acts to reduce the effective indefiniteness in this region.
If using the $\Delta_k$-GenEO, then the following  operational guidance should be followed,
\begin{itemize}
  \item For low to moderate frequency ($k\lesssim 40$): choose $\tau\approx 0.4\!-\!0.7$ to get $\sim 20\!-\!50$ iterations with moderate $CS$.
  \item For high frequency (e.g.\ $k\ge 80$): set $\tau\approx 0.6\!-\!0.7$ to keep iterations $\lesssim 50\!-\!80$ across $N$.
  \item Expect $CS\propto N$ and $CS_{\text{loc}}$ to grow with $H$ and $k$. 
\end{itemize}
}

\subsection{Comparison with $H_k$-GenEO}

{To further clarify the positioning of the proposed $\Delta_k$-GenEO method, we include a brief numerical comparison with the $H_k$-GenEO coarse space briefly mentioned in the introduction, where the local spectral problem involves the full (indefinite) operator. That is, for each $i\in\{1,\ldots,N\}$ we solve the local eigenvalue problem consisting in finding the eigenpairs $(p,\lambda)$ such that
\[
	b_{\Omega_i}(p,v) = \lambda\bigl(\Xi_i p,\Xi_i v\bigr)_{1,k,\Omega_i}
	\qquad\text{for all} \ v\in\widetilde{V}_i
\]
 The discrete versions of eigenfunctions $p$ are added then to the coarse space in the same fashion as for the $\Delta_k$-GenEO method.}

{The problem being studied for the comparison between the $\Delta_k$-GenEO and $H_k$-GenEO coarse space is that of the homogenous problem as outlined in Section \ref{sec:HomProb}. Table \ref{table:HkDeltaHomog} displays the results comparing the one-level AS and the two-level AS with the $H_k$-GenEO and the $\Delta_k$-GenEO. The indefiniteness of the problem is controlled by $k$, and the number of subdomains is given by $N$. All other entries are as previously defined. Figure~\ref{fig:HkDeltaHomCSIT} shows the number of iterations to achieve the required GMRES tolerance as a function of the total coarse space size, and the selected threshold. Figure~\ref{fig:HkDeltaHomH} shows the number of GMRES iterations and coarse space size as a function of the subdomain diameter in wavelengths.
}

{
We notice that for low to moderate frequencies, both approaches exhibit robust two-level convergence, with iteration counts that depend only weakly on the number of subdomains once a moderate spectral threshold is used. In this regime, $\Delta_k$-GenEO remains comparable to $H_k$-GenEO, although it typically requires a slightly larger coarse space to achieve similar iteration counts. As the frequency increases, however, the behaviour of the two methods diverges more clearly: while $H_k$-GenEO maintains fast convergence for moderate thresholds, $\Delta_k$-GenEO becomes increasingly sensitive to the wavenumber and the subdomain partition, requiring larger coarse spaces and higher iteration counts to converge. }

{This highlights an intrinsic limitation of purely SPD-based coarse spaces in the high-frequency, strongly indefinite regime. Nevertheless, an important advantage of $\Delta_k$-GenEO lies in its conceptual and theoretical simplicity: the associated local spectral problems are symmetric positive definite and admit a relatively complete and transparent analysis. In contrast, the fully indefinite spectral problems underlying $H_k$-GenEO are considerably more involved, and a comprehensive theoretical understanding of their robustness and scalability remains, to date, far from established.
}

\begin{table}[H]
    \centering
    \begin{tabular}{ccccccc|c|ccc|ccc}
\multicolumn{7}{c|}{} & One Level & \multicolumn{3}{c}{$H_k$-GenEO} & \multicolumn{3}{c}{$\Delta_k$-GenEO} \\
$N$ & $k$ & $h^{-1}$ & $n$ & $L$ & $H$ & $n_\text{loc}$ & It & It & CS & $CS_\text{loc}$ & It & CS & $CS_\text{loc}$  \\
\hline
16 & 20 & 240 & 58081 & 3.18 & 0.82 & 3969 & 41 & 16 & 240 & 18 & 21 & 340 & 21 \\
 & 40 & 480 & 231361 &  6.37 & 1.62 & 15129 & 107 & 20 & 628 & 49 & 29 & 784 & 49 \\
 & 60 & 720 & 519841 &  9.55 & 2.41 & 33489 & 153 & 17 & 1160 & 86 & 31 & 1312 & 82 \\
 & 80 & 960 & 923521 &  12.73 & 3.21 & 59049 & $>$200 & 27 & 1876 & 132 & 75 & 1936 & 121 \\
 & 100 & 1200 & 1442401 &  15.92 & 4.01 & 91809 & $>$200 & 36 & 2660 & 190 & 68 & 2672 & 167 \\
\hline
36 & 20 & 240 & 58081 & 3.18 & 0.56 & 1849 & 79 & 19 & 344 & 12 & 22 & 516 & 14 \\
 & 40 & 480 & 231361 &  6.37 & 1.09 & 6889 & $>$200 & 31 & 844 & 28 & 56 & 1168 & 32 \\
 & 60 & 720 & 519841 &  9.55 & 1.62 & 15129 & $>$200 & 18 & 1528 & 49 & 34 & 1932 & 54 \\
 & 80 & 960 & 923521 &  12.73 & 2.15 & 26569 & $>$200 & 23 & 2312 & 72 & 68 & 2740 & 76 \\
 & 100 & 1200 & 1442401 &  15.92 & 2.68 & 41209 & $>$200 & 51 & 3260 & 99 & 90 & 3684 & 102 \\
\hline
64 & 20 & 240 & 58081 & 3.18 & 0.42 & 1089 & 92 & 19 & 448 & 8 & 23 & 736 & 12 \\
 & 40 & 480 & 231361 &  6.37 & 0.82 & 3969 & 175 & 24 & 1056 & 18 & 56 & 1604 & 25 \\
 & 60 & 720 & 519841 &  9.55 & 1.22 & 8649 & $>$200 & 18 & 1924 & 34 & 72 & 2460 & 38 \\
 & 80 & 960 & 923521 &  12.73 & 1.62 & 15129 & $>$200 & 23 & 2820 & 49 & 76 & 3584 & 56 \\
 & 100 & 1200 & 1442401 &  15.92 & 2.02 & 23409 & $>$200 & 21 & 3840 & 65 & 88 & 4616 & 72 \\
\hline
100 & 20 & 240 & 58081 & 3.18 & 0.34 & 729 & 96 & 18 & 684 & 8 & 21 & 980 & 10 \\
 & 40 & 480 & 231361 &  6.37 & 0.66 & 2601 & $>$200 & 24 & 1440 & 16 & 87 & 1936 & 19 \\
 & 60 & 720 & 519841 &  9.55 & 0.98 & 5625 & $>$200 & 18 & 2200 & 24 & 104 & 3020 & 30 \\
 & 80 & 960 & 923521 &  12.73 & 1.30 & 9801 & $>$200 & 24 & 3224 & 34 & 144 & 4324 & 43 \\
 & 100 & 1200 & 1442401 &  15.92 & 1.62 & 15129 & $>$200 & 24 & 4504 & 49 & 80 & 5740 & 57 \\
\hline
144 & 20 & 240 & 58081 & 3.18 & 0.29 & 529 & 104 & 21 & 672 & 5 & 19 & 1056 & 7 \\
 & 40 & 480 & 231361 &  6.37 & 0.56 & 1849 & $>$200 & 30 & 1544 & 12 & 64 & 2256 & 16 \\
 & 60 & 720 & 519841 &  9.55 & 0.82 & 3969 & $>$200 & 23 & 2448 & 18 & 77 & 3796 & 26 \\
 & 80 & 960 & 923521 &  12.73 & 1.09 & 6889 & $>$200 & 45 & 3700 & 28 & 162 & 5140 & 36 \\
 & 100 & 1200 & 1442401 &  15.92 & 1.35 & 10609 & $>$200 & 25 & 5232 & 39 & 197 & 6576 & 46 \\
\hline
\end{tabular}
    \caption{Results showing the dimension of the fine mesh ($n$), the diameter of the domain measured in wavelengths ($L$), the diameter of the sub-domains measured in wavelengths ($H$), the number of local degrees of freedom ($n_\text{loc}$), the GMRES iteration count (It.) for the one-level and two level methods using the $H_k$-GenEO and $\Delta_k$-GenEO coarse space, the dimension of the coarse space (CS) and averaged number of contributions to the coarse space per subdomain (CS$_\text{loc}$) for the $H_k$-GenEO and $\Delta_k$-GenEO. These results are for the homogeneous media test case. Using $\tau = 0.4$ and $\tau=0.6$ for $H_k$-GenEO and $\Delta_k$-GenEO respectively.}
    \label{table:HkDeltaHomog}
\end{table}

\begin{figure}[H]
    \centering
    \begin{subfigure}[t]{0.48\textwidth}
        \centering
        \begin{tikzpicture}
\begin{axis}[
  xmode=log,
  ymode=log,
  enlarge x limits=0.05,
  enlarge y limits=0.05,
  xmin=36,
  xmax=1540,
  ymin=13,
  ymax=104,
  grid=both,
  minor grid style={opacity=0.3},
  major grid style={opacity=0.3},
  xlabel={CS Size},
  ylabel={Iteration Count},
  width=\textwidth,
  height=\textwidth,
  legend pos=north east,
  legend style={fill=none, draw=none, font=\scriptsize, text opacity=1}
]
\addplot+[mark=square*, mark options={solid, fill = blue}, mark size=1.8pt, line width=1pt, color=blue, opacity=1, style=solid] coordinates {(68,71) (128,29) (172,19) (240,16) (320,14) (532,13)};
\addplot+[mark=square*, mark options={solid, fill = orange}, mark size=1.8pt, line width=1pt, color=orange, opacity=1, style=solid] coordinates {(84,59) (220,22) (260,20) (344,19) (496,15) (808,14)};
\addplot+[mark=square*, mark options={solid, fill = green!60!black}, mark size=1.8pt, line width=1pt, color=green!60!black, opacity=1, style=solid] coordinates {(160,57) (224,32) (420,19) (448,19) (644,17) (1024,14)};
\addplot+[mark=square*, mark options={solid, fill = red!80!black}, mark size=1.8pt, line width=1pt, color=red!80!black, opacity=1, style=solid] coordinates {(196,91) (324,31) (396,29) (684,18) (720,18) (1180,15)};
\addplot+[mark=square*, mark options={solid, fill = cyan}, mark size=1.8pt, line width=1pt, color=cyan, opacity=1, style=solid] coordinates {(344,50) (484,27) (484,27) (672,21) (1012,17) (1540,16)};
\addplot+[mark=*, mark options={solid, fill = blue}, mark size=1.8pt, line width=1pt, color=blue, opacity=1, style=dashed] coordinates {(36,48) (84,51) (132,39) (180,35) (272,22) (340,21) (464,16)};
\addplot+[mark=*, mark options={solid, fill = orange}, mark size=1.8pt, line width=1pt, color=orange, opacity=1, style=dashed] coordinates {(84,62) (140,38) (220,27) (340,24) (428,23) (516,22) (720,16)};
\addplot+[mark=*, mark options={solid, fill = green!60!black}, mark size=1.8pt, line width=1pt, color=green!60!black, opacity=1, style=dashed] coordinates {(136,73) (196,38) (288,33) (420,25) (608,23) (736,23) (920,22)};
\addplot+[mark=*, mark options={solid, fill = red!80!black}, mark size=1.8pt, line width=1pt, color=red!80!black, opacity=1, style=dashed] coordinates {(100,100) (324,34) (360,33) (620,26) (684,22) (980,21) (1180,21)};
\addplot+[mark=*, mark options={solid, fill = cyan}, mark size=1.8pt, line width=1pt, color=cyan, opacity=1, style=dashed] coordinates {(144,99) (384,38) (484,29) (672,26) (1012,19) (1056,19) (1540,18)};
\addplot+[mark=none, black, dashed, line width=1.2pt, opacity=0.2] coordinates {(36,41) (1540,41)};
\addplot+[mark=none, black, dashed, line width=1.2pt, opacity=0.4] coordinates {(36,79) (1540,79)};
\addplot+[mark=none, black, dashed, line width=1.2pt, opacity=0.6] coordinates {(36,92) (1540,92)};
\addplot+[mark=none, black, dashed, line width=1.2pt, opacity=0.8] coordinates {(36,96) (1540,96)};
\addplot+[mark=none, black, dashed, line width=1.2pt, opacity=1.0] coordinates {(36,104) (1540,104)};
\end{axis}
\end{tikzpicture}
    \end{subfigure}   
    \hfill
    \begin{subfigure}[t]{0.48\textwidth}
        \centering
        \begin{tikzpicture}
\begin{axis}[
  xmode=log,
  ymode=log,
  enlarge x limits=0.05,
  enlarge y limits=0.05,
  xmin=0.1,
  xmax=0.7,
  ymin=13,
  ymax=104,
  grid=both,
  minor grid style={opacity=0.3},
  major grid style={opacity=0.3},
  xlabel={$\tau$},
  ylabel={Iteration Count},
  width=\textwidth,
  height=\textwidth,
  legend pos=north east,
  legend style={fill=none, draw=none, font=\scriptsize, text opacity=1}
]
\addplot+[mark=square*, mark options={solid, fill = blue}, mark size=1.8pt, line width=1pt, color=blue, opacity=1, style=solid] coordinates {(0.1,71) (0.2,29) (0.3,19) (0.4,16) (0.5,14) (0.7,13)};
\addplot+[mark=square*, mark options={solid, fill = orange}, mark size=1.8pt, line width=1pt, color=orange, opacity=1, style=solid] coordinates {(0.1,59) (0.2,22) (0.3,20) (0.4,19) (0.5,15) (0.7,14)};
\addplot+[mark=square*, mark options={solid, fill = green!60!black}, mark size=1.8pt, line width=1pt, color=green!60!black, opacity=1, style=solid] coordinates {(0.1,57) (0.2,32) (0.3,19) (0.4,19) (0.5,17) (0.7,14)};
\addplot+[mark=square*, mark options={solid, fill = red!80!black}, mark size=1.8pt, line width=1pt, color=red!80!black, opacity=1, style=solid] coordinates {(0.1,91) (0.2,31) (0.3,29) (0.4,18) (0.5,18) (0.7,15)};
\addplot+[mark=square*, mark options={solid, fill = cyan}, mark size=1.8pt, line width=1pt, color=cyan, opacity=1, style=solid] coordinates {(0.1,50) (0.2,27) (0.3,27) (0.4,21) (0.5,17) (0.7,16)};
\addplot+[mark=*, mark options={solid, fill = blue}, mark size=1.8pt, line width=1pt, color=blue, opacity=1, style=dashed] coordinates {(0.1,48) (0.2,51) (0.3,39) (0.4,35) (0.5,22) (0.6,21) (0.7,16)};
\addplot+[mark=*, mark options={solid, fill = orange}, mark size=1.8pt, line width=1pt, color=orange, opacity=1, style=dashed] coordinates {(0.1,62) (0.2,38) (0.3,27) (0.4,24) (0.5,23) (0.6,22) (0.7,16)};
\addplot+[mark=*, mark options={solid, fill = green!60!black}, mark size=1.8pt, line width=1pt, color=green!60!black, opacity=1, style=dashed] coordinates {(0.1,73) (0.2,38) (0.3,33) (0.4,25) (0.5,23) (0.6,23) (0.7,22)};
\addplot+[mark=*, mark options={solid, fill = red!80!black}, mark size=1.8pt, line width=1pt, color=red!80!black, opacity=1, style=dashed] coordinates {(0.1,100) (0.2,34) (0.3,33) (0.4,26) (0.5,22) (0.6,21) (0.7,21)};
\addplot+[mark=*, mark options={solid, fill = cyan}, mark size=1.8pt, line width=1pt, color=cyan, opacity=1, style=dashed] coordinates {(0.1,99) (0.2,38) (0.3,29) (0.4,26) (0.5,19) (0.6,19) (0.7,18)};
\addplot+[mark=none, black, dashed, line width=1.2pt, opacity=0.2] coordinates {(0.1,41) (0.7,41)};
\addplot+[mark=none, black, dashed, line width=1.2pt, opacity=0.4] coordinates {(0.1,79) (0.7,79)};
\addplot+[mark=none, black, dashed, line width=1.2pt, opacity=0.6] coordinates {(0.1,92) (0.7,92)};
\addplot+[mark=none, black, dashed, line width=1.2pt, opacity=0.8] coordinates {(0.1,96) (0.7,96)};
\addplot+[mark=none, black, dashed, line width=1.2pt, opacity=1.0] coordinates {(0.1,104) (0.7,104)};
\end{axis}
\end{tikzpicture}
    \end{subfigure} 
    \begin{subfigure}[t]{0.48\textwidth}
        \centering
        \begin{tikzpicture}
\begin{axis}[
  xmode=log,
  ymode=log,
  enlarge x limits=0.05,
  enlarge y limits=0.05,
  xmin=440,
  xmax=8164,
  ymin=16,
  ymax=200,
  grid=both,
  minor grid style={opacity=0.3},
  major grid style={opacity=0.3},
  xlabel={CS Size},
  ylabel={Iteration Count},
  width=\textwidth,
  height=\textwidth,
  legend pos=north east,
  legend style={fill=none, draw=none, font=\small, text opacity=1},
  legend columns = 5,
  legend to name = domains
]
\addplot+[mark=square*, mark options={solid, fill = blue}, mark size=1.8pt, line width=1pt, color=blue, opacity=1] coordinates {(800,86) (1224,64) (1700,32) (2260,26) (4280,23)};
\addlegendentry{$H_k$-GenEO (16)}
\addplot+[mark=square*, mark options={solid, fill = orange}, mark size=1.8pt, line width=1pt, color=orange, opacity=1] coordinates {(1000,143) (1584,33) (2168,24) (2900,20) (5224,16)};
\addlegendentry{$H_k$-GenEO (36)}
\addplot+[mark=square*, mark options={solid, fill = green!60!black}, mark size=1.8pt, line width=1pt, color=green!60!black, opacity=1] coordinates {(1284,153) (1892,32) (2660,23) (3552,20) (6372,17)};
\addlegendentry{$H_k$-GenEO (64)}
\addplot+[mark=square*, mark options={solid, fill = red!80!black}, mark size=1.8pt, line width=1pt, color=red!80!black, opacity=1] coordinates {(1444,82) (2400,36) (3124,24) (4144,21) (7212,17)};
\addlegendentry{$H_k$-GenEO (100)}
\addplot+[mark=square*, mark options={solid, fill = cyan}, mark size=1.8pt, line width=1pt, color=cyan, opacity=1] coordinates {(1776,175) (2884,55) (3556,48) (4944,36) (8164,26)};
\addlegendentry{$H_k$-GenEO (144)}
\addplot+[mark=*, mark options={solid, fill = blue}, mark size=1.8pt, line width=1pt, color=blue, opacity=1, style=dashed] coordinates {(440,200) (712,192) (1060,189) (1424,155) (1936,75) (2692,38)};
\addlegendentry{$\Delta_k$-GenEO (16)}
\addplot+[mark=*, mark options={solid, fill = orange}, mark size=1.8pt, line width=1pt, color=orange, opacity=1, style=dashed] coordinates {(672,200) (1052,200) (1516,200) (2028,200) (2740,68) (3696,42)};
\addlegendentry{$\Delta_k$-GenEO (36)}
\addplot+[mark=*, mark options={solid, fill = green!60!black}, mark size=1.8pt, line width=1pt, color=green!60!black, opacity=1, style=dashed] coordinates {(872,200) (1380,200) (2012,200) (2692,200) (3584,76) (4704,51)};
\addlegendentry{$\Delta_k$-GenEO (64)}
\addplot+[mark=*, mark options={solid, fill = red!80!black}, mark size=1.8pt, line width=1pt, color=red!80!black, opacity=1, style=dashed] coordinates {(1112,200) (1864,200) (2520,200) (3440,192) (4324,144) (5804,55)};
\addlegendentry{$\Delta_k$-GenEO (100)}
\addplot+[mark=*, mark options={solid, fill = cyan}, mark size=1.8pt, line width=1pt, color=cyan, opacity=1, style=dashed] coordinates {(1540,200) (2068,200) (2928,200) (3800,200) (5140,162) (6736,101)};
\addlegendentry{$\Delta_k$-GenEO (144)}
\addplot+[mark=none, black, dashed, line width=1.2pt, opacity=0.2] coordinates {(440,200) (8164,200)};
\addlegendentry{One-level (16)}
\addplot+[mark=none, black, dashed, line width=1.2pt, opacity=0.4] coordinates {(440,200) (8164,200)};
\addlegendentry{One-level (36)}
\addplot+[mark=none, black, dashed, line width=1.2pt, opacity=0.6] coordinates {(440,200) (8164,200)};
\addlegendentry{One-level (64)}
\addplot+[mark=none, black, dashed, line width=1.2pt, opacity=0.8] coordinates {(440,200) (8164,200)};
\addlegendentry{One-level (100)}
\addplot+[mark=none, black, dashed, line width=1.2pt, opacity=1.0] coordinates {(440,200) (8164,200)};
\addlegendentry{One-level (144)}
\end{axis}
\end{tikzpicture}
    \end{subfigure} 
    \hfill
    \begin{subfigure}[t]{0.48\textwidth}
        \centering
        \begin{tikzpicture}
\begin{axis}[
  xmode=log,
  ymode=log,
  enlarge x limits=0.05,
  enlarge y limits=0.05,
  xmin=0.2,
  xmax=0.7,
  ymin=16,
  ymax=200,
  grid=both,
  minor grid style={opacity=0.3},
  major grid style={opacity=0.3},
  xlabel={$\tau$},
  ylabel={Iteration Count},
  width=\textwidth,
  height=\textwidth,
  legend pos=north east,
  legend style={fill=none, draw=none, font=\scriptsize, text opacity=1}
]
\addplot+[mark=*, mark options={solid, fill = blue}, mark size=1.8pt, line width=1pt, color=blue, opacity=1] coordinates {(0.2,86) (0.3,64) (0.4,32) (0.5,26) (0.6,25) (0.7,23)};
\addplot+[mark=*, mark options={solid, fill = orange}, mark size=1.8pt, line width=1pt, color=orange, opacity=1] coordinates {(0.2,143) (0.3,33) (0.4,24) (0.5,20) (0.6,18) (0.7,16)};
\addplot+[mark=*, mark options={solid, fill = green!60!black}, mark size=1.8pt, line width=1pt, color=green!60!black, opacity=1] coordinates {(0.2,153) (0.3,32) (0.4,23) (0.5,20) (0.6,18) (0.7,17)};
\addplot+[mark=*, mark options={solid, fill = red!80!black}, mark size=1.8pt, line width=1pt, color=red!80!black, opacity=1] coordinates {(0.2,82) (0.3,36) (0.4,24) (0.5,21) (0.6,19) (0.7,17)};
\addplot+[mark=*, mark options={solid, fill = cyan}, mark size=1.8pt, line width=1pt, color=cyan, opacity=1] coordinates {(0.2,175) (0.3,55) (0.4,48) (0.5,36) (0.6,33) (0.7,26)};
\addplot+[mark=*, mark options={solid, fill = blue}, mark size=1.8pt, line width=1pt, color=blue, opacity=1, style=dashed] coordinates {(0.2,200) (0.3,192) (0.4,189) (0.5,155) (0.6,75) (0.7,38)};
\addplot+[mark=*, mark options={solid, fill = orange}, mark size=1.8pt, line width=1pt, color=orange, opacity=1, style=dashed] coordinates {(0.2,200) (0.3,200) (0.4,200) (0.5,200) (0.6,68) (0.7,42)};
\addplot+[mark=*, mark options={solid, fill = green!60!black}, mark size=1.8pt, line width=1pt, color=green!60!black, opacity=1, style=dashed] coordinates {(0.2,200) (0.3,200) (0.4,200) (0.5,200) (0.6,76) (0.7,51)};
\addplot+[mark=*, mark options={solid, fill = red!80!black}, mark size=1.8pt, line width=1pt, color=red!80!black, opacity=1, style=dashed] coordinates {(0.2,200) (0.3,200) (0.4,200) (0.5,192) (0.6,144) (0.7,55)};
\addplot+[mark=*, mark options={solid, fill = cyan}, mark size=1.8pt, line width=1pt, color=cyan, opacity=1, style=dashed] coordinates {(0.2,200) (0.3,200) (0.4,200) (0.5,200) (0.6,162) (0.7,101)};
\addplot+[mark=none, black, dashed, line width=1.2pt, opacity=0.2] coordinates {(0.2,200) (0.7,200)};
\addplot+[mark=none, black, dashed, line width=1.2pt, opacity=0.4] coordinates {(0.2,200) (0.7,200)};
\addplot+[mark=none, black, dashed, line width=1.2pt, opacity=0.6] coordinates {(0.2,200) (0.7,200)};
\addplot+[mark=none, black, dashed, line width=1.2pt, opacity=0.8] coordinates {(0.2,200) (0.7,200)};
\addplot+[mark=none, black, dashed, line width=1.2pt, opacity=1.0] coordinates {(0.2,200) (0.7,200)};
\end{axis}
\end{tikzpicture}
    \end{subfigure} 
    \ref{domains}
    \caption{Influence of the coarse space size (left) and threshold choice (right) on the iteration count for the homogeneous media test case with $k$ = 20 (top) and $k$ = 100 (bottom). The number in brackets indicates the number of subdomains.}
    \label{fig:HkDeltaHomCSIT}
\end{figure}

\begin{figure}[H]
    \centering
    \begin{subfigure}[t]{0.48\textwidth}
        \centering
        \begin{tikzpicture}
\begin{axis}[
  xmode=log,
  ymode=log,
  enlarge x limits=0.05,
  enlarge y limits=0.05,
  xmin=0.29178406,
  xmax=0.82230054,
  ymin=84,
  ymax=1540,
  grid=both,
  minor grid style={opacity=0.3},
  major grid style={opacity=0.3},
  xlabel={$H$ (subdomain side length in wavelengths)},
  ylabel={CS Size},
  width=\textwidth,
  height=\textwidth,
  legend pos=north east,
  legend style={fill=none, draw=none, font=\scriptsize, text opacity=1}
]
\addplot+[mark=square*, mark options={solid, fill = blue}, mark size=1.8pt, line width=1pt, color=blue, opacity=1, style=solid] coordinates {(0.29178406,484) (0.34483571,324) (0.42441318,224) (0.5570423,220) (0.82230054,128)};
\addplot+[mark=square*, mark options={solid, fill = orange}, mark size=1.8pt, line width=1pt, color=orange, opacity=1, style=solid] coordinates {(0.29178406,672) (0.34483571,684) (0.42441318,448) (0.5570423,344) (0.82230054,240)};
\addplot+[mark=square*, mark options={solid, fill = cyan}, mark size=1.8pt, line width=1pt, color=cyan, opacity=1, style=solid] coordinates {(0.29178406,1540) (0.34483571,1180) (0.42441318,1024) (0.5570423,808) (0.82230054,532)};
\addplot+[mark=*, mark options={solid, fill = blue}, mark size=1.8pt, line width=1pt, color=blue, opacity=1, style=dashed] coordinates {(0.82230054,84) (0.5570423,140) (0.42441318,196) (0.34483571,324) (0.29178406,384)};
\addplot+[mark=*, mark options={solid, fill = orange}, mark size=1.8pt, line width=1pt, color=orange, opacity=1, style=dashed] coordinates {(0.82230054,180) (0.5570423,340) (0.42441318,420) (0.34483571,620) (0.29178406,672)};
\addplot+[mark=*, mark options={solid, fill = cyan}, mark size=1.8pt, line width=1pt, color=cyan, opacity=1, style=dashed] coordinates {(0.82230054,464) (0.5570423,720) (0.42441318,920) (0.34483571,1180) (0.29178406,1540)};
\end{axis}
\end{tikzpicture}
    \end{subfigure}   
    \hfill
    \begin{subfigure}[t]{0.48\textwidth}
        \centering
        \begin{tikzpicture}
\begin{axis}[
  xmode=log,
  ymode=log,
  enlarge x limits=0.05,
  enlarge y limits=0.05,
  xmin=0.29178406,
  xmax=0.82230054,
  ymin=13,
  ymax=104,
  grid=both,
  minor grid style={opacity=0.3},
  major grid style={opacity=0.3},
  xlabel={$H$ (subdomain side length in wavelengths)},
  ylabel={Iteration Count},
  width=\textwidth,
  height=\textwidth,
  legend pos=north east,
  legend style={fill=none, draw=none, font=\scriptsize, text opacity=1}
]
\addplot+[mark=square*, mark options={solid, fill = blue}, mark size=1.8pt, line width=1pt, color=blue, opacity=1, style=solid] coordinates {(0.29178406,27) (0.34483571,31) (0.42441318,32) (0.5570423,22) (0.82230054,29)};
\addplot+[mark=square*, mark options={solid, fill = orange}, mark size=1.8pt, line width=1pt, color=orange, opacity=1, style=solid] coordinates {(0.29178406,21) (0.34483571,18) (0.42441318,19) (0.5570423,19) (0.82230054,16)};
\addplot+[mark=square*, mark options={solid, fill = cyan}, mark size=1.8pt, line width=1pt, color=cyan, opacity=1, style=solid] coordinates {(0.29178406,16) (0.34483571,15) (0.42441318,14) (0.5570423,14) (0.82230054,13)};
\addplot+[mark=*, mark options={solid, fill = blue}, mark size=1.8pt, line width=1pt, color=blue, opacity=1, style=dashed] coordinates {(0.29178406,38) (0.34483571,34) (0.42441318,38) (0.5570423,38) (0.82230054,51)};
\addplot+[mark=*, mark options={solid, fill = orange}, mark size=1.8pt, line width=1pt, color=orange, opacity=1, style=dashed] coordinates {(0.29178406,26) (0.34483571,26) (0.42441318,25) (0.5570423,24) (0.82230054,35)};
\addplot+[mark=*, mark options={solid, fill = cyan}, mark size=1.8pt, line width=1pt, color=cyan, opacity=1, style=dashed] coordinates {(0.29178406,18) (0.34483571,21) (0.42441318,22) (0.5570423,16) (0.82230054,16)};
\addplot+[mark=none, black, dashed, line width=1.2pt, opacity=1] coordinates {(0.29178406,104) (0.34483571,96) (0.42441318,92) (0.5570423,79) (0.82230054,41)};
\end{axis}
\end{tikzpicture}
    \end{subfigure} 
    \begin{subfigure}[t]{0.48\textwidth}
        \centering
        \begin{tikzpicture}
\begin{axis}[
  xmode=log,
  ymode=log,
  enlarge x limits=0.05,
  enlarge y limits=0.05,
  xmin=1.0875588,
  xmax=3.2096247,
  ymin=440,
  ymax=8164,
  grid=both,
  minor grid style={opacity=0.3},
  major grid style={opacity=0.3},
  xlabel={$H$ (subdomain side length in wavelengths)},
  ylabel={CS Size},
  width=\textwidth,
  height=\textwidth,
  legend pos=north east,
  legend style={fill=none, draw=none, font=\small, text opacity=1}
]
\addplot+[mark=*, mark options={solid, fill = blue}, mark size=1.8pt, line width=1pt, color=blue, opacity=1] coordinates {(1.0875588,1776) (1.2997654,1444) (1.6180753,1284) (2.1485917,1000) (3.2096247,800)};
\addplot+[mark=*, mark options={solid, fill = orange}, mark size=1.8pt, line width=1pt, color=orange, opacity=1] coordinates {(1.0875588,3556) (1.2997654,3124) (1.6180753,2660) (2.1485917,2168) (3.2096247,1700)};
\addplot+[mark=*, mark options={solid, fill = cyan}, mark size=1.8pt, line width=1pt, color=cyan, opacity=1] coordinates {(1.0875588,8164) (1.2997654,7212) (1.6180753,6372) (2.1485917,5224) (3.2096247,4280)};
\addplot+[mark=*, mark options={solid, fill = blue}, mark size=1.8pt, line width=1pt, color=blue, opacity=1, style=dashed] coordinates {(3.2096247,440) (2.1485917,672) (1.6180753,872) (1.2997654,1112) (1.0875588,1540)};
\addplot+[mark=*, mark options={solid, fill = orange}, mark size=1.8pt, line width=1pt, color=orange, opacity=1, style=dashed] coordinates {(3.2096247,1060) (2.1485917,1516) (1.6180753,2012) (1.2997654,2520) (1.0875588,2928)};
\addplot+[mark=*, mark options={solid, fill = cyan}, mark size=1.8pt, line width=1pt, color=cyan, opacity=1, style=dashed] coordinates {(3.2096247,2692) (2.1485917,3696) (1.6180753,4704) (1.2997654,5804) (1.0875588,6736)};
\end{axis}
\end{tikzpicture}
    \end{subfigure} 
    \hfill
    \begin{subfigure}[t]{0.48\textwidth}
        \centering
        \begin{tikzpicture}
\begin{axis}[
  xmode=log,
  ymode=log,
  enlarge x limits=0.05,
  enlarge y limits=0.05,
  xmin=1.0875588,
  xmax=3.2096247,
  ymin=16,
  ymax=200,
  grid=both,
  minor grid style={opacity=0.3},
  major grid style={opacity=0.3},
  xlabel={$H$ (subdomain side length in wavelengths)},
  ylabel={Iteration Count},
  width=\textwidth,
  height=\textwidth,
  legend pos=north east,
  legend style={fill=none, draw=none, font=\small, text opacity=1},
  legend columns = 4,
  legend to name = Hwaves
]
\addplot+[mark=*, mark options={solid, fill = blue}, mark size=1.8pt, line width=1pt, color=blue, opacity=1] coordinates {(1.0875588,175) (1.2997654,82) (1.6180753,153) (2.1485917,143) (3.2096247,86)};
\addlegendentry{$H_k$-GenEO (0.2)}
\addplot+[mark=*, mark options={solid, fill = orange}, mark size=1.8pt, line width=1pt, color=orange, opacity=1] coordinates {(1.0875588,48) (1.2997654,24) (1.6180753,23) (2.1485917,24) (3.2096247,32)};
\addlegendentry{$H_k$-GenEO (0.4)}
\addplot+[mark=*, mark options={solid, fill = cyan}, mark size=1.8pt, line width=1pt, color=cyan, opacity=1] coordinates {(1.0875588,26) (1.2997654,17) (1.6180753,17) (2.1485917,16) (3.2096247,23)};
\addlegendentry{$H_k$-GenEO (0.7)}
\addplot+[mark=*, mark options={solid, fill = blue}, mark size=1.8pt, line width=1pt, color=blue, opacity=1, style=dashed] coordinates {(1.0875588,200) (1.2997654,200) (1.6180753,200) (2.1485917,200) (3.2096247,200)};
\addlegendentry{$\Delta_k$-GenEO (0.2)}
\addplot+[mark=*, mark options={solid, fill = orange}, mark size=1.8pt, line width=1pt, color=orange, opacity=1, style=dashed] coordinates {(1.0875588,200) (1.2997654,200) (1.6180753,200) (2.1485917,200) (3.2096247,189)};
\addlegendentry{$\Delta_k$-GenEO (0.4)}
\addplot+[mark=*, mark options={solid, fill = cyan}, mark size=1.8pt, line width=1pt, color=cyan, opacity=1, style=dashed] coordinates {(1.0875588,101) (1.2997654,55) (1.6180753,51) (2.1485917,42) (3.2096247,38)};
\addlegendentry{$\Delta_k$-GenEO (0.7)}
\addplot+[mark=none, black, dashed, line width=1.2pt, opacity=1] coordinates {(1.0875588,200) (1.2997654,200) (1.6180753,200) (2.1485917,200) (3.2096247,200)};
\addlegendentry{One-level}
\end{axis}
\end{tikzpicture}
    \end{subfigure} 
    \ref{Hwaves}
    \caption{Influence of the subdomain diameter on the iteration count (left) and coarse space size (right) for the homogeneous media test case $k$ = 20 (top) and $k$ = 100 (bottom). The number in brackets indicates $\tau$ used.}
    \label{fig:HkDeltaHomH}
\end{figure}

\section{Conclusion and discussion}

In this work, we have introduced the $\Delta_k$-GenEO coarse space. We then derived and rigorously analysed the required conditions 
to achieve robust convergence of GMRES when applied to the Helmholtz problem preconditioned by the $\Delta_k$-GenEO coarse space. 
Our main theoretical results deliver a refinement of the bounds on the subdomain diameter, $H$, and the eigenvalue tolerance, $\tau$, 
when compared to the $\Delta$-GenEO method.  Reducing the wavenumber dependence from $H \lesssim k^{-2}$ and $\tau \gtrsim k^{8}$ 
to $H \lesssim k^{-1}$ and $\tau \gtrsim k^{2}$, respectively. 
These improvements stem from strengthened estimates on the projection and stability operators used 
in the spectral construction of the coarse space.  The theoretical improvements are also reflected in the numerical results. 
It has been demonstrated that the $\Delta_k$-GenEO method outperforms the $\Delta$-GenEO method across all tested scenarios. 
Not only has the $\Delta_k$-GenEO shown to have improved robustness with respect to more indefinite problems, 
it is also able to achieve convergence in fewer iterations and using a smaller coarse space over that of the $\Delta$-GenEO method. 
While the improved conditions of the $\Delta_k$-GenEO are still conservative, they represent a meaningful step towards 
more practical and scalable preconditioners for indefinite wave problems.     

{Beyond the specific bounds obtained, this work provides insight into why GenEO-type coarse spaces based on 
symmetric positive definite operators can remain effective for indefinite Helmholtz problems. 
By introducing a $k$-dependent modification of the GenEO eigenvalue problem, we clarify how much of the observed robustness 
can be explained within an SPD-based framework, and where its fundamental limitations lie. 
In particular, the analysis helps to narrow the long-standing gap between conservative theoretical guarantees 
and the much milder coarse-space growth and iteration counts observed in practice.}

{Several directions for future research naturally follow from these results. 
These include the development and rigorous analysis of coarse spaces built directly from indefinite operators, 
multilevel extensions of the $\Delta_k$-GenEO framework, and optimised implementations suitable for large-scale parallel simulations. 
The present work provides a theoretical and numerical baseline against which such developments can be assessed.
}

\bibliographystyle{elsarticle-num}
\bibliography{paper_delta}

\end{document}